\documentclass[reqno,a4paper,11pt]{amsart}
\usepackage{color}
\usepackage{amsmath,amssymb,amsfonts,amsthm,enumerate}
\usepackage{mathrsfs,accents}
\usepackage[T1]{fontenc}
\usepackage[active]{srcltx}
\usepackage{epsfig}
\usepackage{graphicx}
\usepackage{cite}
\usepackage{hyperref} 
\usepackage{bookmark}
\usepackage[noabbrev,capitalize]{cleveref}

\usepackage[initials,
]{amsrefs}
\renewcommand{\MR}[1]{} 

\catcode`@=11 \@addtoreset{equation}{section} \catcode`@=12

\newtheorem{theorem}{Theorem}[section]
\newtheorem{lemma}[theorem]{Lemma}
\newtheorem{definition}[theorem]{Definition}
\newtheorem{remark}[theorem]{Remark}
\newtheorem{proposition}[theorem]{Proposition}
\newtheorem{corollary}[theorem]{Corollary}
\newtheorem{example}[theorem]{Example}

\newcommand{\N}{\mathbb{N}}
\newcommand{\R}{\mathbb{R}}

\renewcommand{\prec}[2][*]{#2^{#1}}
\def\bbbr{\mathbb{R}}
\def\R{\mathbb{R}}

\def\N{\mathbb{N}}

\def\pair#1{\left(#1\right)}

\newcommand{\eps}{\varepsilon}

\newcommand{\weakto}{\rightharpoonup}

\newcommand{\Haus}[1]{{\mathscr H}^{#1}} 
\newcommand{\Leb}[1]{{\mathscr L}^{#1}} 

\newcommand{\DM}{\mathcal{DM}}

\newcommand{\ban}[1]{\left\langle  #1 \right\rangle}  

\def\radon{\mathcal{M}(\Omega)}
\newcommand{\Borel}{\mathscr{B}(\Omega)}

\def\DM{{\mathcal{DM}}}

\def\cinftio{C^{\infty}_c(\Omega)}

\DeclareMathOperator{\Div}{div}
\def\pscal#1#2{\left\langle #1,\, #2 \right\rangle}
\DeclareMathOperator{\supp}{supp}

\def\upiu{u^+}
\def\umeno{u^-}
\def\uint{{u^{i}}}
\def\uext{{u^{e}}}

\newcommand{\quotemarks}[1]{``{#1}''}

\newcommand*{\chiut}[1][]{\chi^{#1}_{\{u > t\}}}

\usepackage[dvipsnames]{xcolor}

\definecolor{grey}{rgb}{.7,.7,.7}
\definecolor{evidGP}{rgb}{0,0,1}
\definecolor{evidG}{rgb}{0,0.5,0}

\newcommand{\Lip}{\mathrm{Lip}}

\AtBeginDocument{
  \let\div\relax
  \DeclareMathOperator{\div}{div}
}
\newcommand{\res}{\mathop{\hbox{\vrule height 7pt width .5pt depth 0pt
\vrule height .5pt width 6pt depth 0pt}}\nolimits}

\newcommand{\A}{\boldsymbol{A}}

\begin{document}

\title[Beyond $BV$: a new pairing and Gauss-Green formulas]
{Beyond $BV$: new pairings and Gauss-Green formulas \\ for measure fields with divergence measure}

\author[G.~E.~Comi]{Giovanni E. Comi}
\address[G.~E.~Comi]{Dipartimento di Matematica, Università di Bologna, Piazza di Porta San Donato 5, 40126 Bologna (BO), Italy}
\email{giovannieugenio.comi@unibo.it}
\author[V.~De Cicco]{Virginia De Cicco}
\address[V.~De Cicco]{Dipartimento di Scienze di Base  e Applicate per l'Ingegneria, Sapienza Università di Roma\\
	Via A.\ Scarpa 10 -- I-00185 Roma (Italy)}
\email{virginia.decicco@uniroma1.it}
\author[G.~Scilla]{Giovanni Scilla}
\address[G.~Scilla]{Dipartimento di Matematica ed Applicazioni ``R.Caccioppoli'', Università degli Studi di Napoli Federico II\\
	Via Cintia, Monte Sant'Angelo -- I-80126 Napoli (Italy)}
\email{giovanni.scilla@unina.it}

\thanks{\textit{Acknowledgments}. 
The authors are members of  the Istituto Nazionale di Alta Matematica (INdAM), Gruppo Nazionale per l'Analisi Matematica, la Probabilità e le loro Applicazioni (GNAMPA), and were partially supported by the INdAM--GNAMPA 2022 Project \textit{Alcuni problemi associati a funzionali integrali: riscoperta strumenti classici e nuovi sviluppi}, codice CUP\_E55F22000270001, the INdAM--GNAMPA 2023 Project \textit{Problemi variazionali degeneri}, codice CUP\_E53C22001930001 and the INdAM--GNAMPA 2023 Project \textit{Problemi variazionali per funzionali e operatori non-locali}, codice CUP\_E53\-C22\-001\-930\-001.
Part of this work was undertaken while the first author was visiting the Universit\`a Sapienza and SBAI Department in Rome. He would like to thank these institutions for the support and warm hospitality during the visits. The authors would like to warmly thank Professor Gianni Dal Maso for his insightful comments on a preliminary version of this paper. The authors would also like to thank Doctor Anna Nikishova for Figure \ref{fig-4.1}. Finally, the authors would like to express their gratitude to the reviewer for their helpful comments.
}


\begin{abstract}
A new notion of pairing between measure vector fields with divergence measure and scalar functions, which are not required to be weakly differentiable, is introduced. In particular, in the case of essentially bounded divergence-measure fields, the functions may not be of bounded variation. This naturally leads to the definition of $BV$-like function classes on which these pairings are well defined. Despite the lack of fine properties for such functions, our pairings surprisingly preserve many features of the recently introduced $\lambda$-pairings \cite{CDM}, as coarea formula, lower semicontinuity, Leibniz rules, and Gauss-Green formulas. Moreover, in a natural way new anisotropic ``degenerate'' perimeters are defined, possibly allowing for sets with fractal boundary. 
\end{abstract}

\maketitle

\tableofcontents

\section{Introduction}
Starting from the seminal paper of Anzellotti \cite{Anz}, in the last decades there has been a growing interest
in giving a well posed definition of {\it pairing}, i.e. a scalar product between a vector field and the weak gradient of a suitably regular function, motivated by the integration by parts.
Indeed, many mathematicians have made effort to look for formulas of this type in increasingly general contexts and 
to establish the validity of the Gauss-Green formula under very weak regularity assumptions 
(see \cite{AmbCriMan, Anz2, ChCoTo, ChIrTo, CCDM, ComiLeo, ComiPayne, CD3, crasta2019extension, CDM, Ir, Silh, MR2532602, Silhavy19}).
 The motivation of this large interest is that the notion of pairing is a fundamental tool and is widely used in several contexts, e.g. in the study of Dirichlet problems for the Prescribed Mean Curvature and $1$-Laplace equations, as well as related topics (see \cite{MR1814993, AVCM, Cas, MR2348842, LeoComi, LeoSar, LeoSar2, MR2502520, MR3501836, MR3813962}). In addition, we point to some recent extensions of the notion of pairing to non-Euclidean frameworks \cite{ComiMagna, BuffaComiMira, MR4381316, gorny2023anzellotti}.

The aim of this paper is to introduce a very general pairing $(\A, Du)$ which entails two significant features at the same time: a measure (instead of $L^\infty$) vector field $\A$ with  divergence measure and 
a function $u$ which can be beyond the class of functions of bounded variation. It extends all the notions of pairings given previously in the literature, possibly paving the way to integration by parts formulas on sets with fractal boundary.\\
\\
\noindent
\emph{A brief literature review about the classical pairing.} We briefly review the classical notion of pairing. To this purpose, we recall the definition of {\it{divergence-measure fields}}.
We consider vector fields $\A \in L^{p}(\Omega; \R^{N})$, for some open set $\Omega \subset \R^{N}$ and $p\in [1,+\infty]$, with distributional divergences $\div \A$ represented by scalar valued measures belonging to $\mathcal{M}(\Omega)$; that is, the space of finite Radon measures on $\Omega$. Their space is denoted by $\DM^{p}(\Omega)$, and its local version by $\DM^p_{\rm loc}(\Omega)$. Then,  the pairing between $\A\in\DM^{\infty}_{\rm loc}(\Omega)$ and $D u$, for a given function $u \in BV_{\rm loc}(\Omega) \cap L^{\infty}_{\rm loc}(\Omega)$ is the distribution given by
\begin{equation}\label{def_pairing}
 \varphi \in C^{\infty}_{c}(\Omega) \to \ban{(\A, Du), \varphi} := - \int_{\Omega}\varphi u^* \, d \div \A - \int_{\Omega} u \A \cdot \nabla \varphi \, dx,
 \end{equation}
where $u^{*}$ is the precise representative of the $BV$-function $u$. Since the measure $\div \A$ is absolutely continuous with respect to the $(N-1)$-Hausdorff measure $\Haus{N-1}$, and $u^*$ is well defined $\Haus{N-1}$-almost everywhere, this definition is well-posed.
As proven in  \cite{ChenFrid}, the vector field $u \A$ belongs to $\DM^{\infty}_{\rm loc}(\Omega)$ and the pairing is a measure satisfying the following Leibniz-type formula:
\begin{equation*} (\A, Du) = \div (u \A) - u^{*} \div \A. \end{equation*}

A decisive step in the study of the pairing was done in \cite{CD3}: here the authors, redefining the pairing distribution 
for functions $u \in BV_{\rm loc}(\Omega)$ such that $u^* \in L^1_{\rm loc}(\Omega, |\div \A|)$ (as it was firstly done in \cite{MR3939259}, thus requiring a dependence between the scalar functions and the vector field), 
proved a coarea formula and a Leibniz rule, and finally achieved generalizations of the Gauss--Green formulas.

Later on, motivated by obstacle problems in $BV$ \cite{MR3501836,MR3813962} and semicontinuity issues, for which the standard pairing is not adequate, in \cite{CDM} the authors introduced a new family of pairings, depending on the choice of the pointwise representative of $u$ (afterwards, this tool has proven to be very useful for the Prescribed Mean Curvature problem with measure data, see \cite{ComiLeo, LeoComi}). More precisely, they proved that, given $\A \in \DM^{\infty}_{\rm loc}(\Omega)$ and $u \in BV_{\rm loc}(\Omega)$ such that $u^* \in L^1_{\rm loc}(\Omega, |\div \A|)$, for every Borel function $\lambda\colon\Omega \to [0,1]$ there exists a measure
$\pair{\A, Du}_\lambda$ defined as
\begin{equation} \label{def_lambdapairing}
 \varphi \in C^{1}_{c}(\Omega) \to \ban{(\A, Du)_\lambda, \varphi} := - \int_{\Omega}\varphi u^\lambda \, d \div \A - \int_{\Omega} u \A \cdot \nabla \varphi \, dx, \end{equation}
which then satisfies the Leibniz rule
\begin{equation}\label{f:genp}
\Div(u \A) = {u}^\lambda\, \Div\A + \pair{\A, Du}_\lambda,
\end{equation}
where
\begin{equation} \label{eq:intro_u_lambda_classic}
{u}^\lambda := (1-\lambda)\umeno+\lambda\upiu,
\end{equation} 
and $u^\mp$ are the approximate liminf and limsup of $u$. Indeed, the approximate liminf and limsup of a functions with bounded variation are well defined $\Haus{N-1}$-almost everywhere. In particular, if $\lambda(x) \equiv \frac{1}{2}$, we get $u^{\frac{1}{2}}(x) = u^*(x)$ for $\Haus{N-1}$-a.e. $x \in \Omega$. In addition, the flexibility of choice of $\lambda$ allows to select a suitable convex combination of $u^-$ and $u^+$, with a possibly non constant weight, and this proved to be useful in the applications above mentioned (for further details, we refer to \cite{CDM,LeoComi}).

As for divergence-measure fields not in $L^\infty$, problems related to the foundations of continuum mechanics (especially concerning the representation of Cauchy fluxes) naturally led to weaker versions of the Gauss--Green formulas for $\DM^1_{\rm loc}$-fields and sets satisfying suitable measure-theoretic assumptions (see \cite{ChCoTo, ChFr1, DGMM, Schu, Silh, MR2532602, Silhavy19}). As a further step in this direction, it was introduced the larger space $\DM(\Omega)$ of measure-valued vector fields $\A \in \mathcal M(\Omega; \R^{N})$ whose distributional divergences $\div \A$ belong to $\mathcal{M}(\Omega)$; and analogously the local version $\DM_{\rm loc}(\Omega)$. These vector fields are sometimes called {\it{extended}} divergence-measure fields.
In this general setting, another definition of pairing may be found in \cite{MR2532602}: given $\A \in \DM_{\rm loc}(\Omega)$ and $u \in \Lip_{\rm loc}(\Omega)$, we have $u\A\in\DM_{\rm loc}(\Omega)$ and
\begin{equation}
\Div(u\A) = u\Div \A + \left\langle\left\langle\nabla u,\A\right\rangle\right\rangle, 
\label{eq:intro_pairingShilavy}
\end{equation}
where $\left\langle\left\langle\nabla u,\A\right\rangle\right\rangle$ is a Radon measure. In particular, thanks to \cite{ChFr1}, if $u \in C^1(\Omega)$, we see that
\begin{equation*}
\left\langle\left\langle\nabla u,\A\right\rangle\right\rangle = \nabla u \cdot \A,
\end{equation*}
thus showing that this is a generalization of the scalar product between a continuous function and a vector valued measure. We refer also to \cite{Ir} for more recent developments on the study of this pairing.\\
\\
\noindent
\emph{Description of our results.} The initial goal of our research was to further weaken the regularity assumption on the scalar function $u$, also in view of possible applications to conservation laws (see \cite{ChenFrid,Frid}), where the solutions may not belong to $BV$. We started with the case $\A \in \DM^\infty(\Omega)$, in order to check whether the pairing could be defined outside of the natural class of $BV$ functions. Then, we noticed that our assumptions did not actually depend on the summability of the vector fields, and that there was no reason not to directly consider measure valued divergence-measure fields. Hence, we were able to develop a unifying approach to the theory of pairings, which takes into account all the previously known results and particular cases. In addition, we provided a functional framework where to study classical problems which not necessarily admit solutions in the usual Sobolev and $BV$ spaces.

Given a vector field $\A\in \DM_{\rm loc}(\Omega)$ and a Borel function $\lambda\colon \Omega \to [0,1]$, we define the class $X^{\A,\lambda}(\Omega)$ of those equivalence classes of Borel functions $u$ such that
\begin{equation*}
u^{\lambda} \in L^1(\Omega, |\A|) \cap L^1(\Omega, |\Div\A|) \,,
\end{equation*}
and we denote by $X^{\A,\lambda}_{\rm loc}(\Omega)$ its local version.
Here, $u^\lambda(x)$ is defined as in \eqref{eq:intro_u_lambda_classic} for all $x \in \Omega \setminus Z_u$, where
$Z_u$ is the set where both $u^+$ and $u^-$ are infinite with different sign. This definition means that $u^\lambda$ is the convex combination of $u^+$ and $u^-$: in order to extend it even on $Z_u$, we set 
\begin{equation*}
u^\lambda(x) = (2 \lambda(x) - 1) \cdot (+ \infty) \ \text{ if } x \in Z_u
\end{equation*} 
where we convene that $0 \cdot (\pm \infty) = 0$ (see \eqref{f:pr} and the subsequent comments).
Note that, differently from the previous approach in \cite{CDM}, $u$ does not need to be an $L^1$ function, and $u^+(x)$, $u^-(x)$ are defined \emph{for every $x\in\Omega$} through the densities of sublevels of $u$, for which we refer to \eqref{def:traces_Maz} and the remarks below it. However, $u^\lambda$ has to be a Borel function, in order to ensure that its integration against the Radon measures $\A$ and $\Div \A$ is well posed. We refer to Remark \ref{rem:vecchio_approccio} for a detailed comparison between our approach and the previous one.

Then, for $u\in X^{\A,\lambda}_{\rm loc}(\Omega)$, we define the $\lambda$-pairing between $\A$ and $u$ as the distribution
$\pair{\A,Du}_\lambda$ 
acting as
\begin{equation} \label{eq:def_lambda_pairingintro}
\pscal{\pair{\A,Du}_\lambda}{\varphi}:= -\int_\Omega {u^\lambda\,\varphi} \, d \div \A - \int_\Omega u^\lambda \nabla \varphi \cdot d \A
\quad \text{ for } \
\varphi\in C^\infty_c(\Omega)\,.
\end{equation}
We note that, unless otherwise specified, $Du$ is only a distributional gradient. To the best of our knowledge, \eqref{eq:def_lambda_pairingintro} is the first definition of pairing which extends to a more general setting all the notions of pairing available in the literature, often not comparable each other (see Section \ref{ss:div}).

In analogy with the classical theory of functions of bounded variation, we introduce a version of $BV$-type classes, with respect to the variation given by this generalized $\lambda$-pairing, by setting
\begin{equation*}
BV^{\A,\lambda}(\Omega) := \left\{
u\in X^{\A,\lambda}(\Omega): (\A,Du)_{\lambda}\in {\mathcal{M}(\Omega)}
\right\}\,.
\end{equation*}
Then, if $u\in BV^{\A,\lambda}(\Omega)$, we prove that $u^\lambda \A \in \DM(\Omega)$ and the following Leibniz rule holds
\begin{equation*}
(\A,Du)_\lambda=-u^\lambda\,\Div\A+\Div(u^\lambda \A) \ \text{ on } \ \Omega,
\end{equation*}
in the sense of Radon measures (Proposition~\ref{prop:main_inclusions}). We then investigate many of the fundamental properties of $BV^{\A,\lambda}(\Omega)$, and we introduce its related Sobolev-type class
\begin{equation*}
W^{\A,\lambda}(\Omega) := \left\{
u\in BV^{\A,\lambda}(\Omega): (\A,Du)_{\lambda} \ll \Leb{N}
\right\}\,.
\end{equation*}

Motivated by the classical theory of the calculus of variations, we investigated whether our objects enjoy the necessary topological properties to obtain the existence of minimizers of functionals involving the total variation of the $\lambda$-pairings. We show in Theorem~\ref{t:lscnuovo11} that a suitable notion of convergence in $BV^{\A,\lambda}(\Omega)$ which ensures the lower semicontinuity of the pairing is the \lq\lq $(\A, \lambda)$-convergence'' (see Definition \ref{d:convergence}). 
That said, we have to notice that $BV^{\A, \lambda}(\Omega)$ is not a linear space in general, since the map
\begin{equation*}
u\to (\A,Du)_\lambda
\end{equation*} 
is not linear (see Remarks \ref{rem:not_lin_pairing_short} and \ref{rem:not_lin_pairing}). 
However, 
we identify a subclass of $BV^{\A,\lambda}(\Omega)$ 
which is a linear space and, 
under some additional assumptions, a Banach space (Proposition \ref{prop:property}) with respect to a suitable norm. 
Within this functional setting, by using the direct method we prove the existence of minimizers for functionals of the type 
\begin{equation*}
\mathcal{E}(u):= |(\A,Du)_\lambda|(\Omega)\,,
\end{equation*}
up to adding a suitable fidelity term 
(see Theorem \ref{thm:min_funct_E}).

In analogy with the results established in \cite{CDM} for the $\lambda$-pairing, we obtain some coarea-type formula (see Theorem \ref{t:coarea}).

The notion of $\lambda$-pairing suggests in a natural way the definition of an \lq\lq anisotropic degenerate'' perimeter, the \emph{$(\A, \lambda)$-perimeter}, as the total variation of the pairing; i.e., 
$$
P_{\A,\lambda}(E,\Omega)=|(\A,D\chi_E)_\lambda|(\Omega).
$$
This perimeter could be relevant for many applications in which one is indeed not interested in the growth of some quantity $u$ outside of a selected field of directions.
We immediately notice that a set of finite (Euclidean) perimeter has finite $(\A, \lambda)$-perimeter, but the converse is, in general, not true (see Remark~ \ref{rem:sizesupp}). 

The $(\A, \lambda)$-perimeter enjoys some typical properties of the classical perimeter: absolute continuity properties (Proposition \ref{prop:abs_cont_per} and Proposition \ref{prop:abs_cont_per_unbdd}); 
locality and additivity on suitably \lq \lq disjoint'' sets (Proposition \ref{prop:per_prop}) and the lower semicontinuity with respect to the $(\A, \lambda)$-convergence (Proposition \ref{prop:per_lssc}). 

Another natural question is to ask whether the $(\A, \lambda)$-perimeter is concentrated on some type of generalized boundary of the set. To this purpose, we recall the definition of another type of Lebesgue measure-invariant boundary of a Borel set $E$:
\begin{equation*}
\partial^- E := \{ x \in \R^N : 0 < \Leb{N}(E \cap B_r(x)) < \Leb{N}(B_r(x)) \ \text{ for all } r > 0\}.
\end{equation*}
Then we have ${\rm supp}((\A,D\chi_E)_\lambda) \subseteq \partial^- E$ (see Proposition \ref{prop:supp_A_lambda_per}), even though such a control on the size of the support of the pairing distribution is in general too large (see Remark \ref{rem:sizesupp}).

In particular, if $E$ is a set of finite perimeter in $\Omega$, we retrieve the representation formula
\begin{equation*}
(\A, D \chi_E)_\lambda = \left ( (1 - \lambda) {\rm Tr}^i (\A,\partial^*E) + \lambda {\rm Tr}^e (\A,\partial^*E) \right ) \Haus{N-1} \res \partial^{*} E,
\end{equation*}
where $\partial^*E$ is the {\em measure theoretic boundary} of $E$\footnote{We point out that in the literature \cite{evans2015measure, Maggi} it is sometimes De Giorgi's {\em reduced boundary} to be denoted as $\partial^* E$; however, given that we do not employ it in this paper, we use this notation for Federer's measure theoretic boundary.} and $${\rm Tr}^i (\A,\partial^*E), {\rm Tr}^e (\A,\partial^*E) \in L^\infty(\partial^* E, \Haus{N-1})$$ are the {\em interior and exterior normal traces} of $\A$ on $\partial^* E$, see \eqref{eq:def_normal_traces}.

For sets $E$ with finite $(\A,\lambda)$-perimeter, we prove a general Gauss-Green formula, Theorem~\ref{thm:GG}. In order to make a comparison with the classical results, we only mention here that, for Borel sets $E \Subset \Omega$, if $\lambda\equiv0$ and $\chi_E \in BV^{\A, 0}(\Omega)$, then
\begin{equation*}
\div \A(E^1) = - \int_{\partial^- E} \, d (\A, D \chi_E)_0 \,,
\end{equation*}
while if $\lambda\equiv1$ and $\chi_E \in BV^{\A, 1}(\Omega)$, then
\begin{equation*}
\div \A(E^1 \cup \partial^* E) = - \int_{\partial^- E}
 \, d (\A, D \chi_E)_1 \, ,
\end{equation*}
where the pairings on the right hand sides play the role of the traces on $\partial^-E$. 
This is a simplified version ($u=1$) of the following main result (see Theorem \ref{thm:IBP_lambda_gen}). 
\begin{theorem} \label{thm:IBP_lambda_genintro}
Let $\A \in \DM_{\rm loc}(\Omega)$, $\lambda : \Omega \to [0, 1]$ be a Borel function and $u \in BV^{\A, \lambda}(\Omega)$. Let $E \subset \Omega$ be a Borel set. If $\chi_E \in BV^{u^\lambda \A, 0}(\Omega)$ and ${\rm supp}(\chi_E^{-} u^{\lambda} |\A|) \Subset \Omega$, then we have
\begin{equation} \label{eq:IBP_iintro}
\int_{E^1} u^\lambda \, d \div \A + \int_{E^1} d (\A, Du)_\lambda = - \int_{\partial^- E} \, d (u^\lambda \A, D \chi_E)_0 \, ;
\end{equation}
while, if $\chi_E \in BV^{u^\lambda \A, 1}(\Omega)$ and ${\rm supp}(\chi_E^{+} u^{\lambda} |\A|) \Subset \Omega$, then we have
\begin{equation} \label{eq:IBP_eintro}
\int_{E^1 \cup \partial^*E} u^\lambda \, d \div \A + \int_{E^1 \cup \partial^*E} d (\A, Du)_\lambda = - \int_{\partial^- E} \, d (u^\lambda \A, D \chi_E)_1 \,.
\end{equation}
\end{theorem}
In particular, we can integrate on sets which do not have local finite perimeter, given that the only assumptions are $(\A, D \chi_E)_0, (\A, D \chi_E)_1 \in \mathcal{M}(\Omega)$. 

Finally, in Section~\ref{sec:1D} we give an insight in the one-dimensional case $N=1$, which is related to research carried out in the recent papers \cite{ChDCMe1, ChDCMe2}. 
\\

\noindent
\emph{Outline of the paper.} In Section~\ref{s:prelim} we fix the basic notation and recall some definitions and preliminary results about measures and distributions, approximate limits and representatives, and divergence-measure fields. In Section~\ref{s:spaceBVA} we introduce the class $BV^{\A,\lambda}$ and the new notion of $\lambda$-pairing, then we investigate its basic properties, in particular the absolute continuity, and we define a related \lq\lq degenerate'' Sobolev class. In Section~\ref{s:lsc} we prove the lower semicontinuity of the $\lambda$-pairing functional with respect to the $(\A,\lambda)$-convergence and discuss some smooth approximation results. In Section~\ref{s:banach} we identify a linear space contained in $BV^{\A, \lambda}$ and briefly investigate its properties. Section~\ref{s:coarea} focuses on a coarea-type formula for the $\lambda$-pairing, while in Section~\ref{s:perimeter} we introduce a corresponding notion of perimeter and prove some of its features. In Section~\ref{s:gauss-green} we establish general Gauss-Green and integration by parts formulas. Eventually, the last Section~\ref{sec:1D} deals with the case $N=1$.

\section{Notation and preliminary results}
\label{s:prelim}

In the following {we denote by} \(\Omega\) a nonempty open subset of 
\(\R^N\), and for every {set} $E\subset \R^N$ {we denote by} $\chi_{E}$ its 
characteristic function. For $x \in \R^N$ and $r > 0$, we denote by $B_{r}(x)$ the ball centered in $x$ with radius $r$, and we set $B_1 := B_1(0)$. Given a set $U$, we denote its closure by $\overline{U}$. We say that a set $E$ is compactly contained in $\Omega$, and we write $E \Subset \Omega$, if $\overline{E}$ is a bounded set and $\overline{E}\subset \Omega$.

\subsection{Measures and distributions}

The following definitions and basic facts about measures can be found, e.g., in \cite[Chapter 1]{AFP}.

We denote by
$\Leb{N}$ 
and $\Haus{\alpha}$
the Lebesgue measure 
and the $\alpha$-dimensional 
Hausdorff measure in $\R^N$ for some $\alpha \in [0, N]$, respectively. Unless otherwise stated, a measurable set is a $\Leb{N}$-measurable set. We set $\omega_N := \Leb{N}(B_1) = \frac{\pi^{\frac{N}{2}}}{\Gamma\left (\frac{N}{2} + 1\right )}$, where $\Gamma$ is Euler's Gamma function, so that $\Haus{N-1}(\partial B_1) = N \omega_N$.

Following the notation of \cite{AFP}, we denote by $\mathcal{M}_{\rm loc}(\Omega)$ the space of Radon measures on $\Omega$, and by $\radon$ the space of finite Radon measures on $\Omega$. Analogously, we denote by $\mathcal{M}_{\rm loc}(\Omega; \R^N)$ and $\mathcal{M}(\Omega; \R^N)$ the spaces of $\R^N$-valued measures whose components belong to $\mathcal{M}_{\rm loc}(\Omega)$ and $\mathcal{M}(\Omega)$, respectively.

Given $\mu\in\mathcal{M}_{\rm loc}(\Omega)$ and a $\mu$-measurable set $E$, the {\em restriction} 
$\mu\res E$
is the Radon measure defined by
\[
\mu\res E(B)=\mu(E\cap B), \qquad \forall\ B\ \text{$\mu$-measurable},\ 
B\subset\Omega.
\]

The {\em total variation} $|\mu|$ of $\mu \in \mathcal{M}_{\rm loc}(\Omega)$ is 
the nonnegative Radon measure defined by
\[
|\mu|(E) := \sup\left\{ \sum_{h=0}^\infty |\mu(E_h)| \colon \ E_h\ 
\text{$\mu$-measurable sets, pairwise disjoint},\ E=\bigcup_{h=0}^\infty E_h 
\right\},
\]
for every $\mu$-measurable set $E \Subset \Omega$. If $\mu \in \mathcal{M}(\Omega)$, then $|\mu|(\Omega) < \infty$.

\medskip  

A measure $\mu \in \mathcal{M}_{\rm loc}(\Omega)$ is {\em absolutely continuous} with respect to a given
nonnegative measure $\nu$ (notation: $\mu \ll \nu$) 
if $|\mu|(B)=0$ for every Borel set $B$ such that $\nu(B)=0$. 
Two positive measures $\nu_1, \nu_2\in \mathcal{M}_{\rm loc}(\Omega)$ are {\em mutually 
singular}
(notation: $\nu_1 \perp \nu_2$) 
if there exists a Borel set $E$ such that $|\nu_1|(E)=0$ and 
$|\nu_2|(\Omega\setminus E) = 0$. Given $\mu\in \mathcal{M}_{\rm loc}(\Omega)$, the \emph{Lebesgue decomposition} of $\mu$ with respect to $\Leb{N}$ is 
\begin{equation*}
\mu=\mu^a + \mu^s\,,
\end{equation*}
where $\mu^a$ is the {\em absolutely continuous part}, satisfying $\mu^a \ll \Leb{N}$, and $\mu^s$ is the {\em singular part}, satisfying $\mu^s \perp \Leb{N}$.

A measure $\mu$ in $\Omega$ is \emph{concentrated on} $E\subset\Omega$ if $|\mu|(\Omega\backslash E)=0$. The intersection of the closed sets $E\subset\Omega$ such that $\mu$ is concentrated on $E$ is called the \emph{support of $\mu$} and is denoted by ${\rm supp}(\mu)$. In particular,
\begin{equation*}
\Omega\backslash {\rm supp}(\mu) = \left\{x\in\Omega:\,\, \mu(B_r(x))=0 \mbox{ for some $r>0$} \right\}\,.
\end{equation*}

We define a space of equivalence classes of Borel measure functions in the following way:
\begin{equation*}
\Borel = \{ u : \Omega \to \R \text{ Borel measurable} \}/ \sim,
\end{equation*}
where $\sim$ is the equivalence relation given by the almost everywhere equality; that is,  
\begin{equation*}
u \sim v \ \iff \ \Leb{N}(\{ x \in \Omega : u(x) \neq v(x)  \}) = 0.
\end{equation*}

Now, we recall some basic definitions about distributions. We refer the interested reader to \cite[Chapter 6]{Rudin} for a detailed treatment of this topic. We denote by $\mathscr{D}(\Omega)$ the space of \emph{test functions}; i.e., $\varphi\in \mathscr{D}(\Omega)$ if and only if $\varphi\in C^\infty(\Omega)$ and the support of $\varphi$ is a compact subset of $\Omega$. We consider the norms
\begin{equation*}
\|\varphi\|_n:=\max\{|D^\alpha\varphi(x)|:\,\, x\in\Omega,\, \alpha \in \N^n, \, |\alpha|\leq n\}
\end{equation*}
for $\varphi\in \mathscr{D}(\Omega)$ and $n \in \N$, where, corresponding to the multi-index $\alpha=(\alpha_1,\alpha_2,\dots,\alpha_n)$, $D^\alpha$ denotes the differential operator of order $|\alpha|:=\alpha_1+\alpha_2+ \dots + \alpha_n$ defined by $D^\alpha :=(\partial_{x_1})^{\alpha_1}(\partial_{x_2})^{\alpha_2}\dots (\partial_{x_n})^{\alpha_n}$. It is well-known (see \cite[Theorem 6.4 and 6.5]{Rudin}) that $\mathscr{D}(\Omega)$ is a topological vector space when equipped with a suitable topology $\tau$ for which all Cauchy sequences do converge. A linear functional $\Lambda$ on $\mathscr{D}(\Omega)$ which is continuous (with respect to the topology $\tau$) is called a \emph{distribution} in $\Omega$. The space of all distributions in $\Omega$ is denoted by $\mathscr{D}'(\Omega)$. The smallest integer $n \in \N$ such that
\begin{equation*}
|\Lambda(\varphi)|\leq C \|\varphi\|_n
\end{equation*}
for every $\varphi\in \mathscr{D}(\Omega)$, if exists, is called the \emph{order} of $\Lambda$.

\begin{remark}
Given a $\mu\in \mathcal{M}_{\rm loc}(\Omega)$, setting
\begin{equation*}
\Lambda_\mu(\varphi):=\int_\Omega \varphi \, d\mu\,,\,\, \,\,\,\,\varphi\in \mathscr{D}(\Omega)\,,
\end{equation*}
defines a distribution $\Lambda_\mu$ in $\Omega$ of order zero. Conversely, given a distribution $\Lambda$ in $\Omega$ of order zero,
by the Riesz Representation Theorem there exists $\mu\in \mathcal{M}_{\rm loc}(\Omega)$ such that $\Lambda=\Lambda_\mu$.
\end{remark}

Let $\Lambda\in \mathscr{D}'(\Omega)$, and set
\begin{equation*}
W:=\bigcup \left\{U\subset\Omega\,,\,\, U \mbox{ open: } \Lambda(\varphi)=0\,,\,\,\mbox{ for all }\varphi\in \mathscr{D}(U) \right\} \,.
\end{equation*}
Then, the \emph{support of $\Lambda$} is defined as ${\rm supp}(\Lambda):=\Omega\backslash W$. In the case of distributions of order zero, this definition coincides with the one of the support of a measure; that is, ${\rm supp}(\Lambda_\mu) = {\rm supp}(\mu)$.

Finally, by exploiting the relation between measures and distributions just recalled, we say that a sequence $(\mu_k) \subset \mathcal{M}(\Omega)$ \emph{weakly converges} to some $\mu \in \mathcal{M}(\Omega)$, and we write $\mu_k \weakto \mu$ in $\mathcal{M}(\Omega)$, if
\begin{equation*}
\lim_{k \to + \infty} \int_\Omega \varphi \, d \mu_k = \int_\Omega \varphi \, d \mu \ \text{ for all } \varphi \in C_c(\Omega).
\end{equation*}

\medskip

\subsection{Approximate limits and $\lambda$-representatives}
\label{ss:def}

The following basic definitions and results can be found, e.g., in \cite[Sections 3.6 and 4.5]{AFP}.

We say that a function \(u\in L^1_{{\rm loc}}(\Omega)\) has an {\em approximate limit} 
\(z\in\R\) at
$x\in\Omega$ if
\begin{equation}
\lim_{r\rightarrow0^{+}}\frac{1}{\Leb{N}\left(  B_r(x)\right)}\int_{B_r\left(  
x\right)
}\left|  u(y)  -z  \right|  \,dy=0\,;
\label{eq:approxlim1}
\end{equation}
in this case we say that $x$ is a {\em Lebesgue point} of $u$.
The set $S_u\subset\Omega$ of points where this property does not hold is called the
{\em approximate discontinuity set} of $u$, and, thanks to Lebesgue's differentiation theorem, we know that $\Leb{N}(S_u) = 0$.
For any $x\in \Omega \setminus S_u$ the approximate limit $z$ is uniquely 
determined and is denoted by $z=:\tilde{u}(x)$. Note that Chebychev inequality and \ref{eq:approxlim1} imply
\begin{equation}
\lim_{r\rightarrow0^{+}} \frac{\Leb{N}\left(\left\{y\in B_r(x):\,\, |u(y)-\tilde{u}(x)|>\varepsilon \right\}\right)}{\Leb{N}(B_r(x))}=0
\label{eq:approxlim2}
\end{equation}
for every $\varepsilon>0$. In fact, \eqref{eq:approxlim2} provides an alternative (weaker) definition of approximate limit for a Borel measurable (even non locally summable) function (see \cite[Remark 4.29]{AFP} and (see \cite[\S 2.9.12]{FED})), and the two definitions are equivalent for locally bounded functions (see \cite[Proposition 3.65]{AFP}).

If $u = \chi_E$, for a measurable set $E\subset\R^N$, then the approximate limit at a point $x\in\R^N$ is also called {\em density} of $E$ at $x$, and it is given by
\[
D(E;x) := \lim_{r \to 0^+} \frac{\Leb{N}(E\cap B_r(x))}{\Leb{N}(B_r(x))}
\] 
whenever this limit exists. We call {\em measure theoretic interior} of $E$ the set of points with density 1, and we denote it by 
$$E^1 := \{ x \in \R^N : D(E; x) = 1 \}.$$
We call {\em measure theoretic boundary} of $E$ the approximate discontinuity set of $\chi_E$, and we denote it by $\partial^*E := S_{\chi_E}$, which also satisfies $\partial^* E = \R^N \setminus (E^1 \cup (\R^N \setminus E)^1)$.

Set $\overline{\mathbb{R}}:=\mathbb{R}\cup\{\pm \infty\}$. Given a Borel measurable function $u : \Omega \to \R$, we denote the {\em sublevel and superlevel sets} of $u$ as
\begin{equation*}
\{u < t\} = \{ x \in \Omega : u(x) < t \} \ \text{ and } \ \{u > t\} = \{ x \in \Omega : u(x) > t \},
\end{equation*}
and we recall the definition of the {\em approximate liminf and limsup} at a point $x \in \Omega$ for $u \in \Borel$:  
\begin{equation} \label{def:traces_Maz}
u^-(x) :=
\sup\left\{t\in\overline{\mathbb{R}}\colon
D(\{u < t\};x) = 0\right\},
\quad
u^+(x) :=
\inf\left\{t\in\overline{\mathbb{R}}\colon
D(\{u > t\};x) = 0\right\}
\end{equation}
(see \cite[Definition 4.28]{AFP}).
We notice that $u^+, u^- : \Omega \to [- \infty, + \infty]$ are Borel measurable functions and the set $S_u^* := \{ x \in \Omega : u^-(x) < u^+(x) \}$ satisfies 
\begin{equation} \label{eq:S_u_*_negl}
\Leb{N}(S_u^*) = 0,
\end{equation} 
so that $u^+(x) = u^-(x)$ for $\Leb{N}$-a.e. $x \in \Omega$, by \cite[Definition 4.28]{AFP} and the comments below. 
In the particular case $u \in L^{1}_{\rm loc}(\Omega)$, \eqref{eq:approxlim2} implies that 
\begin{equation} \label{eq:u_+_-_tilde}
u^+(x) = u^-(x) = \tilde{u}(x) \ \text{ for all } x \in \Omega \setminus S_u,
\end{equation} 
which implies $S_u^* \subset S_u$.
Therefore, in $\Omega \setminus S_u^*$ we shall write $\tilde{u}(x) := u^+(x) = u^-(x)$, with a little abuse of notation. If $u \in \Borel$, we tacitly identify the canonical representative of the class with $\tilde{u}$, by setting $\tilde{u} = 0$ on $S_u^*$. However, by definition of $\Borel$ it is clear that $u = \tilde{u}$ with respect to the Lebesgue measure, so that, whenever dealing with $\Leb{N}$, we shall simply write $u$ as the representative of the class.

Given $u \in L^1_{\rm loc}(\Omega)$, we say that \(x\in\Omega\) is an {\em approximate jump point} of \(u\) if
there exist \(a,b\in\R\), \(a\neq b\), and a unit vector \(\nu\in\R^N\) such that 
\begin{equation}\label{f:disc}
\begin{gathered}
\lim_{r \to 0^+} \frac{1}{\Leb{N}(B_r^i(x))}
\int_{B_r^i(x)} |u(y) - a|\, dy = 0,
\\
\lim_{r \to 0^+} \frac{1}{\Leb{N}(B_r^e(x))}
\int_{B_r^e(x)} |u(y) - b|\, dy = 0,
\end{gathered}
\end{equation}
where \(B_r^i(x) := \{y\in B_r(x):\ (y-x)\cdot \nu > 0\}\), and 
\(B_r^e(x) := \{y\in B_r(x):\ (y-x)\cdot \nu < 0\}\).
The triplet \((a,b,\nu)\), uniquely determined by \eqref{f:disc} 
up to a permutation
of \((a,b)\) and a change of sign of \(\nu\),
is denoted by \((\uint(x), \uext(x), \nu_u(x))\).
The set of approximate jump points of \(u\) is denoted by \(J_u\), and it clearly satisfies $J_u \subset S_u$. We notice that, even for a function $u \in L^1_{\rm loc}(\Omega)$, the jump set $J_u$ is $(N-1)$-rectifiable (see \cite{DELNIN}), while $S_u$ can have high Hausdorff dimension, even equal to the space dimension $N$.

We point out that (cfr. \cite[Definition 4.30]{AFP} and the comments below therein)
\begin{equation} \label{eq:u_+_-_i_e}
\umeno(x) =\min\{\uint(x),\uext(x)\} \ \text{ and } \ 
\upiu(x) =\max\{\uint(x),\uext(x)\}
\quad \text{ for all } x\in J_u.
\end{equation}

Finally, for $u \in L^1_{\rm loc}(\Omega)$ we define the {\em precise representative} of $u$ in $x \in \Omega$ as
\begin{equation} \label{def:precise_repr}
u^{*}(x) := \lim_{r\to0^+}\frac{1}{\Leb{N}\left(  B_r(x)\right)} \int_{B_r(x)}u(y) \, d y,
\end{equation}
whenever the limit exists. It is therefore evident that 
\begin{equation}\label{f:pr2}
u^*(x)=
\begin{cases}
\tilde{u}(x) & \text{ if } x\in \Omega \setminus S_u, \\
\displaystyle \frac{u^i(x)+ u^e(x)}{2} & \text{ if } x\in J_u.
\end{cases}
\end{equation} 
A priori, it is not clear whether $u^*$ is well posed in $S_u \setminus J_u$, in general. However, for sufficiently regular functions it is known that $S_u \setminus J_u$ is suitably small. To this purpose, we recall that $u$ is a {\em function of bounded variation}, and we write $u \in BV(\Omega)$, if $u \in L^1(\Omega)$ and its distributional gradient $Du$ belongs to $\mathcal{M}(\Omega; \R^N)$; and we denote by $BV_{\rm loc}(\Omega)$ the local version of the space. For a detailed treatment of the theory of $BV$ functions, we refer the reader to the monograph \cite{AFP}. Therefore, for $u \in BV_{\rm loc}(\Omega)$, it is well known that we have $\Haus{N-1}(S_u \setminus J_u) = 0$, so that $u^*(x)$ exists for $\Haus{N-1}$-a.e $x \in \Omega$ and, up to a $\Haus{N-1}$-negligible set, is given by \eqref{f:pr2}. However, in the rest of the paper we are mostly dealing with functions not in $BV_{\rm loc}(\Omega)$ nor even locally summable, so that we will have to consider another kind of representative: this is one of the crucial technicality and novelty of our paper.

For every function $u \in \Borel$ and every
Borel function $\lambda : \Omega \to [0, 1]$, we define the {\em $\lambda$--representative} $u^\lambda: \Omega \to \overline{\R}$ as
\begin{equation}\label{f:pr}
u^\lambda(x):= \begin{cases} (1-\lambda(x))\umeno(x)+\lambda(x)\upiu(x) & \text{ if } x \in \Omega \setminus Z_u \\
+ \infty & \text{ if } x \in Z_u \text{ and } \lambda(x) > \frac{1}{2} \\
0 & \text{ if } x \in Z_u \text{ and } \lambda(x) = \frac{1}{2} \\
- \infty & \text{ if } x \in Z_u \text{ and } \lambda(x) < \frac{1}{2} \\
\end{cases} 
\end{equation}
where $Z_u := \{ x \in \Omega:  u^+(x) = + \infty \text{ and } u^-(x) = - \infty \}$.  In the particular cases $\lambda \equiv 1$ and $\lambda \equiv 0$, we simply have $u^1 := u^+$ and $u^0 := u^-$, respectively. 

We notice that in the classical settings the set $Z_u$ is either negligible or empty. Indeed, if $u \in L^1_{\rm loc}(\Omega)$, then $Z_u \subseteq S_u \setminus J_u$ and so $\Leb{N}(Z_u) = 0$, while, if $u \in L^\infty_{\rm loc}(\Omega)$, then $Z_u = \emptyset$. 
Moreover, this definition of $u^\lambda$ can be seen as a generalization of \eqref{eq:intro_u_lambda_classic}, first introduced in \cite[Eq. (2.4)]{CDM}, where the authors considered $BV$ functions (see Remark \ref{rem:vecchio_approccio} below). Indeed, in the particular case of $u \in BV_{\rm loc}(\Omega)$, we see that $\Haus{N-1}(Z_u) \le \Haus{N-1}(S_u \setminus J_u) = 0$ (see \cite[Theorem 3.78]{AFP}), so that the two definitions coincide up to an $\Haus{N-1}$-negligible set.

Furthermore, if $u \in L^1_{\rm loc}(\Omega)$, we notice that $u^\lambda(x) = \tilde{u}(x)$ for all $x \in \Omega \setminus S_u$. More in general, if $u \in \Borel$, we get 
\begin{equation} \label{eq:lambda_u_tilde_S_u_*}
u^\lambda(x) = \tilde{u}(x) \ \text{ for all } x \in \Omega \setminus S_u^*,
\end{equation} 
so that, by \eqref{eq:S_u_*_negl}, we deduce that
\begin{equation} \label{eq:u_lambda_u_tilde_Leb_N}
u^{\lambda}(x) = u(x) \ \text{ for } \Leb{N}\text{-a.e. } x \in \Omega.
\end{equation}
In addition, if $u \in L^1_{\rm loc}(\Omega)$ and $\lambda\equiv \frac{1}{2}$ on $J_u$, we get $u^\frac{1}{2}(x) = u^*(x)$ for all $x \in \Omega \setminus (S_u\setminus J_u)$, but we might have $u^*(x) \neq u^{\frac{1}{2}}(x)$ for some $x \in S_u \setminus J_u$, where the limit in \eqref{def:precise_repr} exists. In particular, if $x \in Z_u$, we might have $u^*(x) \neq 0$: consider for instance the case $N = 1$, $\Omega = (-1, 1)$ and 
\begin{equation*}
u(x) = \begin{cases} \frac{1}{\sqrt{x}} & \text{ if } x > 0 \\
0 & \text{ if } x = 0 \\
\frac{1}{\sqrt[3]{x}} & \text{ if } x < 0 
\end{cases}, 
\end{equation*}
for which we have $u^+(0) = + \infty$, $u^-(0) = - \infty$ and $u^*(0) = + \infty$.

Arguing as in the proof of \cite[Lemma 2.2]{DCFV2}, we can characterize the approximate liminf and limsup of the characteristic functions of the superlevel sets of $u$, outside of some essential discontinuity set. This result will play a fundamental role in the proof of the coarea formula (Theorem \ref{t:coarea}).

\begin{lemma} \label{funzcaracter}
Let $u\colon \Omega \to\R$ be a measurable function, $t \in \R$ and let
$$N_t := \{u^-\leq t < u^+\} \setminus\{u>t\}^{1/2}.$$
Then
\begin{equation}\label{funzcaracter0}
\chi_{\{u^{\pm}>t\}}(x)= \chi_{\{u>t\}}^\pm(x)
\qquad
\forall x\in\Omega\setminus N_t.
\end{equation}
In particular, this statement holds true for each representative of an equivalence class $u \in \Borel$.
\end{lemma}

\begin{remark}
We recall that, if $u \in BV(\Omega)$, then formula \eqref{funzcaracter0} holds $\Haus{N-1}$-a.e. in $\Omega$, since
\begin{equation*}
\qquad\qquad \Haus{N-1}(N_t)=0
\qquad \hbox{\rm for $\Leb{1}$-a.e. $t\in\bbbr,$}
\end{equation*}
see \cite[Theorems 3.40 and 3.61]{AFP}.
\end{remark}

\begin{proof}
Thanks to \cite[Lemma 2.2, eq. (2.11)]{DCFV2}, we know that for all $t \in \R$ and $x \in \Omega$
\begin{equation*}
\qquad
u^-(x)>t\,\,\Longrightarrow\,\,\chi_{\{u > t \}}^*(x)=1,
\qquad
u^+(x)\le t\,\,\Longrightarrow\,\,\chi_{\{u > t \}}^*(x)=0\,,
\end{equation*}
which easily implies
\begin{equation} \label{proof:Fusco}
\qquad
\chi_{\{u^->t\}}(x) = 1 \,\,\Longrightarrow\,\,\chi_{\{u > t \}}^-(x)=1,
\qquad
\chi_{\{u^+>t\}}(x) = 0 \,\,\Longrightarrow\,\,\chi_{\{u > t \}}^+(x)=0\,.
\end{equation}
If instead $u^-(x)\leq t < u^+(x)$ and $x\not\in N_t$, then necessarily
$x\in\{u > t \}^{1/2}$, which means $\chi_{\{u > t \}}^*(x)=1/2$. Hence, we get
$$\chi_{\{u > t \}}^+(x) + \chi_{\{u > t \}}^-(x) = 1, $$
and, since $\chi_{\{u > t \}}^\pm(x) \in \{0,1\}$, we obtain
$$\chi_{\{u > t \}}^+(x) = 1 \quad \text{ and } \quad \chi_{\{u > t \}}^-(x) = 0,$$
which means that
\begin{equation} \label{proof:Fusco_2}
\qquad
\chi_{\{u^+>t\}}(x) = 1 \,\,\Longrightarrow\,\,\chi_{\{u > t \}}^+(x)=1,
\qquad
 \chi_{\{u^->t\}}(x) = 0 \,\,\Longrightarrow\,\,\chi_{\{u > t \}}^-(x)=0\,.
\end{equation}
On the other hand, if $\chi_{\{u > t \}}^+(x)=1$, then we must have $u^+(x) > t$, because otherwise \eqref{proof:Fusco} would imply $\chi_{\{u > t \}}^+(x)=0$; and analogously, if $\chi_{\{u > t \}}^-(x)=0$, then we must have $u^-(x) \le t$, because otherwise \eqref{proof:Fusco} would imply $\chi_{\{u > t \}}^-(x)=1$. Hence, by noticing that $u^-(x)>t$ implies $u^+(x)>t$ and 
\begin{equation*}
1 = \chi_{\{u^->t\}}(x)  \le \chi_{\{u^+ > t \}}(x) \le 1\,,
\end{equation*}
we conclude from \eqref{proof:Fusco} and \eqref{proof:Fusco_2} that
\begin{equation*}
\chi_{\{u>t\}}^+(x) = 1 \,\,\Longleftrightarrow\,\,\chi_{\{u^+ > t \}}(x)=1.
\end{equation*}
Arguing analogously with the condition $u^+(x) \le t$, we deduce that
\begin{equation*}
\chi_{\{u>t\}}^-(x) = 0 \,\,\Longleftrightarrow\,\,\chi_{\{u^- > t \}}(x)=0\,.
\end{equation*}
Let now $x \in \Omega\setminus N_t$ such that $\chi_{\{u > t \}}^-(x)=1$. We notice that, for $\tau \in \R$, we have
\begin{equation*}
\left \{ \chi_{\{u > t \}} < \tau \right \} = 
\begin{cases} \Omega & \text{ if } \tau > 1 \\
 \{ u \le t\} & \text{ if } 0 < \tau \le 1 \\
\emptyset & \text{ if } \tau \le 0
\end{cases}.
\end{equation*}
Hence, thanks to the definitions \eqref{def:traces_Maz}, $\chi_{\{u > t \}}^-(x)=1$ implies 
\begin{equation*}
D\bigl(\{u\le t\};x\bigr)=0, \text{ and so } D\bigl(\{u > t \};x\bigr)=1.
\end{equation*}
All in all, this yields $u^-(x)\geq t$. However, if $u^-(x) = t$, then $x$ should satisfy $u^-(x)\leq t < u^+(x)$, which, given that $x \notin N_t$, would imply $x \in \{u > t \}^{1/2}$, but this contradicts the fact that $D\bigl(\{u > t \};x\bigr)=1$. Therefore, we conclude that
$u^-(x)>t$, and so we obtain 
$$\chi_{\{u>t\}}^-(x) = 1 \,\,\Longrightarrow\,\,\chi_{\{u^- > t \}}(x)=1\,.$$
In a similar way it is possible to prove that 
$$\chi_{\{u>t\}}^+(x) = 0 \,\,\Longrightarrow\,\,\chi_{\{u^+ > t \}}(x)=0\,,$$
so that, by \eqref{proof:Fusco}, we conclude that
\begin{equation*} 
\qquad
\chi_{\{u^->t\}}(x) = 1 \,\,\Longleftrightarrow\,\,\chi_{\{u > t \}}^-(x)=1,
\qquad
\chi_{\{u^+>t\}}(x) = 0 \,\,\Longleftrightarrow\,\,\chi_{\{u > t \}}^+(x)=0\,.
\end{equation*}
Finally, we notice that, if $u \in \Borel$, then $u^\pm$ do not depend on the choice of a representative from the equivalence class, and so the argument above works as well.
\end{proof}

\subsection{Divergence-measure fields }
\label{ss:div}

Given a measure valued vector field \(\A\in \mathcal M(\Omega; \R^N)\), we say that $\A$ is a {\em divergence-measure field} if its divergence in the sense of distributions is a finite Radon measure in 
\(\Omega\), acting as
\begin{equation}\label{f:defdiv}
\int_\Omega \varphi\, d\Div\A = -\int_{\Omega} \nabla\varphi\cdot d\A \qquad
\forall \varphi\in\cinftio.
\end{equation}
We denote by \(\DM(\Omega)
\) the space of all such vector fields.
Analogously, we define the local spaces $\DM_{\rm loc}(\Omega)$, as the sets of all vector fields $\A \in \mathcal{M}_{\rm loc}(\Omega; \R^N)$, such that $\div \A \in \mathcal{M}_{\rm loc}(\Omega)$. 

We exhibit here a family of $\DM$-fields that are not absolutely continuous with respect to $\Leb{N}$ and whose divergence is a non-trivial measure.
We fix $y \in \Omega$ and we consider
\begin{equation*}
\A :=(a_1, a_2,\dots, a_{N-1},a_{N}\Leb{N})\,,
\end{equation*}
where 
$$a_j = \Haus{1} \res \{ x \in \Omega : x_k = y_k \text{ for all } k \in \{1, 2, \dots, N \}, k \neq j \} \ \text{ for } j \in \{1, \dots, N - 1\}$$
and $a_N\in BV_{\rm loc}(\Omega)$. By construction, $\A\in \mathcal M_{\rm loc}(\Omega;\R^N)$, $a_j \perp \Leb{N}$ for every $j \in \{1, \dots, N-1\}$ and $\div\A=D_{x_N}a_N\in \mathcal M_{\rm loc}
(\Omega)$. 

Given $p \in [1, +\infty]$, by \(\DM^p
(\Omega)\) we denote the space of all
vector fields 
 \(\A\in L^p
(\Omega; \R^N)\)
whose divergence in the sense of distributions is a finite Radon measure in 
\(\Omega\), again acting as in \eqref{f:defdiv}. Similarly, we define the local spaces $\DM^p_{\rm loc}(\Omega)$, as the sets of all vector fields $\A \in L^p_{\rm loc}(\Omega; \R^N)$, such that $\div \A \in \mathcal{M}_{\rm loc}(\Omega)$. With a little abuse of notation, even in the case $\A \in L^p_{\rm loc}(\Omega; \R^N)$ for some $p \in [1, + \infty]$, we shall sometimes write $\A$ to denote the measure $\A \Leb{N}$.

For a more detailed exposition on the properties of these vector fields, we refer the reader to \cite{Anz,ChenFrid,ChFr1,ComiPayne,CCDM,CD3,crasta2019extension,CD5,CDM,Silh,MR2532602,Silhavy19}.

We recall that the divergence measure of a field $\A \in \DM^p_{\rm loc}(\Omega)$ enjoys absolute continuity properties with respect to suitable Hausdorff measures, depending on the value of $p \in [1, + \infty]$. More precisely, we have the following cases:
\begin{enumerate}
\item if $p=+\infty$, then $\div \A \ll \Haus{N-1}$ (\cite[Proposition 3.1]{ChenFrid} and \cite{Silh}*{Theorem 3.2}) and, if $\A \in L^\infty(\Omega; \R^N)$, then $|\Div \A| \le c_N \|\A\|_{L^\infty(\Omega; \R^N)} \Haus{N-1}$, where $c_N > 0$ is a constant depending only on the space dimension (\cite[Proposition 3.1]{Silhavy19});
\item if $p \in \left [\frac{N}{N - 1}, + \infty \right )$, then $|\Div \A|(B) = 0$ for all Borel sets $B$ of $\sigma$-finite $\Haus{N - \frac{p}{p - 1}}$ measure (\cite{Silh}*{Theorem 3.2});
\item if $p \in \left [1, \frac{N}{N - 1} \right )$, then the singularities of $\Div\A$ may be arbitrary (\cite{Silh}*{Example 3.3}).
\end{enumerate}
Actually, in case (2) a slightly stronger result holds; that is, $|\div \A|$ vanishes on Borel sets with zero $\frac{p}{p - 1}$-Sobolev capacity (\cite[Theorem 2.8]{MR3676052}), but this goes beyond the scope of our paper.

For the purposes of calculus, many different Leibniz rules for $\DM$ vector fields and suitably regular scalar functions have been discovered in the past years (see \cite{ChenFrid, ChFr1, Comi, Frid, MR2532602}). We collect here a list of such rules. 
In the most general case of $\A \in \DM_{\rm loc}(\Omega)$, if $u \in C^1_c(\Omega)$, then we have $u\A\in\DM_{\rm loc}(\Omega)$ and
\begin{equation}\label{eq:leibnitzC1}
\Div(u\A) = u\Div \A + \nabla u\cdot \A \ \text{ on } \Omega,
\end{equation}
which is a particular case of \cite[Theorem 3.2]{ChFr1}. 
More in general, even if $u \in \Lip_{\rm loc}(\Omega)$, it is possible to prove that $u\A\in\DM_{\rm loc}(\Omega)$, with
\begin{equation}
\Div(u\A) = u\Div \A + \left\langle\left\langle\nabla u,\A\right\rangle\right\rangle \ \text{ on } \Omega,
\label{eq:pairingShilavy}
\end{equation}
where $\left\langle\left\langle\nabla u,\A\right\rangle\right\rangle$ is the weak$^*$-limit of the family of measures $\nabla u_\rho \cdot \A$, for any standard mollification $u_\rho$ of $u$, which is a Radon measure satisfying 
\begin{equation} \label{eq:abs_cont_pairing_irreg}
\left | \left\langle\left\langle\nabla u,\A\right\rangle\right\rangle \right | \le \|\nabla u\|_{L^{\infty}(\Omega'; \R^N)} |\A| \ \text{ on } \Omega'
\end{equation}
for every open set $\Omega' \Subset \Omega$ (see \cite[Proposition 2.2]{MR2532602}).

If instead $\A \in \DM^p_{\rm loc}(\Omega)$ for some $p \in [1, +\infty]$, then, under the following set of assumptions:
\begin{enumerate}
\item if $p\in[1,+\infty)$ and $u\in L^\infty_{\rm loc}(\Omega)\cap W^{1,q}_{\rm loc}(\Omega)$, where $q = \frac{p}{p-1}$ is the conjugate exponent of $p$,
\item if $p = + \infty$ and $u\in L^\infty_{\rm loc}(\Omega)\cap BV_{\rm loc}(\Omega)$,
\end{enumerate}
it holds that $u\A\in\DM^r_{\rm loc}(\Omega)$ for all $r\in[1,p]$, with
\begin{equation} \label{eq:Leibniz_Sobolev_classic}
\Div(u\A) = u^*\Div \A +
\begin{cases}
 \A\cdot\nabla u\,\Leb{N}, & \mbox{if $p< + \infty$, or $p =+ \infty$ and $u\in L^\infty_{\rm loc}(\Omega)\cap W^{1,1}_{\rm loc}(\Omega)$,} \\
\\
 (\A, Du), & \mbox{if $p=+\infty$ and $u\in L^\infty_{\rm loc}(\Omega)\cap (BV_{\rm loc}(\Omega)\backslash W^{1,1}_{\rm loc}(\Omega))$,}
\end{cases} 
\end{equation}
where $u^*$ is the precise representative of $u$, which satisfies $u^*(x) = \tilde{u}(x)$ for $|\div \A|$-a.e. $x \in \Omega$ if $u$ is a Sobolev function, and $(\A, Du)$ is the standard {\em pairing measure} between $\A$ and $Du$ introduced in \cite{Anz} (see \cite{CD3} for more details). For such results, we refer the reader to \cite[Theorem 3.1]{ChenFrid}, \cite[Theorem 3.2.3]{Comi}, \cite[Theorem 2.1]{Frid}. The pairing measure can be also characterized as the weak$^*$-limit of the measures $\A \cdot \nabla (u * \rho_\eps) \Leb{N}$ as $\eps \to 0$, for any standard mollifier $\rho$. In addition, for every open set $\Omega' \Subset \Omega$ we have
\begin{equation} \label{eq:abs_cont_classic_pairing}
|(\A, Du)| \le \|\A\|_{L^{\infty}(\Omega'; \R^N)} |D u| \ \text{ on } \Omega'
\end{equation}
by \cite[Proposition 3.4]{CCDM}.

Finally, in \cite[Proposition 4.4]{CDM} the notion of $\lambda$-pairings is introduced in the case $\A \in \DM^\infty_{\rm loc}(\Omega)$. More precisely, given $u\in BV_{\rm loc}(\Omega)$ such that $u^* \in L^1(\Omega; |\div \A|)$ and a Borel function $\lambda:\Omega\to[0,1]$, then \eqref{def_lambdapairing} defines a family of $\lambda$-pairings, which provide the following Leibniz rule
\begin{equation}\label{f:geninfty}
\Div(u \A) = {u}^\lambda\, \Div\A + \pair{\A, Du}_\lambda,
\end{equation}
where $u^\lambda$ is given by \eqref{f:pr} (see the subsequent discussion for the particular case of $BV$ functions). In addition, analogously to \eqref{eq:abs_cont_classic_pairing}, for every open set $\Omega' \Subset \Omega$ we have
\begin{equation} \label{eq:abs_cont_classic_pairing_lambda}
|(\A, Du)_\lambda| \le \|\A\|_{L^{\infty}(\Omega'; \R^N)} |D u| \ \text{ on } \Omega'
\end{equation}
If $\lambda \equiv \frac{1}{2}$, then we retrieve the standard pairing given by \eqref{def_pairing} (see also \cite[Theorem 4.12]{CD3}). If instead $\lambda \equiv 0$ or $\lambda \equiv 1$, we simply write $(\A, Du)_0$ and $(\A, Du)_1$, respectively.

We state now a technical result which can be seen as a basic version of a Gauss--Green formula. It has been proved already in some special cases (see for instance \cite[Lemma 3.1]{ComiPayne}), while we were not able to find it in literature in the most general form.

\begin{lemma} \label{lem:comp_supp_div_A_0}
Let $\A \in \DM(\Omega)$ be such that ${\rm supp}(|\A|) \Subset \Omega$. Then $\div \A(\Omega) = 0$.
\end{lemma}
\begin{proof}
Let $V \Subset \Omega$ be an open set such that ${\rm supp}(|\A|) \subset V$. It is well known (see, e.g., \cite[Remark 2.21]{ComiPayne}) that ${\rm supp}(\div \A) \subset {\rm supp}(|\A|)$. Therefore, $\div \A = 0$ on $\Omega \setminus \overline{V}$. Then, we can apply the definition of weak divergence \eqref{f:defdiv} to a function $\eta \in C^{\infty}_c(\Omega)$ such that $\eta \equiv 1$ on a neighborhood of $V$: we get
\begin{equation*}
 \div \A (\overline{V}) = \int_{\overline{V}} \eta \, d \div \A = \int_{\Omega} \eta \, d \div \A = - \int_{\Omega} \nabla \eta \cdot \, d \A = - \int_{\Omega \setminus V } \nabla \eta \cdot \, d \A = 0.
\end{equation*}
Thus, we get $\div\A(\Omega) =  \div \A (\overline{V}) +  \div \A (\Omega \setminus \overline{V}) = 0$.
\end{proof}

We conclude this section with the Gauss--Green formulas for essentially bounded divergence-measure fields and sets of finite perimeter, for which we refer the reader to \cite[Theorem 3.2]{ComiPayne}. To this purpose, we recall that a measurable $E \subset \Omega$ is a {\em set of finite perimeter} if $\chi_E \in BV(\Omega)$; and, in such a case, we have $|D\chi_E| = \Haus{N-1} \res \partial^* E$, due to De Giorgi's and Federer's Theorems (see for instance \cite[Theorems 3.59 and 3.61]{AFP}). Then, given $\A \in \DM^{\infty}(\Omega)$ and $E \Subset \Omega$ of finite perimeter, we have
\begin{equation} \label{eq:classical_Gauss_Green}
\div \A(E^1) = 
- \int_{\partial^* E} {\rm Tr}^i (\A,\partial^*E) \, d \Haus{N-1}
 \ \text{ and } \
 \div \A(E^1 \cup \partial^* E) = - \int_{\partial^* E} {\rm Tr}^e (\A,\partial^*E) \, d \Haus{N-1},
\end{equation}
where ${\rm Tr}^i (\A,\partial^*E), {\rm Tr}^e (\A,\partial^*E) \in L^\infty(\partial^*E, \Haus{N-1})$ are the {\em interior and exterior normal traces} of $\A$ on $\partial^* E$. In particular, the normal traces can be characterized as the densities of the measures $(\A, D \chi_E)_0$ and $(\A, D \chi_E)_1$, respectively, with respect to the measure $|D\chi_E|$, thanks to \cite[Proposition 4.7]{CDM} applied to the case $u = \chi_E$. More precisely, if $\A \in \DM^\infty_{\rm loc}(\Omega)$ and $E \subset \Omega$ is a set of locally finite perimeter, we have
\begin{equation} \label{eq:def_normal_traces}
(\A, D \chi_E)_0 =  {\rm Tr}^i (\A,\partial^*E) \, \Haus{N-1} \res \partial^* E \ \text{ and } \ (\A, D \chi_E)_1 = {\rm Tr}^e (\A,\partial^*E) \, \Haus{N-1} \res \partial^* E
\end{equation}
and, for all open sets $\Omega' \Subset \Omega$,
\begin{equation} \label{eq:normal_traces_L_infty}
\begin{split}
\|{\rm Tr}^i (\A,\partial^*E)\|_{L^\infty(\Omega' \cap \partial^* E, \Haus{N-1})} & \le \|\A\|_{L^\infty(\Omega' \cap E; \R^N)}, \\
\|{\rm Tr}^e (\A,\partial^*E)\|_{L^\infty(\Omega' \cap \partial^* E, \Haus{N-1})} & \le \|\A\|_{L^\infty(\Omega' \setminus E; \R^N)}.
\end{split}
\end{equation}
For these results we refer the reader to \cite[Proposition 4.7]{CDM} and \cite[Theorem 4.2]{ComiPayne}.

\section{Pairing-related $BV$-type functions}
\label{s:spaceBVA}

In this section we introduce our new general definition of $\lambda$-pairing.

\subsection{Ambient class and general pairing}
First of all, we define good ambient classes of summable functions, inevitably depending on the chosen field and the Borel function which defines the representatives. 

Given a vector field $\A\in \DM_{\rm loc}(\Omega)$ and a Borel function $\lambda\colon \Omega \to [0,1]$, we set
\begin{equation*}
X^{\A,\lambda}(\Omega):=\{u \in \Borel : u^{\lambda} \in L^1(\Omega, |\A|) \cap L^1(\Omega, |\Div\A|)
\},
\end{equation*}
\begin{equation*}
X^{\A,\lambda}_{\rm loc}(\Omega):=\{u\in \Borel : u^{\lambda} \in L^1_{\rm loc}(\Omega, |\A|) \cap L^1_{\rm loc}(\Omega, |\Div\A|)
\}.
\end{equation*}
It is interesting to notice that these sets of functions are not linear spaces, in general, due to the fact that the $\lambda$-representative of a sum is not the sum of $\lambda$-representatives, see Remark \ref{rem:not_lin_pairing} for some examples.

We define now a $\lambda$-pairing for functions in $X^{\A, \lambda}_{\rm loc}(\Omega)$.

\begin{definition}[General $\lambda$-pairing]\label{d:lambdapairbis}
Let $\A\in \DM_{\rm loc}(\Omega)$, $\lambda\colon \Omega \to [0,1]$ be a Borel function
and $u\in X^{\A,\lambda}_{\rm loc}(\Omega)$. We define the $\lambda$-pairing between $\A$ and $u$ as the distribution
$$\pair{\A,Du}_\lambda\colon C^\infty_c(\Omega) \to \R$$ 
acting as
\begin{equation} \label{eq:def_lambda_pairing}
\pscal{\pair{\A,Du}_\lambda}{\varphi}:= -\int_\Omega {u^\lambda\,\varphi} \, d \div \A - \int_\Omega u^\lambda \nabla \varphi \cdot d \A
\quad \text{ for } \
\varphi\in C^\infty_c(\Omega)\,.
\end{equation}
\end{definition}

We emphasize that, for a given function $u \in X^{\A,\lambda}_{\rm loc}(\Omega)$, the derivative $Du$ may not exist as a Radon measure, but, if in addition $u \in L^1_{\rm loc}(\Omega)$, $Du$ is well defined as a distribution (of order one), and so, with a little abuse of notation, we choose the standard $\lambda$-pairing notation as in the case of $BV$ functions, \cite[Definition 4.1]{CDM}.

In addition, if $\A \in \DM^{1}_{\rm loc}(\Omega)$, then \eqref{eq:def_lambda_pairing} becomes
\begin{equation} \label{eq:def_lambda_pairing_abs_cont}
\pscal{\pair{\A,Du}_\lambda}{\varphi}:= -\int_\Omega {u^\lambda\,\varphi} \, d \div \A - \int_\Omega u \nabla \varphi \cdot \A \, dx
\quad \text{ for } \
\varphi\in C^\infty_c(\Omega)\,,
\end{equation}
since $u^\lambda(x) = \tilde{u}(x) = u(x)$ for $\Leb{N}$-a.e. $x \in \Omega$, where we denote by $u$ one of the elements of its equivalence class, with a little abuse of notation, following the convention adopted in Section \ref{ss:def}.

We provide an alternative equivalent formulation of the definition of $\lambda$-pairing \eqref{eq:def_lambda_pairing}, which is often used in the following.

\begin{lemma} \label{eq:further_def_pairing}
Let $\A\in \DM_{\rm loc}(\Omega)$, $\lambda\colon \Omega \to [0,1]$ be a Borel function and $u\in X^{\A,\lambda}_{\rm loc}(\Omega)$. Then, the distribution $\pair{\A,Du}_\lambda$ is of order 1 and satisfies
\begin{equation}\label{eq:further}
\pscal{\pair{\A,Du}_\lambda}{\varphi}= -\int_\Omega u^\lambda \, d \div (\varphi \A)
\quad \text{ for all } \
\varphi\in C^1_c(\Omega)\,.
\end{equation}
In particular, if $u \equiv c$ for some $c \in \R$, then $\pair{\A,D c}_\lambda = 0$.
\end{lemma}

\begin{proof}
It is clear that 
\begin{equation*}
\left | \pscal{\pair{\A,Du}_\lambda}{\varphi} \right | \le \|\varphi\|_{L^\infty(\Omega)} \int_\Omega |u^\lambda| \, d |\div \A| + \|\nabla \varphi \|_{L^\infty(\Omega; \R^N)} \int_\Omega |u^\lambda| d |\A|,
\end{equation*}
so that the distribution $\pair{\A,Du}_\lambda$ can be extended to $C^1_c$ test functions. Then, \eqref{eq:further} is a straightforward consequence of \eqref{eq:leibnitzC1}. Indeed, for all $\varphi \in C^1_c(\Omega)$, we can rewrite \eqref{eq:def_lambda_pairing} as
\begin{equation*} 
\pscal{\pair{\A,Du}_\lambda}{\varphi}:= -\int_\Omega u^\lambda \left ( \varphi \, d \div \A + \nabla \varphi \cdot d \A \right ),
\end{equation*}
and now we apply the Leibniz rule \eqref{eq:leibnitzC1} to $\A$ and $\varphi$. Finally, if $u$ is constant and equal to $c \in \R$, then we have
\begin{equation*} 
\pscal{\pair{\A,D c}_\lambda}{\varphi}:= -\int_\Omega c \, d \div (\varphi \A) = - c \div (\varphi \A)(\Omega) = 0,
\end{equation*}
thanks to Lemma \ref{lem:comp_supp_div_A_0}.
\end{proof}

\subsection{The class $BV^{\A, \lambda}(\Omega)$} \label{subsec:BV_A_lambda}

Arguing in analogy with the classical theory of functions of bounded variation, we introduce a version of $BV$-type classes with respect to the variation given by this generalized $\lambda$-pairing.

\begin{definition}
Given $\A \in \DM_{\rm loc}(\Omega)$ and a Borel function $\lambda\colon \Omega \to [0,1]$, we define the classes
\begin{align*}
BV^{\A,\lambda}(\Omega) & := \left\{
u\in X^{\A,\lambda}(\Omega): (\A,Du)_{\lambda}\in \radon
\right\}\,,
\\
BV^{\A,\lambda}_{\rm{loc}}(\Omega) & := \left\{
u\in X^{\A,\lambda}_{\rm loc}(\Omega): (\A,Du)_{\lambda}\in {\mathcal{M}_{\rm loc}(\Omega)}
\right\}\,.
\end{align*}
\end{definition}

\begin{remark} \label{rem:not_lin_pairing_short}
We point out that $BV^{\A, \lambda}(\Omega)$ is not a linear space, in general, due to the fact that the $\lambda$-pairing is not linear in the second component. This was already noticed in \cite[Remark~4.6]{CDM} in the classical case of $\A \in \DM^\infty(\Omega)$ and $u \in BV(\Omega)$ with $u^* \in L^1(\Omega, |\div \A|)$, whenever $\lambda(x) \neq \frac{1}{2}$ for all $x \in B$, for some Borel set $B \subseteq J_u$ with $|\div \A|(B) > 0$. In our setting, however, the map $u\to (\A,Du)_\lambda$ may fail to be linear even in the case $\lambda \equiv \frac{1}{2}$, as showed by the second example in Remark \ref{rem:not_lin_pairing}. It is also relevant to notice that the linearity may be ensured in some particular cases, which are explored in Section \ref{s:banach}.
\end{remark}

Whenever $\lambda$ is constant; that is, $\lambda \equiv t$ for some $t \in [0, 1]$, we simply write $BV^{\A, t}(\Omega)$ (and analogously for the local classes). In the extreme cases $\lambda \equiv 0$ and $\lambda \equiv 1$, we therefore have $BV^{\A,0}(\Omega)$ and $BV^{\A,1}(\Omega)$, respectively.
In the particular case $\lambda \equiv \frac{1}{2}$, we use the following shorthand notation:
\begin{equation} \label{eqdef:BV_A}
X^{\A}(\Omega) := X^{\A,\frac{1}{2}}(\Omega), \quad BV^{\A}(\Omega) := BV^{\A, \frac{1}{2}}(\Omega), \quad (\A,Du) := (\A,Du)_{\frac{1}{2}},
\end{equation}
and analogously for the local classes. This case is indeed relevant due to the fact that, in the classical theory of pairings, for $\A \in \DM^{\infty}_{\rm loc}(\Omega)$ and $u \in BV_{\rm loc}(\Omega)$, we have $u^{\frac{1}{2}}(x) = u^*(x)$ for $|\div \A|$-a.e. $x \in \Omega$, and, as long as $u^* \in L^1_{\rm loc}(\Omega, |\div \A|)$, the pairing which we obtain coincides with the one defined by Anzellotti. We discuss this equivalence more in detail in the following proposition.

\begin{proposition} \label{prop:main_inclusions}
Let $\A \in \DM_{\rm loc}(\Omega)$ and $\lambda : \Omega \to [0, 1]$ be a Borel function. 
\begin{enumerate}
\item For all $u \in X^{\A, \lambda}_{\rm loc}(\Omega)$ we have
\begin{equation} \label{eq:pairing_div_lambda_distrib}
\pscal{\pair{\A,Du}_\lambda}{\varphi} = -\int_\Omega \varphi \, u^\lambda \, d \div \A + \pscal{\div(u^\lambda \A)}{\varphi} \ \text{ for } \ \varphi \in C^{1}_c(\Omega)
\end{equation}
in the sense of distributions.
Hence, $u\in BV^{\A,\lambda}(\Omega)$ if and only if $u \in X^{\A, \lambda}(\Omega)$ and $\div(u^\lambda \A) \in \mathcal{M}(\Omega)$, and, for all $u\in BV^{\A,\lambda}_{\rm loc}(\Omega)$, we get
\begin{equation}\label{f:senseofmeasures2}
(\A,Du)_\lambda=-u^\lambda\,\Div\A+\Div(u^\lambda \A) \ \text{ on } \ \Omega,
\end{equation}
in the sense of Radon measures. Therefore, 
\begin{align*}
BV^{\A,\lambda}(\Omega) & = \left\{
u\in X^{\A,\lambda}(\Omega): \Div(u^\lambda \A)\in \radon
\right\}\,, \\
BV^{\A,\lambda}_{\rm loc}(\Omega) & = \left\{
u\in X^{\A,\lambda}_{\rm loc}(\Omega): \Div(u^\lambda \A)\in \mathcal{M}_{\rm loc}(\Omega)
\right\}\,.
\end{align*}

\item If $\A \in \DM^1_{\rm loc}(\Omega)$, then for all $u \in X^{\A, \lambda}_{\rm loc}(\Omega)$ we have $\Div(u^\lambda \A) = \Div(u \A)$ in the sense of distributions. Hence, $u \in BV^{\A, \lambda}(\Omega)$ if and only if $u \in X^{\A, \lambda}(\Omega)$ and $\Div(u \A) \in \mathcal{M}(\Omega)$, and, for all $u\in BV^{\A,\lambda}_{\rm loc}(\Omega)$, we get
\begin{equation}\label{eq:Leibniz_A_vector_field}
(\A,Du)_\lambda=-u^\lambda\,\Div\A+\Div(u \A) \ \text{ on } \ \Omega,
\end{equation}
in the sense of Radon measure. Therefore, 
\begin{align*}
BV^{\A,\lambda}(\Omega) & = \left\{
u\in X^{\A,\lambda}(\Omega): \Div(u \A)\in \radon
\right\}\,, \\
BV^{\A,\lambda}_{\rm loc}(\Omega) & = \left\{
u\in X^{\A,\lambda}_{\rm loc}(\Omega): \Div(u \A)\in \mathcal{M}_{\rm loc}(\Omega)
\right\}\,.
\end{align*}

\item If $\A \in \DM^{\infty}(\Omega)$, then we have
$$ 
BV(\Omega)\cap BV^{\A,\lambda}(\Omega) = BV(\Omega) \cap X^{\A, \lambda}(\Omega) = \{ u \in BV(\Omega) : u^* \in L^1(\Omega, |\div \A|) \},
$$
and, if $u \in BV(\Omega)$ and $u^* \in L^1(\Omega, |\div \A|)$, then the $\lambda$-pairing $(\A, Du)_\lambda$ coincides with the one defined in \cite[Definition 4.1]{CDM}.

\item If $\A \in \DM^1(\Omega)$, then for every two Borel functions $\lambda_1, \lambda_2 \colon\Omega\to [0,1]$, we have 
$$
BV^{\A, \lambda_1}(\Omega) \cap L^\infty(\Omega) = BV^{\A, \lambda_2}(\Omega) \cap L^\infty(\Omega).
$$
In particular, if $u \in BV^{\A, \lambda_1}(\Omega) \cap L^\infty(\Omega)$, then
$$(\A, Du)_{\lambda_1} - (\A, Du)_{\lambda_2} = (u^{\lambda_2} - u^{\lambda_1}) \div \A \res S_u^* = (\lambda_2 - \lambda_1) (u^+ - u^-)  \div \A \res S_u^* .$$
\end{enumerate}
\end{proposition}

\begin{proof}
In order to prove (1), we observe that, given $\A \in \mathcal{M}_{\rm loc}(\Omega; \R^N)$ and $v \in L^1_{\rm loc}(\Omega, |\A|)$, we clearly have
\begin{equation*}
\pscal{\div(v \A)}{\varphi} = - \int_\Omega v \nabla \varphi \cdot d \A \ \text{ for } \ \varphi \in C^{\infty}_c(\Omega)
\end{equation*}
in the sense of distributions. Hence, for $u \in X^{\A, \lambda}_{\rm loc}(\Omega)$ \eqref{eq:def_lambda_pairing} can be rewritten as \eqref{eq:pairing_div_lambda_distrib}. Given that both distributions are clearly of order 1, we can test them against functions in $C^1_c(\Omega)$. Therefore, \eqref{f:senseofmeasures2} follows as long as either $(\A,Du)_\lambda$ or $\Div(u^\lambda \A)$ are Radon measures. This concludes the proof of point (1).

Point (2) follows immediately from point (1) by noticing that, if $\A \in L^1_{\rm loc}(\Omega; \R^N)$, then $u^\lambda(x) = u(x)$ for $|\A| \Leb{N}$-a.e. $x \in \Omega$. 

As for point (3), we start by observing that Definition \ref{d:lambdapairbis} extends the standard definition of $\lambda$-pairing \cite[Definition 4.1]{CDM}, given that the equation satisfied by $(\A, Du)_\lambda$ is the same in both cases.
Now, if $u \in BV(\Omega)$, in particular $u$ is a class of Lebesgue measurable functions, and therefore, as an easy consequence of Lusin's Theorem, it admits a Borel representative. Therefore, we see that $u \in \Borel$. Then, we recall that, if $u\in BV(\Omega)$, we have $\mathcal H^{N-1}(S_u \setminus J_u) = 0$ (for instance, see \cite[Theorem 3.78]{AFP}) and so $|\div \A|(S_u \setminus J_u) = 0$, since $\div \A \ll \Haus{N-1}$ for $\A \in \DM^{\infty}_{\rm loc}(\Omega)$, see Section \ref{ss:div}. This means that $u^*(x) = u^{\frac{1}{2}}(x)$ for $|\div \A|$-a.e. $x \in \Omega$ (see Section \ref{ss:def}).
Moreover, as proven in \cite[Lemma 3.2]{CDM} for every $u\in BV(\Omega)$ and 
for every two Borel functions $\lambda_1, \lambda_2 \colon\Omega\to [0,1]$,
\(
\prec[\lambda_1]{u}\in L^1_{\rm{loc}}(\Omega,|\Div \A|)
\)
if and only if
\(
\prec[\lambda_2]{u}\in L^1_{\rm{loc}}(\Omega,|\Div \A|).
\)
Hence, if $u \in BV(\Omega)\cap X^{\A,\lambda}(\Omega)$ for some Borel function $\lambda: \Omega \to [0, 1]$, then we get $u^{\frac{1}{2}} \in L^1(\Omega,|\Div \A|)$, which allows to conclude that $u^* \in L^1(\Omega,|\Div \A|)$. This proves the inclusions
\begin{equation*}
BV(\Omega)\cap BV^{\A,\lambda}(\Omega) \subseteq BV(\Omega)\cap X^{\A,\lambda}(\Omega) \subseteq \{ u \in BV(\Omega) : u^* \in L^1(\Omega, |\div \A|) \}.
\end{equation*}
As for the second inclusion, the opposite one can be proved by reversing this argument. Then, we notice that \cite[Lemma 3.2 and Proposition 4.4]{CDM} imply that, if $u \in BV(\Omega)$ is such that $u^* \in L^1(\Omega, |\div \A|)$, then $(\A, Du)_\lambda \in \mathcal{M}(\Omega)$: this proves the remaining opposite inclusion.

Finally, in dealing with point (4) we start by observing that, if $u\in L^\infty(\Omega)$, we have 
\begin{equation*}
u^{\lambda_1},u^{\lambda_2} \in L^1(\Omega, |\A| \Leb{N}) \cap L^1(\Omega, |\Div\A|).
\end{equation*} 
Hence, $u \in X^{\A, \lambda_1}(\Omega) \cap X^{\A, \lambda_2}(\Omega)$. Then, we notice that 
\begin{equation*}
u^{\lambda_1} \A = u^{\lambda_2} \A = u \A \text{ up to a }\Leb{N}\text{-negligible set},
\end{equation*} 
so that, if $u \in BV^{\A, \lambda_1}(\Omega) \cap L^\infty(\Omega)$, then point (2) implies $\div(u\A) \in \mathcal{M}(\Omega)$, and thus $u \in BV^{\A, \lambda_2}(\Omega)$. The opposite inclusion follows by exchanging the roles of $\lambda_1$ and $\lambda_2$. Then, if $u \in BV^{\A, \lambda_1}(\Omega) \cap L^\infty(\Omega)$, we know that it also belongs to $BV^{\A, \lambda_2}(\Omega) \cap L^\infty(\Omega)$, so that we obtain \eqref{eq:Leibniz_A_vector_field} in both cases, and we subtract the equation with $\lambda = \lambda_2$ from the one with $\lambda = \lambda_1$. Thus, it is enough to exploit \eqref{eq:lambda_u_tilde_S_u_*} and the fact that $Z_u = \emptyset$ to end the proof.
\end{proof}

\begin{remark}
We point out that the \quotemarks{if and only if} condition in point (1) of Proposition \ref{prop:main_inclusions} is actually equivalent to saying that, given $\A \in \DM_{\rm loc}(\Omega)$, $\lambda : \Omega \to [0, 1]$ Borel and $u \in X^{\A, \lambda}_{\rm loc}(\Omega)$, we have
\begin{equation*}
u \in BV^{\A, \lambda}_{\rm loc}(\Omega) \iff u^\lambda \A \in \DM_{\rm loc}(\Omega).
\end{equation*}
Analogously, the necessary and sufficient condition in point (2), for which we have $\A \in\DM^1_{\rm loc}(\Omega)$, can be rewritten as
\begin{equation*}
u \in BV^{\A, \lambda}_{\rm loc}(\Omega) \iff u \A \in \DM^1_{\rm loc}(\Omega).
\end{equation*}
We stress the fact that no equivalent statement holds for higher summability of the field $\A$, since $u \in X^{\A, \lambda}_{\rm loc}(\Omega)$ only implies that $u \A \in L^1_{\rm loc}(\Omega; \R^N)$.
\end{remark}

As a consequence of the Leibniz-type rule \eqref{f:senseofmeasures2}, we obtain an integration by parts formula for $BV^{\A, \lambda}$ functions with compact support.

\begin{lemma} \label{lem:gausslambda_u}
Let $\A \in \DM_{\rm loc}(\Omega)$, $\lambda : \Omega \to [0, 1]$ be a Borel function and $u\in BV^{\A, \lambda}(\Omega)$ be such that ${\rm supp}(u) \Subset \Omega$. Then we have 
\begin{equation}\label{eq:gausslambda2}
\int_{\Omega} u^\lambda\,d\div  \A= - \int_{\Omega}
 \, d (\A, Du)_\lambda \,.
\end{equation}
\end{lemma}

\begin{proof}
Since ${\rm supp}(u) \Subset \Omega$, it is clear that $\supp (u^\lambda |\A|) \Subset \Omega$, and therefore, since $u \in BV^{\A, \lambda}(\Omega)$, we conclude that $u^\lambda \A \in \DM(\Omega)$. Hence, thanks to Lemma \ref{lem:comp_supp_div_A_0}, we see that $\div(u^\lambda \A)(\Omega) = 0$. Thus, \eqref{eq:gausslambda2} follows by evaluating \eqref{f:senseofmeasures2} over $\Omega$. 
\end{proof}

\begin{remark} \label{rem:Silhavy_DM}
If $\A \in \DM_{\rm loc}(\Omega)$ and $u \in \Lip_{\rm loc}(\Omega)$, then $u \in BV^{\A, \lambda}_{\rm loc}(\Omega)$ for all Borel function $\lambda : \Omega \to [0, 1]$, and
\begin{equation*}
(\A, Du)_\lambda = \left\langle\left\langle\nabla u,\A\right\rangle\right\rangle \quad \mbox{ on } \, \Omega,
\end{equation*} 
where $\left\langle\left\langle\nabla u,\A\right\rangle\right\rangle$ was recalled in \eqref{eq:pairingShilavy}. Indeed, $u$ is continuous, so that 
\begin{equation*}
u^\lambda(x) = \tilde{u}(x) = u(x)  \ \text{ for all } x \in \Omega.
\end{equation*}
The assertion then follows combining \eqref{f:senseofmeasures2} and \eqref{eq:pairingShilavy}. Furthermore, if $\A \in \DM(\Omega)$ and $u \in \Lip(\Omega) \cap X^{\A}(\Omega)$, then $u \in BV^{\A, \lambda}(\Omega)$ for all Borel functions $\lambda : \Omega \to [0, 1]$.
\end{remark}

\begin{remark} \label{rem:A_dm_pure_no_equiv_L_infty}
Point (4) in Proposition \ref{prop:main_inclusions} does not hold if $\A \in \DM(\Omega) \setminus \DM^1_{\rm loc}(\Omega)$, in general. Indeed, let us consider $N \ge 2$, $\Omega = (-1, 1)^N$, $\A = \left ( \Haus{1} \res J , 0, \dots, 0 \right )$, where $J$ is the segment
$$J = \{ x \in \Omega : x_j = 0 \text{ for all } j = 2, \dots, N \} = \{ x \in \Omega : x = (t, 0, \dots, 0) \text{ for } t \in (-1, 1) \},$$
so that $\div \A = 0$ and $\A \in \DM(\Omega) \setminus \DM^1_{\rm loc}(\Omega)$. Then, we choose $u = \chi_{(-1, 1) \times (0, 1)^{N-1}}$: it is clear that $u \in L^\infty(\Omega)$ and
\begin{equation*}
u^-(x) = 0 \ \text{ and } \ u^+(x) = 1 \ \text{ for all } x \in J.
\end{equation*}
In particular, we see that $u \in BV^{\A, 0}(\Omega)$. Indeed, $u^0(x) = u^-(x) = 0$ for $|\A|$-a.e. $x \in \Omega$, so that $u^- \in L^1(\Omega, |\A|)$, obviously $u^- \in L^1(\Omega, |\div \A|)$, and, thanks to \eqref{f:senseofmeasures2},
\begin{equation*}
(\A, Du)_0 = -u^-\,\Div\A+\Div(u^- \A) = 0.
\end{equation*} 
Let now $\overline{\lambda}(x) = \chi_{F}(x_1)$, for some Borel set $F \subset (-1, 1)$ such that $\chi_F \notin BV_{\rm loc}((-1, 1))$ (for instance, $F$ could be any fat Cantor set). Then we have
\begin{equation*}
u^{\overline{\lambda}}(x) = (1 - \overline{\lambda}(x)) u^-(x) + \overline{\lambda}(x) u^+(x) = \chi_F(x_1) \ \text{ for all } x \in J.
\end{equation*}
Hence, $u \in X^{\A, \overline{\lambda}}(\Omega)$, and for all $\varphi \in C^1_c(\Omega)$ we get
\begin{align*}
\pscal{\pair{\A,Du}_{\overline{\lambda}}}{\varphi} & = - \int_\Omega \varphi u^{\overline{\lambda}} \, d \div \A - \int_\Omega u^{\overline{\lambda}} \nabla \varphi \cdot d \A = - \int_J \chi_F(x_1) \frac{\partial \varphi(x)}{\partial x_1} \, d \Haus{1} \\
& = - \int_{-1}^{-1} \chi_F(x_1) \frac{\partial \varphi(x)}{\partial x_1} \, d x_1.
\end{align*}
Therefore, $\pair{\A,Du}_{\overline{\lambda}} \notin \mathcal{M}_{\rm loc}(\Omega)$, given that $\chi_F \notin BV_{\rm loc}((-1, 1))$, and so $u \notin BV^{\A, \overline{\lambda}}_{\rm loc}(\Omega)$. Thus, we conclude that
\begin{equation*}
BV^{\A, 0}(\Omega) \cap L^\infty(\Omega) \neq BV^{\A, \overline{\lambda}}(\Omega) \cap L^\infty(\Omega).
\end{equation*}
\end{remark}

\subsection{Comparison with $BV$} An immediate consequence of point (3) Proposition \ref{prop:main_inclusions} is that, if $\A \in \DM^\infty(\Omega)$, then 
\begin{equation} \label{eq:BV_inclusion_BV_A_lambda}
\{ u \in BV(\Omega) : u^* \in L^1(\Omega, |\div \A|) \} \subseteq BV^{\A,\lambda}(\Omega)
\end{equation}
for all Borel functions $\lambda : \Omega \to [0, 1]$.
In general, the inclusion is strict, as shown in the following remark. 

\begin{remark} \label{eq:examplev}
Let $N \ge 2$ and $\lambda : \Omega \to [0, 1]$ be a Borel function.
There exist a field $\A \in \DM^\infty(\Omega)$ and a function $v\in BV^{\A, \lambda}(\Omega)$ such that $v\not\in BV_{\rm loc}(\Omega)$. Let us consider 
\begin{equation*}
\A(x):=(a_1(\hat{x}_1), a_2(\hat{x}_2),\dots, a_{N-1}(\hat{x}_{N-1}),0)\,,
\end{equation*}
where 
\begin{equation*}
\hat{x}_j:=(x_1,\dots,x_{j-1},x_{j+1},\dots,x_N), \ a_j\in L^\infty(\Omega_j) \ \text{ and } \ \Omega_j = \{ y \in \R^{N-1} : y = \hat{x}_j \text{ for some } x \in \Omega \}
\end{equation*}
for every $j\in\{1,2,\dots,N-1\}$. By construction, $\div\A=0$. Then, thanks to \eqref{eq:Leibniz_A_vector_field}, a function $u\in X^{\A,\lambda}(\Omega)$ belongs to $BV^{\A, \lambda}(\Omega)$ if and only if $(\A,Du)_\lambda =\div(u\A)\in\mathcal{M}(\Omega)$. 

Now, let $y \in \Omega$, $r > 0$ be such that $B_{2r}(y) \subset \Omega$ and $\eta_{r, y} \in C^\infty_c(B_{2r}(y))$ be such that $\eta_{r, y} \equiv 1$ on $B_r(y)$. Then we set
$$v(x):=\log(|x_N - y_N|) \eta_{r, y}(x).$$ 
It is clear that $v\in L^1(\Omega)$, which implies $v \in L^1(\Omega, |\A| \Leb{N})$, since $|\A| \in L^\infty(\Omega)$, while trivially $v^\lambda \in L^{1}(\Omega, |\div \A|)$, being $\div\A=0$. Hence, $v \in X^{\A, \lambda}(\Omega)$. On the other hand, $v\not\in BV_{\rm loc}(\Omega)$ given that
\begin{equation*}
D_{x_N}v = {\rm p.v.} \frac{1}{( (\cdot)_N - y_N)} \, \text{ on } B_r(y),
\end{equation*}
so that the distributional gradient of $v$ cannot be a vector valued Radon measure on $\Omega$.
However, we have 
\begin{equation*}
\div(v \A)= \sum_{j = 1}^{N-1} D_j(v a_j) = \left ( \log(|(\cdot)_N - y_N|) a_j \sum_{j=1}^{N-1} \partial_{x_j} \eta_{r, y} \right ) \Leb{N}
\end{equation*} 
since $a_j$ does not depend on $x_j$, for every $j\in\{1,2,\dots,N-1\}$. Thus, as remarked above, $$(\A,Dv)_{\lambda}= \div(v \A) \in \mathcal{M}(\Omega),$$ so that $v\in BV^{\A,\lambda}_{\rm loc}(\Omega)$.
\end{remark}

In addition, we point out that, if $\A \in \DM(\Omega) \setminus \DM^{\infty}(\Omega)$, then the inclusion \eqref{eq:BV_inclusion_BV_A_lambda} may fail to hold true for every $\lambda$. We provide now an example of such a case.

\begin{example}
Let $N = 2$, $\Omega = (-1, 1)^2$, 
$$\A(x_1, x_2) = \frac{(-x_2, x_1)}{x_1^2 + x_2^2}$$
and $u = \chi_E$, where $E = (-1, 1) \times (-1, 0)$. It is clear that $\A \in \DM^p(\Omega)$ for all $p \in [1, 2)$, $u \in BV(\Omega)$ and $u^{\lambda} \in L^1(\Omega, |\div \A|)$ trivially, since $\div \A = 0$. In particular, we easily deduce that $u \in X^{\A, \lambda}(\Omega)$ for all Borel functions $\lambda : \Omega \to [0, 1]$. However, as it was showed in \cite[Remark 4.8]{ChCoTo} (see also \cite[Example 2.5]{MR2532602}), we have
$$\div(u \A) = {\rm p. v.} \left ( \frac{1}{x_1} \right ) \Haus{1} \res (-1, 1) \times \{0\},$$
where ${\rm p.v.}$ stands for principal value integral. Hence, $\div(u \A)$ is not a Radon measure, and therefore, by Proposition \ref{prop:main_inclusions}, this means that $u \notin BV^{\A, \lambda}(\Omega)$ for any Borel function $\lambda : \Omega \to [0, 1]$.
\end{example}

\subsection{Interpolation between $\lambda \equiv 0$ and $\lambda \equiv 1$} Intuitively, the cases $\lambda \equiv 0$ and $\lambda \equiv 1$ represent the extreme choices for $\lambda$, and therefore it is natural to expect that the $\lambda$-pairing is a convex combination of the $0$-pairing and the $1$-pairing. We prove this idea in the following lemma. To this purpose, we notice that, if $u \in BV^{\A,0}(\Omega)$ or $u \in BV^{\A,1}(\Omega)$, then by \eqref{f:senseofmeasures2} we have
\begin{equation} \label{eq:pairing_0_1_expl}
(\A,Du)_0=-u^- \Div\A+\Div(u^- \A) \ \text{ or } \ (\A,Du)_1=-u^+ \Div\A+\Div(u^+ \A),
\end{equation}
respectively.

\begin{lemma} \label{lem:pairings01}	
Let $\A \in \DM_{\rm loc}(\Omega)$. 
\begin{enumerate}[(i)]
\item $u \in X^{\A, 0}(\Omega)\cap X^{\A, 1}(\Omega)$ if and only if $u \in X^{\A, \lambda}(\Omega)$ for every Borel function $\lambda : \Omega \to [0, 1]$.
\item If $u \in BV^{\A, 0}(\Omega)\cap BV^{\A, 1}(\Omega)$ and furthermore $\A \in \DM^1_{\rm loc}(\Omega)$, then we have $u \in BV^{\A, \lambda}(\Omega)$ for every Borel function $\lambda : \Omega \to [0, 1]$, and it holds
\begin{equation} \label{eq:convex_lambda_pairing_0_1}
(\A, Du)_\lambda = (1 - \lambda) (\A, Du)_0 + \lambda (\A, Du)_1 \ \text{ on } \Omega.
\end{equation}
\item If $u \in BV^{\A, 0}(\Omega)\cap BV^{\A, 1}(\Omega)$, then we have $u \in BV^{\A, t}(\Omega)$ for every $t \in [0, 1]$, and \eqref{eq:convex_lambda_pairing_0_1} holds for $\lambda \equiv t$.
\end{enumerate}
\end{lemma}

\begin{proof}
\begin{enumerate}[(i)]
\item If $u \in X^{\A, 0}(\Omega)\cap X^{\A, 1}(\Omega)$, by the boundedness of $\lambda$ we have
\begin{equation}  \label{eq:lambda_u_+_-_summ}
	\lambda u^{+}, (1- \lambda)u^- \in L^1(\Omega, |\A|) \cap L^1(\Omega, |\div \A|)
	\end{equation}
	and so
	\begin{equation*} 
	u^\lambda \in L^1(\Omega, |\A|) \cap L^1(\Omega, |\div \A|).
	\end{equation*}
	The opposite implication is trivial.
\item If $u \in BV^{\A, 0}(\Omega)\cap BV^{\A, 1}(\Omega)$, by (i) we have \eqref{eq:lambda_u_+_-_summ} and $u\in X^{\A,\lambda}(\Omega)$ for every Borel function $\lambda : \Omega \to [0, 1]$.
Therefore, since $\A \in \DM^1_{\rm loc}(\Omega)$, by \eqref{eq:Leibniz_A_vector_field} and \eqref{eq:pairing_0_1_expl} we obtain
\begin{align*}
(\A,Du)_\lambda & = -u^\lambda \Div\A+\Div(u \A) \\
& = - (1- \lambda) u^- \div \A - \lambda u^+ \div \A + (1 - \lambda) \, \Div(u \A) +\lambda \, \Div(u \A) \\
& = (1 - \lambda) \left (- u^- \div \A + \div(u \A) \right ) + \lambda \left (- u^+ \div \A + \div(u \A) \right ) \\
& = (1 - \lambda) (\A, Du)_0 + \lambda (\A, Du)_1 \ \text{ on } \ \Omega.
\end{align*}
Hence, the conclusion follows, since $\pair{\A,Du}_\lambda\in \mathcal{M}(\Omega)$ if
$\pair{\A,Du}_0,\pair{\A,Du}_1\in \mathcal{M}(\Omega)$.
\item Thanks to (i) with $\lambda \equiv t$ for $t \in [0, 1]$, we have $u\in X^{\A,t}(\Omega)$. Then, we argue as in the proof of point (ii), this time exploiting \eqref{f:senseofmeasures2} and noticing that, for $\lambda \equiv t$, we have $\div(u^\lambda \A) = (1 - \lambda) \div(u^- \A) + \lambda \div(u^+ \A)$, since $\lambda$ is constant.
\end{enumerate}
\end{proof}

\begin{remark}
We notice that in point (ii) of Lemma \ref{lem:pairings01} we cannot drop the assumption $\A \in \DM^1_{\rm loc}(\Omega)$. Indeed, if we choose $\A \in \DM(\Omega) \setminus \DM^1_{\rm loc}(\Omega)$ and $u$ as in Remark \ref{rem:A_dm_pure_no_equiv_L_infty}, we know that $u \in BV^{\A, 0}(\Omega)$, with $(\A, Du)_0 = 0$, and it is not difficult to check that actually $u \in BV^{\A, 0}(\Omega) \cap BV^{\A, 1}(\Omega)$, with 
$$(\A, Du)_1 = \div(u^+ \A) = \div \A = 0.$$
However, whenever $\overline{\lambda}(x) = \chi_F(x_1)$, for some Borel set $F \subset (-1, 1)$ such that $\chi_F \notin BV((-1, 1))$, Remark \ref{rem:A_dm_pure_no_equiv_L_infty} shows that $(\A, Du)_{\overline{\lambda}}$ is not a Radon measure, and so $u \notin BV^{\A, \overline{\lambda}}_{\rm loc}(\Omega)$.

Instead, point (iii) of Lemma \ref{lem:pairings01} ensures that, under the general assumption $\A \in \DM_{\rm loc}(\Omega)$, for all $t \in [0, 1]$ the $t$-pairing is the convex combination of the $0$-pairing and the $1$-pairing. In particular, if $u \in BV^{\A, 0}(\Omega)\cap BV^{\A, 1}(\Omega)$, then $u \in BV^{\A}(\Omega) = BV^{\A, \frac{1}{2}}(\Omega)$, with
\begin{equation*}
(\A, Du) = \frac{(\A, Du)_0 + (\A, Du)_1}{2} \text{ on } \Omega.
\end{equation*}
\end{remark}

\subsection{Absolute continuity properties} As recalled in Section \ref{ss:div}, the $\lambda$-pairing measure between a vector field $\A \in \DM^{\infty}_{\rm loc}(\Omega)$ and $u \in BV_{\rm loc}(\Omega)$ with $u^* \in L^1_{\rm loc}(\Omega, |\div \A|)$ enjoys the quite natural absolute continuity property given by \eqref{eq:abs_cont_classic_pairing_lambda} (see also \cite[Proposition 4.4]{CDM}). Indeed, given that, by Proposition \ref{prop:main_inclusions}, under these assumptions our $\lambda$-pairing coincides with the one introduced in \cite{CDM}, we conclude that it enjoys the same absolute continuity. However, it is interesting to see if a similar result can be proved in the case $\A \in \DM^{\infty}_{\rm loc}(\Omega)$, even for functions $u$ which do not have bounded variation.

\begin{proposition} \label{prop:absol} 
Let $\A\in \DM^{\infty}_{\rm loc}(\Omega)$, $\lambda : \Omega \to [0, 1]$ be a Borel function and $\Omega' \Subset \Omega$ be an open set.
If $u\in BV^{\A,\lambda}_{\rm loc}(\Omega)\cap L^\infty_{\rm loc}(\Omega)$, then we have 
$$|(\A,Du)_\lambda| \le 2 c_N \|u\|_{L^\infty(\Omega')} \|\A\|_{L^\infty(\Omega'; \R^N)} \Haus{N-1} \ \text{ on } \Omega',$$ where 
\begin{equation} \label{eq:c_N}
c_N = N \left (\frac{2N}{N+1} \right )^{\frac{N-1}{2}} \frac{\omega_N}{\omega_{N-1}}.
\end{equation}
In addition, if $\A\in \DM^{\infty}(\Omega)$ and $u\in BV^{\A,\lambda}(\Omega)\cap L^\infty(\Omega)$, then the statement holds true with $\Omega$ instead of $\Omega'$.
\end{proposition} 
\begin{proof}
It suffices to note that, if $u\in L^\infty(\Omega')$, then $u\A\in L^\infty(\Omega';\R^N)$. By equation \eqref{eq:Leibniz_A_vector_field} and by point (2) of Proposition \ref{prop:main_inclusions}, we deduce that $\div(u \A) \in \mathcal{M}(\Omega')$, and so $u\A \in \DM^{\infty}(\Omega')$. Hence, thanks to \cite[Proposition 3.1]{Silhavy19}, we get 
\begin{equation*}
|\div \A| \le c_N \|\A\|_{L^\infty(\Omega'; \R^N)} \Haus{N-1} \ \text{ on } \Omega' \ \text{ and } \ |\div(u\A)| \le c_N \|u \A\|_{L^\infty(\Omega'; \R^N)} \Haus{N-1} \ \text{ on } \Omega',
\end{equation*}
where $c_N$ is as in \eqref{eq:c_N}.
Thanks to \eqref{eq:Leibniz_A_vector_field} and the fact that $|u^{\lambda}| \le \|u\|_{L^\infty(\Omega')}$ on $\Omega'$, the conclusion follows. Finally, it is clear that, under global assumptions, we can simply repeat on $\Omega$ the argument above.
\end{proof}

Arguing analogously, we can characterize the absolute continuity properties of the pairing measure in the case $\A \in \DM^p_{\rm loc}(\Omega)$ for $p \in [1, + \infty]$.

\begin{proposition} \label{prop:abs_cont_p}
Let $p, q \in [1, + \infty]$ be conjugate exponents; that is, $\frac{1}{p} + \frac{1}{q} = 1$.
Let $\A\in \DM^{p}_{\rm loc}(\Omega)$ and $\lambda : \Omega \to [0, 1]$ be a Borel function.
\begin{enumerate}[(i)]
\item If $u \in L^\infty_{\rm loc}(\Omega) \cap W^{1, q}_{\rm loc}(\Omega)$, then $$(\A, Du)_\lambda = \A \cdot \nabla u \, \Leb{N}.$$
\item If $N \ge 2$, $p \in \left [\frac{N}{N-1}, + \infty \right )$ and $u \in BV^{\A, \lambda}_{\rm loc}(\Omega) \cap L^\infty_{\rm loc}(\Omega)$, then $|(\A, D u)_\lambda|(B) = 0$ for every Borel set $B \subset \Omega$ which is $\sigma$-finite with respect to the measure $\Haus{N - q}$.
\end{enumerate}
\end{proposition}
\begin{proof}
Assertion (i) follows from the Leibniz rules recalled in Section \ref{ss:div}, see \eqref{eq:Leibniz_Sobolev_classic} and the comments afterwards. Indeed, it is clear that $u \in X^{\A, \lambda}_{\rm loc}(\Omega)$. In addition, if $u \in W^{1, q}_{\rm loc}(\Omega)$, we know that either $u$ admits a continuous representative and so $S_u$ is empty for $q > N$ (by Morrey's inequality), or $u^*(x) = \tilde{u}(x)$ for $|\div \A|$-a.e. $x \in \Omega$, by \cite[Theorem 1.1.24 and Theorem 3.2.2]{Comi}. All in all, this implies that $|\div \A|(S_u) = 0$, and so $u^\lambda(x) = \tilde{u}(x) = u^*(x)$ for $|\div \A|$-a.e. $x \in \Omega$. Therefore, we exploit \eqref{eq:Leibniz_Sobolev_classic} and \eqref{eq:Leibniz_A_vector_field} to obtain
\begin{equation*}
\A \cdot \nabla u \, \Leb{N} = \div(u \A) - u^* \div \A = \div(u \A) - u^\lambda \div \A = (\A, Du)_\lambda
\end{equation*}
As for (ii), we notice that $u\A \in L^p_{\rm loc}(\Omega; \R^N)$. Since $u \in BV^{\A, \lambda}_{\rm loc}(\Omega)$, then $(\A, Du)_\lambda \in \mathcal{M}_{\rm loc}(\Omega)$, and so, thanks to \eqref{eq:Leibniz_A_vector_field}, we conclude that $u\A \in \DM^p_{\rm loc}(\Omega)$. Hence, thanks to the absolute continuity properties recalled in Section \ref{ss:div} (see also \cite{Silh}*{Theorem 3.2}), we exploit again \eqref{eq:Leibniz_A_vector_field} to conclude.
\end{proof}

\begin{remark}
In point (ii) of Proposition \ref{prop:abs_cont_p}, we can exploit \cite[Theorem 2.8]{MR3676052} to conclude that $(A, Du)_\lambda$ vanishes on Borel sets with zero $q$-Sobolev capacity, and the proof is completely analogous.
\end{remark}

\begin{remark} \label{rem:Dirac_delta}
As for the subcritical case $\A \in \DM^p_{\rm loc}(\Omega)$ for $p \in \left [1, \frac{N}{N-1} \right )$, which is not covered by point (ii) of Proposition \ref{prop:abs_cont_p}, we notice that we cannot expect any absolute continuity property for the $\lambda$-pairing. Indeed, we provide an example of a field $\A$, a Borel function $\lambda : \Omega \to [0, 1]$ and a function $u \in BV^{\A, \lambda}(\Omega)$ such that $(\A, Du)_\lambda$ involves a Dirac delta. Let $N \ge 2$, $\Omega = \R^N$ and 
$$\A(x) = \frac{1}{N \omega_N} \frac{x}{|x|^N}.$$ 
Clearly, $\A \in L^p_{\rm loc}(\R^N; \R^N)$ for all $p \in \left [1, \frac{N}{N-1} \right )$, and, by a standard calculation, we see that
$\div \A = \delta_0$, which is the Dirac's delta measure centered in the origin . We choose $u = \chi_{(0, 1)^N}$. Arguing as it was done in \cite[Example 3.1]{ChCoTo} for the special case $N = 2$, we can show that $\div( \chi_{(0, 1)^N} \A) \in \mathcal{M}(\R^N)$ and 
\begin{equation} \label{eq:prod_rule_Dirac_delta}
\div( \chi_{(0, 1)^N} \A) = \frac{1}{2^N} \delta_0 + \overline{(\A, D \chi_{(0, 1)^N})},
\end{equation}
where $\overline{(\A, D \chi_{(0, 1)^N})}$ is the measure acting as 
\begin{align}
\int_{\R^N} \varphi \, d \overline{(\A, D \chi_{(0, 1)^N})}&  = \frac{1}{N \omega_N} \int_{\partial (0, 1)^N \cap \{ x_j > 0, \ \forall j \in \{1, \dots , N\} \}} \varphi(x) \frac{x \cdot \nu_{(0, 1)^N}(x)}{|x|^N} \, d \Haus{N-1}(x) \nonumber \\
& = - \frac{1}{N \omega_N} \sum_{j = 1}^N \int_{\partial (0, 1)^N \cap \{ x_j = 1 \} } \varphi(x) \frac{1}{(1 + |\hat{x}_j|^2)^{\frac{N}{2}}} \, d \Haus{N-1}(x)  \label{eq:def_overline_pairing}
\end{align}
for all $\varphi \in C_c(\R^N)$, where $\hat{x}_j = (x_1, \dots, x_{j-1}, x_{j+1}, \dots, x_N)$.
Thanks to point (2) in Proposition \ref{prop:main_inclusions}, this fact implies that $\chi_{(0, 1)^N} \in BV^{\A, \lambda}(\R^N)$ for all Borel functions $\lambda : \R^N \to [0, 1]$, with
\begin{equation} \label{eq:prod_rule_Dirac_delta_lambda}
(\A, D \chi_{(0, 1)^N})_{\lambda} = \left ( \frac{1}{2^N} - \lambda (0) \right ) \delta_0 + \overline{(\A, D \chi_{(0, 1)^N})}.
\end{equation}
Hence, as long as $\lambda(0) \neq 2^{-N}$, the measure $(\A, D \chi_{(0, 1)^N})_{\lambda}$ cannot be absolutely continuous with respect to $\Haus{\alpha}$ for all $\alpha \in (0, N]$.
In order to prove \eqref{eq:prod_rule_Dirac_delta}, we take $\varphi \in C^1_c(\R^N)$, and we notice that $\A \in C^1(\R^N \setminus B_\delta(0))$ and $\div \A = 0$ on $\R^N \setminus B_\delta(0)$ for all $\delta > 0$. We set $H^+ = \{ x_j > 0, \ \forall j \in \{1, \dots , N\} \}$ and $Q = (0, 1)^N$ for brevity. Hence, we can integrate by parts in the following way:
\begin{align*}
& \int_{\R^N} \chi_{Q} \A \cdot \nabla \varphi \, dx = \lim_{\eps \to 0^+} \int_{Q \setminus B_\eps(0)} \A \cdot \nabla \varphi \, dx \\
& = - \lim_{\eps \to 0^+} \Big ( \int_{\partial Q  \setminus B_\eps(0)} \varphi \A \cdot \nu_{Q} \, d \Haus{N-1} + \int_{\partial B_\eps(0) \cap H^+} \varphi \A \cdot \nu_{\R^N \setminus B_\eps(0)} \, d \Haus{N-1} \Big ) \\
& = - \frac{1}{N\omega_N} \lim_{\eps \to 0^+} \Big (  \int_{\left (\partial Q \setminus B_\eps(0) \right ) \cap H^+} \varphi(x) \frac{x \cdot \nu_{Q}(x)}{|x|^N} \, d \Haus{N-1}(x) + \int_{\partial B_\eps(0) \cap H^+} \frac{\varphi(x)}{|x|^{N-1}} \, d \Haus{N-1}(x) \Big ) \\
& = - \frac{1}{N\omega_N} \int_{\partial Q \cap H^+} \varphi(x) \frac{x \cdot \nu_{Q}(x)}{|x|^N} \, d \Haus{N-1}(x) - \frac{1}{N\omega_N} \lim_{\eps \to 0^+} \int_{\partial B_1(0) \cap H^+} \varphi(\eps x) \, d \Haus{N-1}(x) \\
& = - \int_{\R^N} \varphi \, d \overline{(\A, D \chi_{(0, 1)^N})} - \frac{1}{N\omega_N} \Haus{N-1}(\partial B_1(0) \cap H^+) \varphi(0),
\end{align*}
and this implies \eqref{eq:prod_rule_Dirac_delta}, due to the fact that 
$$\Haus{N-1}(\partial B_1(0) \cap \{ x_j > 0, \ \forall j \in \{1, \dots , N\} \}) = \frac{N\omega_N}{2^N}.$$
All in all, \eqref{eq:Leibniz_A_vector_field} and \eqref{eq:prod_rule_Dirac_delta} imply \eqref{eq:prod_rule_Dirac_delta_lambda}.
\end{remark}

\subsection{The case of $L^1$-divergence}

\begin{proposition} \label{rem:div_A_abs_cont_Leb}
Let $\A\in \DM^1_{\rm loc}(\Omega)$ be such that $\div\A \ll \Leb{N}$. Then, for every Borel function $\lambda:\Omega\to[0,1]$, we have
\begin{equation} \label{eq:div_A_Leb_n_id}
BV^{\A,\lambda}(\Omega) = BV^{\A}(\Omega)= \{u\in \Borel :\  u \in L^1(\Omega, |\A|  + |\Div\A|), \, \Div(u\A)\in \mathcal{M}(\Omega)\}
\end{equation}
and, given $u \in BV^{\A}(\Omega)$ 
\begin{equation} \label{eq:pairing_lambda_ide}
(\A, Du)_\lambda = (\A, Du) = - u \div \A + \div (u \A) \ \text{ on } \Omega.
\end{equation}
In particular, the operator $BV^{\A}(\Omega) \ni u \to (\A, Du) \in \mathcal{M}(\Omega)$ is linear.
\end{proposition}

\begin{proof}
Thanks to \eqref{eq:S_u_*_negl} and \eqref{eq:u_lambda_u_tilde_Leb_N}, we know that, for every $u\in \Borel$ and every Borel function $\lambda : \Omega \to [0, 1]$,
\begin{equation*}
|\div \A|(S_u^*) = 0, \quad \mbox{ and } \quad u^\lambda(x)=u^{\frac{1}{2}}(x) = u(x) \quad \mbox{for $|\Div\A|$-a.e $x\in\Omega$,} 
\end{equation*}
Thus, due to \eqref{eq:Leibniz_A_vector_field}, we get
\begin{equation*}
(\A,Du)_{\lambda} = - u \div \A + \div (u \A) = (\A,Du) \quad \mbox{ for every Borel function $\lambda : \Omega \to [0, 1]$,}
\end{equation*}
and this implies \eqref{eq:div_A_Leb_n_id}, which in turn easily yields the linearity of the operator $u \to (\A, Du)$.
\end{proof}

\begin{remark}
Due to \eqref{eq:lambda_u_tilde_S_u_*}, we notice that, under the assumption $|\div\A|(S_u^*)=0$ for a given $u \in \Borel$, we have $u \in BV^{\A, \lambda_1}(\Omega)$ if and only $u \in BV^{\A, \lambda_2}(\Omega)$ for all Borel functions $\lambda_1, \lambda_2 : \Omega \to [0, 1]$, with
\begin{equation*}
(\A, Du)_\lambda = (\A, Du) = - \tilde{u} \div \A + \div (u \A) \ \text{ on } \Omega,
\end{equation*}
for every Borel function $\lambda : \Omega \to [0, 1]$, which is \eqref{eq:pairing_lambda_ide} with the representative $\tilde{u}$. However, in general we cannot weaken the absolute continuity assumption $\div \A \ll \Leb{N}$, since the set $S_u^*$ could have Hausdorff dimension equal to $N$.
\end{remark}

\begin{remark} \label{rem:A_nu}
If $\A\in \DM^1_{\rm loc}(\Omega)$ satisfies $\div\A=0$, then the results of Proposition \ref{rem:div_A_abs_cont_Leb} can be simplified, since the condition $u \in L^1(\Omega, |\Div\A|)$ is always satisfied, so that \eqref{eq:div_A_Leb_n_id} is improved to
\begin{equation*}
BV^{\A,\lambda}(\Omega) = BV^{\A}(\Omega)= \{u\in \Borel :\  u \in L^1(\Omega, |\A| ), \ \Div(u\A)\in \mathcal{M}(\Omega)\}
\end{equation*}
for every Borel function $\lambda : \Omega \to [0, 1]$, and 
\begin{equation}
(\A,Du)_{\lambda} = (\A,Du) = \Div(u\A).
\label{eq:equalitypairing}
\end{equation}
As a particular case of divergence-free vector fields, let us choose a fixed direction $\nu\in\mathbb{S}^{N-1}$, and consider the constant field $\A:= \nu$. Then, by virtue of \eqref{eq:equalitypairing}, for all $u \in BV^{\nu}(\Omega)$ we have
\begin{equation*}
 (\A,Du) = \Div(u\nu) = D_\nu u \,,
\end{equation*}
where $D_\nu u$ is nothing else than the distributional derivative of $u$ in the direction $\nu$. Note that, if $N = 1$, then $\nu \in \{\pm 1\}$, and so, up to a sign, we simply get the first distributional derivative of $u$: therefore, in the one-dimensional case, if $\A$ is a nonzero constant, we obtain $BV^{\A}(\Omega) = BV(\Omega)$. If instead $N \ge 2$, by the characterization of the differentiability of the precise representative (see \cite[Theorem 3.107]{AFP}) we then have
\begin{equation}
D_\nu u = \Leb{N-1}\res\Omega_\nu\otimes Du_y^\nu\,,
\label{eq:reprdirdev}
\end{equation}
where $\Omega_\nu=\pi_\nu(\Omega)$, $\pi_\nu : \R^N \to \R^N$ is the projection onto the hyperplane $\{ x \in \R^N : x \cdot \nu = 0 \}$, $u_y^\nu(t)=u(y+t\nu)$ for $t\in\R$ and $y\in \Omega_\nu$ such that $y+t\nu\in\Omega$. Although \cite[Theorem 3.107]{AFP} requires $u\in BV(\Omega)$, an inspection to the proof, based on Fubini's Theorem, shows that $D_\nu u\in\mathcal{M}(\Omega)$; i.e., the weak differentiability of $u$ in a fixed direction $\nu$, will suffice to get \eqref{eq:reprdirdev} (see also the remarks in \cite[Section 3.11]{AFP}). We claim that \eqref{eq:reprdirdev} and general results on the disintegration of measures  (see \cite[Section 2.5]{AFP}) imply
\begin{equation*}
(\nu,Du)=D_\nu u \ll  \mathcal{H}^{N-1} \,.
\end{equation*}
Indeed, we first note that, as a consequence of \cite[Theorem 28]{evans2015measure} applied to $\pi_\nu$, $\mathcal{H}^{N-1}(B)=0$ implies $\Leb{N-1}(\pi_\nu(B))=0$ for every Borel set $B\subset\Omega$.
Now, for any such set $B$, if we set $B_\nu := \pi_\nu(B)$ and, for $y\in B_\nu$, $B_y^\nu := \{t \in \R : y+t\nu\in\Omega \}$, we obtain 
\begin{equation*}
|D_\nu u|(B) = \int_{B_\nu}\left(\int_{B_y^\nu} dDu_y^\nu\right)\, d\Leb{N-1}(y) =0\,,
\end{equation*}
since $\Leb{N-1}(B_\nu)=0$ and $|Du_y^\nu|(B_y^\nu)<+\infty$ for $\Leb{N-1}$-a.e. $y\in B_\nu$. 
\end{remark}

\subsection{The class $W^{\A, \lambda}(\Omega)$} We end this section with the definition of a \quotemarks{degenerate} Sobolev-type class. \\
Let $\A\in \DM_{\rm loc}^1(\Omega)$ and $\lambda : \Omega \to [0, 1]$ be a Borel function. We define
\begin{equation*}
W^{\A, \lambda}(\Omega) := \{u\in BV^{\A, \lambda}(\Omega):\,\, (\A,Du)_\lambda \ll \Leb{N}\}\,,
\end{equation*}
\begin{equation*}
W^{\A, \lambda}_{{\rm loc}}(\Omega) := \{u\in BV^{\A, \lambda}_{\rm loc}(\Omega):\,\, (\A,Du)_\lambda \ll \Leb{N}\}\,.
\end{equation*}
Analogously to the case of $BV^{\A}(\Omega)$ (see \eqref{eqdef:BV_A}), we set 
\begin{equation} \label{eqdef:W_A}
W^{\A}(\Omega) := W^{\A, \frac{1}{2}}(\Omega) = \{ u \in BV^{\A}(\Omega) : (\A,Du) \ll \Leb{N}\}\,.
\end{equation}
and analogously for the local class.

We point out that the assumption $\A \in L^1_{\rm loc}(\Omega; \R^N)$ is not restrictive. Indeed, given $\A \in \DM_{\rm loc}(\Omega)$ and $u \in \Lip_{\rm loc}(\Omega)$, we have only $(\A, Du) \ll |\A|$, thanks to Remark \ref{rem:Silhavy_DM} and \eqref{eq:abs_cont_pairing_irreg} (see also \cite[Propositions 2.1 and 2.2]{MR2532602}), which implies that $(\A, Du)_\lambda$ may not be absolutely continuous with respect to the Lebesgue measure even for a regular function $u$, as long as $\A$ is a singular measure.

In addition, we note that, analogously to the inclusion \eqref{eq:BV_inclusion_BV_A_lambda} for $BV^{\A, \lambda}$, classical Sobolev functions with suitable summability for their precise representatives belong indeed to $W^{\A, \lambda}$, as long as $\A$ satisfies a natural summability assumption.

\begin{proposition} \label{prop:Sobolev_inclusion}
Let $p, q \in [1, + \infty]$ be conjugate exponents; that is, $\frac{1}{p} + \frac{1}{q} = 1$, and $\A\in \DM^{p}(\Omega)$. Then for all Borel functions $\lambda : \Omega \to [0, 1]$ we have
\begin{equation*}
W^{1,q}(\Omega) \cap X^{\A}(\Omega) \subset W^{\A, \lambda}(\Omega)
\end{equation*}
with strict inclusion, and $(\A, Du)_\lambda = (\A \cdot \nabla u) \Leb{N}$ for all $u \in W^{1,q}(\Omega) \cap X^{\A}(\Omega)$.
\end{proposition}

\begin{proof}
Let $\lambda : \Omega \to [0, 1]$ be a Borel function. Arguing as in the proof of Proposition \ref{prop:abs_cont_p}, we see that, if $u \in W^{1, q}(\Omega)$, then $|\div \A|(S_u) = 0$ and $u^\lambda(x) = \tilde{u}(x) = u^*(x)$ for $|\div \A|$-a.e. $x \in \Omega$. Hence, $W^{1,q}(\Omega) \cap X^{\A}(\Omega) = W^{1,q}(\Omega) \cap X^{\A, \lambda}(\Omega)$. Let now $u \in W^{1,q}(\Omega) \cap X^{\A}(\Omega)$. For $k \in \N$, we set $T_k : \R \to \R$ to be the truncation map
\begin{equation} \label{eq:truncation}
T_k(u) = \begin{cases} k & \text{ if } u > k, \\
u & \text{ if } |u| \le k, \\
- k & \text{ if } u < - k. \end{cases}
\end{equation}
It is clear that $T_k(u) \in L^\infty(\Omega) \cap W^{1, q}(\Omega)$ for all $k \in \N$. By Proposition \ref{prop:abs_cont_p}, we have 
\begin{equation*}
(\A, D T_k(u))_\lambda = (\A \cdot \nabla T_k(u)) \Leb{N} = (\A \cdot \nabla u ) \Leb{N} \res \{ |u| < k \},
\end{equation*}
where $\A \cdot \nabla u \in L^1(\Omega)$.
Hence, for all $\varphi \in C^1_c(\Omega)$, thanks to Lebesgue's Dominated Convergence Theorem with respect to the measure $|\div \A|$, we get
\begin{align*}
\int_{\Omega} \varphi (\A \cdot \nabla u ) \, dx & = \lim_{k \to + \infty} \int_{\Omega} \varphi \, d (\A, D T_k(u))_\lambda = - \lim_{k \to + \infty} \int_{\Omega} \widetilde{T_k(u)} \, d \div \A + \int_{\Omega} T_k(u) \A \cdot \nabla \varphi \, dx \\
& = - \int_{\Omega} \tilde{u} \, d \div \A - \int_{\Omega} u \A \cdot \nabla \varphi \, dx.
\end{align*}
Thus, by Definition \ref{d:lambdapairbis}, we conclude that $(\A, Du)_\lambda = (\A \cdot \nabla u) \Leb{N}$, and so $u \in W^{\A, \lambda}(\Omega)$.
Finally, in order to prove that the inclusion is in general strict, we consider the example given in Remark \ref{eq:examplev}.
For such choices of the vector field $\A$ and the function $v$, we clearly have $v\not\in W^{1,q}_{\rm loc}(\Omega)$ and $v \in W^{\A, \lambda}(\Omega)$ for all Borel functions $\lambda : \Omega \to [0, 1]$, given that $v \in BV^{\A, \lambda}(\Omega)$ and it satisfies $(\A,Dv)_\lambda \ll \Leb{N}$. This ends the proof.
\end{proof}

\begin{remark} \label{rem:W_A_div_A_abs_cont}
If $\A\in \DM_{\rm loc}^1(\Omega)$ satisfies $\div\A \ll \Leb{N}$, then, arguing as in the proof of Proposition \ref{rem:div_A_abs_cont_Leb}, we get $W^{\A, \lambda}(\Omega) = W^{\A}(\Omega)$ for every Borel function $\lambda : \Omega \to [0, 1]$. In particular, for $u \in BV^{\A}(\Omega)$ \eqref{eq:pairing_lambda_ide} implies
\begin{equation*}
(\A,Du) \ll \Leb{N} \iff \Div(u\A) \ll \Leb{N} \,,
\end{equation*}
so that we get 
\begin{equation*}
W^{\A}(\Omega) = \{ u \in X^{\A}(\Omega) : \div(u \A) \in \mathcal{M}(\Omega) \text{ and } \Div(u\A) \ll \Leb{N} \}.
\end{equation*}
\end{remark}

\section{Lower semicontinuity properties and approximations results} \label{s:lsc}

Let $\A\in \DM_{\rm loc}(\Omega)$ and $\lambda\colon \Omega \to [0,1]$ be a Borel function. In this section we will study some lower semicontinuity and continuity properties of the $\lambda$-pairing and its total variation in the class $BV^{\A,\lambda}_{\rm{loc}}(\Omega)$, with respect to a suitable notion of convergence.
Since \(\pair{\A,Du}_\lambda\) is affected by the pointwise value of $u^\lambda$, the natural notion of convergence in $BV^{\A,\lambda}_{\rm{loc}}(\Omega)$ involves the function $\lambda$.

\begin{definition}
\label{d:convergence}
Let $\A\in \DM_{\rm loc}(\Omega)$ and $\lambda: \Omega \to [0, 1]$ be a Borel function.
We say that a sequence $(u_n)_n \subset X^{\A, \lambda}_{\rm loc}(\Omega)$ $(\A, \lambda)$-converges to $u\in X^{\A, \lambda}_{\rm loc}(\Omega)$ 
if 
\begin{enumerate}
\item $u_n^\lambda \weakto u^\lambda$ in $L^1_{\rm{loc}}(\Omega,|\A|)$,
\item $u^\lambda_n \weakto u^\lambda$ in $L^1_{\rm{loc}}(\Omega, |\Div \A|)$.
\end{enumerate}
\end{definition}

\begin{remark}\label{newconverg} We list some particular cases in which the $(\A, \lambda)$-convergence is easier to check.
\begin{enumerate}[i)]
\item If $(|\A| + |\div \A|) \ll \Leb{N}$, then $\lambda$ does not play any role in the convergence, given that $u^\lambda(x) = u(x)$ for $(|\A| + |\div \A|)$-a.e. $x \in \Omega$, due to \eqref{eq:u_lambda_u_tilde_Leb_N}. In such a case, we omit the $\lambda$ in the notation for the $(\A, \lambda)$-convergence and simply refer to it as $\A$-convergence.
\item If $\A \in L^\infty(\Omega; \R^N)$ and $|\A| \ge c$ for some $c > 0$, then condition (1) is equivalent to the weak convergence in $L^1_{\rm loc}(\Omega)$.
\item If $\Div\A=0$, then condition (2) can be dropped, so that the $(\A, \lambda)$-convergence reduces to the weak convergence in $L^1_{\rm{loc}}(\Omega,|\A|)$.
\end{enumerate}
\end{remark}

We have the following lower semicontinuity and continuity results.

\begin{theorem}
\label{t:lscnuovo11}
Let $\A\in \DM_{\rm loc}(\Omega)$ and $\lambda:\Omega\to [0,1]$ be a Borel function.
Then for every sequence $(u_n)_n \subset X^{\A, \lambda}_{\rm loc}(\Omega)$ and for every $u\in X^{\A, \lambda}_{\rm loc}(\Omega)$
and such that $(u_{n})_{n}$ $(\A, \lambda)$-converges to $u$, it holds
\begin{equation}\label{f:liminf1}
\pscal{\pair{\A,Du}_{\lambda}}{\varphi}
=
\lim_{n\to+\infty} \pscal{\pair{\A,Du_n}_{\lambda}}{\varphi}
\qquad
\forall \varphi\in C^1_c(\Omega)
\end{equation}
in the sense of distributions. In addition, if $u, u_n \in BV^{\A, \lambda}(\Omega)$ for all $n \in \N$, then
\begin{equation}\label{f:liminf3}
|\pair{\A,Du}_{\lambda}|(\Omega)
\leq
\liminf_{n\to+\infty} |\pair{\A,Du_n}_{\lambda}|(\Omega).
\end{equation}
and, if $\sup_{n \in \N} |(\A, Du_n)_\lambda|(\Omega) < + \infty$, we get 
\begin{equation}\label{f:liminf2}
(\A,Du_n)_{\lambda} \weakto (\A,Du)_{\lambda} \ \text{ in } \, \mathcal{M}(\Omega).
\end{equation}
\end{theorem}

\begin{proof}
We see that \eqref{f:liminf1} is a consequence of the definition of $(\A,\lambda)$-convergence; indeed
\begin{align*}
\pscal{\pair{\A,Du}_{\lambda}}{\varphi}=&
-\int_\Omega \varphi u^\lambda\, d \Div\A
- \int_\Omega u^\lambda\, \nabla\varphi \cdot d \A \nonumber\\
=& \lim_{n\to+\infty}\left\{
-\int_\Omega \varphi u_{n}^\lambda \, d \Div\A
- \int_\Omega u_{n}^\lambda \, \nabla\varphi \cdot d \A \right\}\nonumber\\
=&\lim_{n\to+\infty} \pscal{\pair{\A,Du_n}_{\lambda}}{\varphi}.
\end{align*}
Then, if $u, u_n \in BV^{\A, \lambda}(\Omega)$ for all $n \in \N$, we exploit the fact that the $\lambda$-pairings are all finite Radon measures to take the supremum in $\varphi$ with $\|\varphi\|_{L^\infty(\Omega)} \le 1$ to get \eqref{f:liminf3} from \eqref{f:liminf1}. Finally, if the sequence $(\A,Du_n)_{\lambda}$ is equibounded in total variation, by approximating a test function in $C_c(\Omega)$ with one in $C^1_c(\Omega)$ in a standard way, we obtain \eqref{f:liminf2}.
\end{proof}

We consider now the existence of smooth approximations for the $\lambda$-pairing. While the problem seems difficult to tackle in full generality, in the case $\lambda \equiv \frac{1}{2}$, we adapt a standard mollification technique to our general setting, by asking that the field $\A$ is summable and by adding a suitable assumption on the concentration of the measure $|\div \A|$. 

\begin{theorem}\label{t:lscnuovo}
If $\A \in \DM^1(\Omega)$, for every $u\in BV^{\A}(\Omega) \cap L^1_{\rm loc}(\Omega)$ such that $|\div \A|(S_u \setminus J_u) = 0$ there exists a sequence $(u_k)_k \subset BV^{\A}(\Omega) \cap L^\infty(\Omega) \cap C^\infty(\Omega)$ converging to $u$ in $L^1(\Omega,|\A|\Leb{N})$ and to $u^*$ in $L^1(\Omega, |\div \A|)$, such that
\begin{equation}\label{f:limsup1}
\lim_{k \to + \infty} \int_\Omega\varphi\A\cdot\nabla u_k \,dx = \int_\Omega \varphi \, d \pair{\A,Du}
\qquad
\forall \varphi\in C^1_c(\Omega).
\end{equation}
\end{theorem}

We notice that Theorem \ref{t:lscnuovo} actually implies that $(u_k)_k$ $(\A, \frac{1}{2})$-converges to $u$. Although the proof of this result requires a standard approach, we provide the details for the ease of the reader, given that, differently from the established literature \cite{ComiLeo,CD3, CDM}, we do not assume that the scalar function $u$ belongs to the space $BV$.

\begin{proof}
We first notice that $u^*(x) = u^{\frac{1}{2}}(x)$ for $|\div \A|$-a.e. $x \in \Omega$, thanks to $|\div \A|(S_u \setminus J_u) = 0$, see Section \ref{ss:def}. We start by assuming $u \in L^\infty(\Omega)$, and we set $u_{\eps} = \rho_{\eps} \ast u$, for some standard mollifier $\rho$. It is clear that $u_\eps \to u$ in $L^1(\Omega)$, $u_\eps(x) \to u^*(x)$ for all $x \in \Omega \setminus (S_u \setminus J_u)$, so that $u_\eps(x) \to u(x)$ for $\Leb{N}$-a.e. $x \in \Omega$, and $|u_\eps(x)| \le \|u\|_{L^{\infty}(\Omega)}$ for all $x \in \Omega$ and $\eps > 0$. Hence, thanks to Lebesgue's Dominated Convergence Theorem with respect to the measures $|\A| \Leb{N}$ and $|\div \A|$, we see that
$$\|u_\eps - u\|_{L^1(\Omega, |\A| \Leb{N})} + \| u_\eps - u^*\|_{L^1(\Omega, |\div \A|)} \to 0.$$
Therefore, \eqref{f:limsup1} immediately follows from Definition \ref{d:lambdapairbis} for $\lambda \equiv \frac{1}{2}$. In the general case of $u \in BV^{\A}(\Omega) \cap L^1_{\rm loc}(\Omega)$, we consider the truncation of $u$, $T_k(u)$, defined as in \eqref{eq:truncation}.
Since $u^*(x) = u^{\frac{1}{2}}(x)$ for $|\div \A|$-a.e. $x \in \Omega$, we see that $u^* \in L^1(\Omega, |\div \A|)$ and $T_k(u)^*(x) \to u^*(x)$ as $k \to + \infty$ for $|\div \A|$-a.e. $x \in \Omega$. Due to the fact that $|T_k(u)^*| \le |u^*|$ on $\Omega \setminus (S_u \setminus J_u)$ (see \cite[Equation (2.7)]{ComiLeo} and the preceding comments), by Lebesgue's Dominated Convergence Theorem with respect to the measure $|\div \A|$, we get $T_k(u)^* \to u^*$ in $L^1(\Omega, |\div \A|)$. Analogously, we know that $u \in L^1(\Omega, |\A| \Leb{N})$, so that $T_k(u) \to u$ in $L^1(\Omega, |\A| \Leb{N})$ again by Lebesgue's Dominated Convergence Theorem. Now, we consider $u_{\eps, k} = (T_k(u))_\eps$ and we see that
\begin{equation*}
\lim_{k \to + \infty} \lim_{\eps \to 0} \int_{\Omega} |u_{\eps, k} - u| |\A| \, dx + \int_{\Omega} |u_{\eps, k} - u^*| d |\div \A| = 0.
\end{equation*} 
Hence, via a diagonal argument, we may find a sequence $(u_{\eps_k})_k$ such that $u_{\eps_k} \to u$ in $L^1(\Omega, |\A| \Leb{N})$ and $u_{\eps_k} \to u^*$ in $L^1(\Omega, |\div \A|)$. Thus, \eqref{f:limsup1} immediately follows again from Definition \ref{d:lambdapairbis} in the case $\lambda \equiv \frac{1}{2}$.
\end{proof}

\begin{remark}
As for the general $\lambda$-pairing, we point out that there is an approximation result for $BV$ functions, recently obtained in \cite[Theorem 3.4]{ComiLeo}.
If $\A \in \DM^\infty(\Omega)$ and $\lambda : \Omega \to [0, 1]$ is a Borel function, for every $u\in BV(\Omega)$ with $u^\lambda \in L^1(\Omega, |\div \A|)$ there exists a sequence $(u_k^\lambda)_k \subset BV(\Omega) \cap L^\infty(\Omega) \cap C^\infty(\Omega)$ converging to $u$ in $BV(\Omega)$-strict (that is, $u_k^\lambda \to u$ in $L^1(\Omega)$ and $|D u_k^\lambda|(\Omega) \to |Du|(\Omega)$) such that
\begin{equation*}
(\A \cdot \nabla u_k^\lambda) \, \Leb{N} \weakto \pair{\A, Du}_\lambda \ \text{ in } \, \mathcal{M}(\Omega).
\end{equation*}
We notice that, if we also assume that $u \in L^\infty(\Omega)$, then we obtain $u_k^\lambda \to u^\lambda$ in $L^1(\Omega, |\div \A|)$, so that the sequence $(u_k^\lambda)_k$ $(\A, \lambda)$-converges to $u$. Indeed, the convergence in $L^1(\Omega)$ implies the one in $L^1(\Omega, |\A| \Leb{N})$, given that $\A \in L^\infty(\Omega; \R^N)$. As for the convergence in $L^1(\Omega, |\div \A|)$, it follows by Lebesgue's Dominated Convergence Theorem, since, by \cite[Theorem 3.2 and Theorem 3.4]{ComiLeo}, we have $\|u_k^\lambda\|_{L^\infty(\Omega)} \le 2 \|u\|_{L^\infty(\Omega)}$ and $u_k^\lambda(x) \to u^{\lambda}(x)$ for $\Haus{N - 1}$-a.e. $x \in \Omega$ as $k \to + \infty$, and so $u_k^\lambda(x) \to u^{\lambda}(x)$ for $|\div \A|$-a.e. $x \in \Omega$, given that $\div \A \ll \Haus{N-1}$ (see Section \ref{ss:div}).
\end{remark}

\section{A linear space contained in $BV^{\A, \lambda}(\Omega)$}  \label{s:banach}

Let $\A \in \DM_{\rm loc}(\Omega)$. We define a new subclass of functions summable with respect to the measures $|\A|$ and $|\Div \A|$:
\begin{equation*}
X^{\A,+}(\Omega):=\{u\in \Borel : |u|^+ \in L^1(\Omega, |\A|) \cap L^1(\Omega, |\Div\A|) \}.
\end{equation*}
We notice that 
\begin{equation} \label{eq:X_A_+_inclusion}
X^{\A,+}(\Omega)\subseteq
X^{\A,\lambda}(\Omega)
\end{equation}
for every Borel function $\lambda\colon \Omega \to [0,1]$. Indeed, thanks to the inequality
$$|u|^-=\min\{|u^+|,|u^-|\}\leq \max\{|u^+|,|u^-|\}= |u|^+,$$
it is easy to see that
\begin{equation*}
|u^\lambda|\leq\lambda|u^+|+(1-\lambda)|u^-|\leq\lambda|u|^++(1-\lambda)|u|^+\leq |u|^+
\end{equation*}
for every Borel function $\lambda\colon \Omega \to [0,1]$.

Moreover, we note that $X^{\A,+}(\Omega)$ is a linear space. Indeed, for every $u,v \in X^{\A,+}(\Omega)$ we see that
\begin{equation*}
\begin{split}
\int_\Omega |u+v|^+ \, d \mu \leq \int_\Omega (|u|+|v|)^+ \, d \mu \leq 2 \int_{\Omega\cap \{|v|\leq|u|\}} |u|^+ \, d \mu + 2 \int_{\Omega\cap \{|u|<|v|\}} |v|^+ \, d \mu <+\infty \,,
\end{split}
\end{equation*}
for $\mu = |\A|$ and $\mu = |\Div \A|$, thanks to the fact that $|f|^+\leq|g|^+$, if $|f| \le |g|$.


At this point, we can define the related $BV^{\A}$-like class:
\begin{equation*}
BV^{\A, +}(\Omega) := \left\{ u\in X^{\A, +}(\Omega):\ (\A, Du)_\lambda \in \mathcal{M}(\Omega) \ \text{ for every Borel } \lambda : \Omega \to [0, 1] \right\}.
\end{equation*}
It is not difficult to see that, whenever $\A \in \DM^1_{\rm loc}(\Omega)$, we get
\begin{equation} \label{eq:lin_space}
BV^{\A, +}(\Omega) = \left\{
u\in X^{\A, +}(\Omega):\ \Div(u\A)\in \radon
\right\}.
\end{equation} 
Indeed, given $u \in X^{\A, +}(\Omega)$, by \eqref{eq:X_A_+_inclusion} we know that $u \in X^{\A, \lambda}(\Omega)$ for every Borel function $\lambda : \Omega \to [0, 1]$. Hence, thanks to point (2) of Proposition \ref{prop:main_inclusions}, we see that $(\A, Du)_\lambda \in \mathcal{M}(\Omega)$ if and only if $\Div(u\A)\in \radon$.

We employ all these remarks in order to explore the case in which $BV^{\A, +}(\Omega)$ enjoys a linear structure.

\begin{proposition}\label{prop:property}
Let $\A \in \DM_{\rm loc}(\Omega)$. 
\begin{enumerate}
\item We have $BV^{\A, +}(\Omega) \subseteq BV^{\A,\lambda}(\Omega)$ for every Borel function $\lambda\colon \Omega \to [0,1]$. 
\item If $\A \in \DM^1_{\rm loc}(\Omega)$, then $BV^{\A, +}(\Omega)$ is a linear space. 
\item If $\A \in \DM^1_{\rm loc}(\Omega)$ with $\Div\A = \Div^a\A\Leb{N}$ for some $\div^a \A \in L^1_{\rm loc}(\Omega)$, then 
$$BV^{\A, +}(\Omega) = BV^{\A}(\Omega),$$
so that $BV^{\A}(\Omega)$ is a linear space satisfying 
\begin{equation*}
BV^{\A}(\Omega) = \{ u \in \Borel : \|u\|_{BV^{\A}(\Omega)} < + \infty \},
\end{equation*}
where $\|\cdot\|_{BV^{\A}(\Omega)}$ is a seminorm defined as
\begin{equation}\label{norm}
\|u\|_{BV^{\A}(\Omega)} :=\|u\|_{L^1(\Omega, |\A| \Leb{N})}+\|u\|_{L^1(\Omega,|\Div^a\A|\Leb{N})}
+|\Div(u\A)|(\Omega).
\end{equation}
If in addition $\Leb{N}(\Omega \setminus {\rm supp}(|\A|)) = 0$, then $BV^{\A}(\Omega)$ is a Banach space, endowed with the norm given by \eqref{norm}.
\end{enumerate}
\end{proposition}
\begin{proof} 
Point (1) is a trivial consequence of the definition of $BV^{\A, +}(\Omega)$. As for point (2), the linearity follows from \eqref{eq:lin_space} and the fact $X^{\A,+}(\Omega)$ is a linear space. Concerning point (3), the absolute continuity $(|\A| + |\div \A|) \ll \Leb{N}$ implies
\begin{equation*}
|u|^+(x) = |u(x)| \ \text{ for } (|\A| + |\div \A|)\text{-a.e. } x \in \Omega, 
\end{equation*}
due to \eqref{eq:u_lambda_u_tilde_Leb_N}, and therefore $BV^{\A, +}(\Omega) = BV^{\A}(\Omega)$ (see Proposition \ref{rem:div_A_abs_cont_Leb}), with $u \in BV^{\A}(\Omega)$ if and only if $u \in \Borel$ and $\|u\|_{BV^{\A}(\Omega)} < + \infty$.
Finally, it is easy to check that \eqref{norm} is a norm, under the assumption $\Leb{N}(\Omega \setminus {\rm supp}(|\A|)) = 0$, so that we are left to show the completeness. Let $(u_n)_n$ be a Cauchy sequence in $BV^{\A}(\Omega)$, i.e. for every $\varepsilon>0$ there exists $n_0\in\N$ such that for every $n\geq n_0$ we have
\begin{equation}\label{cauchy}
\|u_n-u_{n+1}\|_{BV^{\A}(\Omega)}<\varepsilon.
\end{equation}
Then $(u_n)_n$ is a Cauchy sequence in $L^1(\Omega, |\A| \Leb{N})$ and in $L^1(\Omega,|\Div^a\A|\Leb{N})$. Since these spaces are Banach, there exists two functions $u\in L^1(\Omega, |\A| \Leb{N})$ and $v\in L^1(\Omega,|\Div^a\A|\Leb{N})$ such that
$$
\lim_{n \to + \infty} \int_\Omega|u_n-u|\,|\A| \,dx = 0\,,\qquad \lim_{n \to +\infty} \int_\Omega|u_n-v||\Div^a\A|\,dx = 0\,.
$$
Hence, there exists a subsequence $(u_{n_k})$ such that $u_{n_k}(x) \to u(x)$ for $|\A| \Leb{N}$-a.e. $x \in \Omega$, and, thanks to the assumption $\Leb{N}(\Omega \setminus {\rm supp}(|\A|)) = 0$, this implies that $u_{n_k}(x) \to u(x)$ for $\Leb{N}$-a.e. $x \in \Omega$. Therefore, we exploit Fatou's Lemma with respect to $\Leb{N}$ to get
\begin{equation*}
\int_{\Omega} |u - v| |\div^a \A| \, dx \le \int_{\Omega} \liminf_{k \to + \infty} |u_{n_k} - v| |\div^a \A| \, dx \le \liminf_{k \to + \infty}  \int_\Omega|u_{n_k}-v||\Div^a\A|\,dx = 0,
\end{equation*}
and so we conclude that $u=v$ for $|\Div^a\A| \Leb{N}$-a.e. $x\in\Omega$. Thus, $u_n \to u$ in $L^1(\Omega, |\A| \Leb{N}) \cap L^1(\Omega, |\div^a \A| \Leb{N})$. Moreover, by \eqref{cauchy} the sequence of measures $\mu_n:=\Div(u_n\A)$ is a Cauchy sequence, and therefore it is uniformly bounded, so that there exists a measure $\mu$ and a subsequence $(\mu_{n_j})$ such that $\mu_{n_j}$ weakly converges to $\mu$ in $\mathcal M(\Omega)$. It remains to prove that $\mu=\Div(u\A)$. We recall that, since $\div \A \ll \Leb{N}$, all pairings are identical, see \eqref{eq:pairing_lambda_ide}, so that we just consider $(\A, Du)$. By \eqref{f:liminf1} for all $\varphi\in C^1_c(\Omega)$ we have
\begin{equation}\label{f:liminf11}
\pscal{\pair{\A,Du}}{\varphi}
=
\lim_{n\to+\infty} \pscal{\pair{\A,Du_n}}{\varphi} = \lim_{n\to+\infty} \int_{\Omega} \varphi \, d \pair{\A,Du_n}.
\end{equation}
Hence, we get
\begin{align*}
\int_\Omega\varphi\,d\div(u\A) & =
-\int_\Omega \varphi u\, \Div^a\A\,dx
+\pscal{\pair{\A,Du}}{\varphi} \\
&= \lim_{j \to+\infty}\left (
-\int_\Omega \varphi u_{n_j}  \Div^a\A\,dx
+ \int_{\Omega} \varphi \, d \pair{\A,Du_{n_j}} \right)\nonumber\\
&= \lim_{j \to+\infty}\int_\Omega\varphi\,d\div(u_{n_j} \A)=\int_\Omega\varphi\,d \mu.
\end{align*}
Therefore $\mu=\Div(u\A)$, so that $\Div(u\A)\in\mathcal M(\Omega)$, and we conclude that $u\in BV^{\A}(\Omega)$.
\end{proof}

\begin{corollary} \label{cor:A_Du_seminorm}
Let $\A \in \DM^1_{\rm loc}(\Omega)$ be such that $\Div\A \ll \Leb{N}$, with $\div \A = \div^a \A \Leb{N}$ for some $\div^a \A \in L^1_{\rm loc}(\Omega)$. Then the functional 
$$BV^{\A}(\Omega) \ni u \to |(\A, Du)|(\Omega)$$ 
is a seminorm, and an equivalent seminorm on $BV^{\A}(\Omega)$ is given by
\begin{equation*}\label{eq:equiv_semi_norm}
BV^{\A}(\Omega) \ni u \to \|u\|_{L^1(\Omega, |\A| \Leb{N})}+\|u\|_{L^1(\Omega,|\Div^a \A|\Leb{N})}
+|(\A, Du)|(\Omega).
\end{equation*}
\end{corollary}

\begin{proof}
Thanks to Proposition \ref{rem:div_A_abs_cont_Leb}, we know that the operator $BV^{\A}(\Omega) \ni u \to (\A, Du)$ is linear. In particular, this immediately implies the total variation of the pairing is 1-homogeneous and satisfies the triangle inequality. Then, given $u \in BV^{\A}(\Omega)$, we exploit \eqref{eq:pairing_lambda_ide} twice, in order to get 
\begin{align*}
|(\A, Du)|(\Omega) & \le \|u\|_{L^1(\Omega,|\Div^a \A|\Leb{N})} + |\div(u\A)|(\Omega), \\
|\div(u\A)|(\Omega) & \le \|u\|_{L^1(\Omega,|\Div^a \A|\Leb{N})} + |(\A, Du)|(\Omega).
\end{align*}
The rest of the proof follows from point (3) of Proposition \ref{prop:property}.
\end{proof}

Under the assumption that both the vector field $\A$ and its divergence are locally summable we prove that also the Sobolev-like class $W^{\A}(\Omega)$ is indeed a linear space (in such a case, $\lambda$ does not play any role, as noticed in Remark \ref{rem:W_A_div_A_abs_cont}).

\begin{corollary}
Let $\A \in \DM^1_{\rm loc}(\Omega)$ be such that $\Div\A \ll \Leb{N}$, with $\div \A = \div^a \A \Leb{N}$ for some $\div^a \A \in L^1_{\rm loc}(\Omega)$. Then $W^{\A}(\Omega)$ is a linear space on which we define the seminorm
\begin{equation}\label{eq:semi_norm_W}
\|u\|_{W^{\A}(\Omega)} := \|u\|_{L^1(\Omega, |\A| \Leb{N})}+\|u\|_{L^1(\Omega,|\Div^a \A|\Leb{N})} +\|(\A, Du)^a \|_{L^1(\Omega)},
\end{equation}
where $(\A, Du) = (\A, Du)^a \Leb{N}$ for some $(\A, Du)^a \in L^1(\Omega)$.
If in addition $\Leb{N}(\Omega \setminus {\rm supp}(|\A|)) = 0$, then $W^{\A}(\Omega)$ is a Banach space, endowed with the norm given by \eqref{eq:semi_norm_W}.
\end{corollary}

\begin{proof}
Thanks to Corollary \ref{cor:A_Du_seminorm}, we know that the pairing is linear in the second component. This implies that $W^{\A}(\Omega)$ is a linear space. Then, it is easy to check that $\|\cdot\|_{W^{\A}(\Omega)}$ is a seminorm, and a norm whenever $\Leb{N}(\Omega \setminus {\rm supp}(|\A|)) = 0$. It remains to prove the completeness. Let $(u_n)_{n \in \N}$ be a Cauchy sequence in $W^{\A}(\Omega)$. Since $W^{\A}(\Omega) \subset BV^{\A}(\Omega)$, and, by Proposition \ref{prop:property}, $BV^{\A}(\Omega)$ is a Banach space, we know that there exists $u \in BV^{\A}(\Omega)$ such that 
$$\lim_{n \to + \infty} \|u_n - u\|_{BV^{\A}(\Omega)} = 0.$$
Hence, by the definition of $W^{\A}(\Omega)$ \eqref{eqdef:W_A}, it remains to check that $(\A, Du) \ll \Leb{N}$. To this purpose, we notice that, for all $n \in \N$, there exists $(\A, Du_n)^a \in L^1(\Omega)$ such that $(\A, Du_n) = (\A, Du_n)^a \Leb{N}$. Hence, $((\A, Du_n)^a)_{n \in \N}$ is a Cauchy sequence in $L^1(\Omega)$, and therefore it admits a limit $\xi_{\A, u} \in L^1(\Omega)$. Due to \eqref{eq:pairing_lambda_ide} and the fact that $u_n \to u$ in $L^1(\Omega, (|\A| + |\Div^a \A|) \Leb{N})$, we see that
\begin{align*}
\int_{\Omega} \varphi \, \xi_{\A, u} \, dx & = \lim_{n \to + \infty} \int_{\Omega} \varphi \, (\A, Du_n)^a \, dx = - \lim_{n \to + \infty} \left ( \int_{\Omega} \varphi \, u_n \, \div^a \A \, dx + \int_{\Omega} u_n \, \A \cdot \nabla \varphi \, dx \right ) \\
& = - \int_{\Omega} \varphi \, u \, \div^a \A \, dx - \int_{\Omega} u \, \A \cdot \nabla \varphi \, dx = \int_{\Omega} \varphi \, d (\A, Du)
\end{align*}
for all $\varphi \in C^1_c(\Omega)$. This proves that $(\A, Du) = \xi_{\A, u} \Leb{N}$, and therefore $u \in W^{\A}(\Omega)$.
\end{proof}

\begin{remark}
We point out the assumption $\Leb{N}(\Omega \setminus {\rm supp}(|\A|)) = 0$ in point (3) of Proposition \ref{prop:property} is necessary in order to avoid the degeneracy of the seminorm $\| \cdot \|_{BV^{\A}(\Omega)}$. Indeed, consider $\A \in C^1_c(\Omega; \R^N)$ such that $V = \Omega \setminus {\rm supp}(|\A|)$ is a non-empty open set. Then, given any nontrivial function $u \in C^1_c(\Omega)$ such that ${\rm supp}(u) \subset V$, we clearly have $u \in BV^{\A}(\Omega)$, with $\| u \|_{BV^{\A}(\Omega)} = 0$, since $u$ and $\A$ are both regular and have disjoint supports, so that, in particular $$(\A, Du) = \A \cdot \nabla u \Leb{N} = 0.$$
\end{remark}

Thanks to Proposition \ref{prop:property}, we know that, if $\A \in \DM^{1}_{\rm loc}(\Omega)$ is such that $\div \A \ll \Leb{N}$, then $BV^{\A}(\Omega)$ is a linear space endowed with a seminorm. Hence, it is natural to ask whether, under such conditions, it enjoys some local compactness with respect to the $\A$-convergence (see Remark \ref{newconverg}). This would be relevant since the seminorm given by the total variation of the pairing is lower semicontinuous with respect to the $\A$-convergence (Theorem \ref{t:lscnuovo11}). However, it is important to point out that $BV^{\A}(\Omega)$ is not locally compact with respect to such convergence, at least in dimension $N \ge 2$. In other words, we can find a field $\A$ and a sequence $(u_k)$ which is uniformly bounded in $BV^{\A}(\Omega)$, but does not admit an $\A$-converging subsequence. The counterexample below shows the occurrence of this pathological phenomenon even for a \quotemarks{smooth} transversal vector field.

\begin{example} \label{exmp:not_compactness}
Let $N \ge 2$, $\Omega = (-1, 1)^N$, $f \in C^1_c(\R)$ such that $f(0) \neq 0$, $\A(x) = (f(x_N), 0, \dots, 0)$. It is clear that $\A \in \DM^{\infty}(\Omega)$ and $\div \A = 0$. For $k \in \N$, $k \ge 1$, we define
\begin{equation*}
u_k(x) = k \chi_{\left (-1, 1 \right )^{N-1} \times \left (0, \frac{1}{k} \right ) }(x),
\end{equation*}
so that $D_j u_k = 0$ for all $j \in \{ 1, \dots, N -1 \}$ and 
\begin{equation*}
D_N u_k = k \left ( \Haus{N-1} \res (-1, 1)^{N-1} \times \{0\} - \Haus{N-1} \res (-1, 1)^{N-1} \times \left \{ \frac{1}{k} \right \} \right ).
\end{equation*}
Therefore, by \eqref{eq:equalitypairing} we get
\begin{equation*}
(\A, Du_k) = \div(u_k \A) = D_1 ( u_k \A_1) = 0,
\end{equation*}
since $u_k$ and $\A_1$ are constant in $x_1$. Hence, we get $u_k \in BV^{\A}(\Omega) = BV^{\A, \lambda}(\Omega)$ for all Borel functions $\lambda : \Omega \to [0, 1]$, due to Proposition \ref{rem:div_A_abs_cont_Leb}. In addition, we have
\begin{align*}
\|u_k\|_{BV^{\A}(\Omega)} & = \|u_k\|_{L^1(\Omega, |\A|)} + \|u_k^{\frac{1}{2}}\|_{L^1(\Omega, |\div \A|)} + |(\A, Du_k)|(\Omega) = 2^{N-1}k \int_0^{\frac{1}{k}} |f(x_N)| \, dx_N \\
& \le 2^{N-1} \|f\|_{L^{\infty}(\R)}
\end{align*}
for all $k \in \N, k \ge 1$. Hence, $(u_k)_{k \in \N}$ is a uniformly bounded sequence in $BV^{\A}(\Omega)$, while it is clearly not so in $BV(\Omega)$. However, it cannot admit a converging subsequence in $L^1(\Omega, |\A|)$, since we have $u_k(x) \to 0$ for all $x \in \Omega$, so that the only limit could be $u = 0$, but we have
\begin{equation*}
\|u_k\|_{L^1(\Omega, |\A|)} = 2^{N-1}k \int_0^{\frac{1}{k}} |f(x_N)| \, dx_N \to 2^{N-1} |f(0)| \neq 0.
\end{equation*}
In particular, this rules out the existence of any $\A$-converging subsequences, even in the case in which $BV^{\A}(\Omega)$ is a Banach space; that is, whenever $\Leb{1}((-1, 1) \setminus {\rm supp}(f)) = 0$ (by point (3) of Proposition \ref{prop:property}).

Furthermore, we notice that we cannot have even the weak convergence $u_k \weakto 0$ in $L^1(\Omega, |\A|)$: indeed, it is not difficult to see that
\begin{equation*}
u_k \Leb{N} \weakto \Haus{N-1} \res (-1, 1)^{N-1} \times \{0\} \ \text{ in } \ \mathcal{M}(\Omega),
\end{equation*}
so that, for all $\phi \in C_c(\Omega)$, we get
\begin{align*}
\int_\Omega \phi(x) u_k(x) |\A(x)| \, dx & \to \int_\Omega \phi(x) |f(x_N)| \, d \Haus{N-1} \res (-1, 1)^{N-1} \times \{0\} \\
& = |f(0)| \int_{(-1, 1)^{N-1}} \phi(x_1, \dots, x_{N-1}, 0) \, dx_1 \dots d x_{N-1}.
\end{align*}
\end{example}

Nevertheless, we can prove the existence of minimizers for functionals involving the pairing and a forcing term. For instance, we consider the following family of functionals:
\begin{equation}\label{functio}
\mathcal{E}_p(u):= |(\A,Du)|(\Omega) +  \|u-g\|_{L^p(\Omega, |\A| \Leb{N})}\,,\quad u\in BV^{\A}(\Omega) \cap L^p(\Omega, |\A| \Leb{N})\,,
\end{equation}
for $1\leq p\leq +\infty$, $g\in L^{p}(\Omega, |\A| \Leb{N})$ and $\A \in \DM^1_{\rm loc}(\Omega)$ with $\div \A \ll \Leb{N}$, and under other suitable additional assumptions on the vector field.

\begin{theorem} \label{thm:min_funct_E}
Let $\A\in \DM^1(\Omega)$ be such that $|\A| \ge c$ for some $c > 0$ and $\Div\A\in L^{p'}(\Omega)$. Let $g\in L^{p}(\Omega, |\A| \Leb{N})$ for some $1\leq p\leq +\infty$.
Then the functional $\mathcal{E}_p$ admits a minimizer in $BV^{\A}_{\rm{loc}}(\Omega)$.
\end{theorem}

\begin{proof} Let $(u_k)_k\subset BV^{\A}_{\rm{loc}}(\Omega)$ be a minimizing sequence. In particular, $(u_k)_k$ is equibounded in $L^p(\Omega, |\A| \Leb{N})$. Therefore, we can find $u\in L^p(\Omega, |\A| \Leb{N})$ such that $u_k\weakto u$ in $L^p(\Omega, |\A| \Leb{N})$, up to a subsequence. Hence, we see that 
 $u_k \weakto u$ in $L^1(\Omega,|\A| \Leb N)$, since $L^\infty(\Omega, |\A| \Leb{N}) \subset L^{r}(\Omega, |\A| \Leb{N})$ for all $r \ge 1$, due to the fact that $|\A| \in L^1(\Omega)$.  In addition, since $|\A| \ge c > 0$, we get 
 \begin{equation*}
 \|u_k-g\|_{L^p(\Omega, |\A| \Leb{N})} \ge c \|u_k-g\|_{L^p(\Omega)},
 \end{equation*}
so that we conclude that $(u_k)_k$ is equibounded in $L^p(\Omega)$ and therefore $u_k\weakto u$ in $L^p(\Omega)$, up to a further subsequence. In particular, this implies $u_k \weakto u$ in $L^1(\Omega, |\Div \A|\Leb N)$, up to another subsequence.
 From Theorem~\ref{t:lscnuovo11} and the lower semicontinuity of the norm in $L^p(\Omega)$ with respect to the weak-$L^p$ convergence, we have
 \begin{equation}\label{semi}
\mathcal{E}_p(u)\leq \liminf_{k\to+\infty}\mathcal{E}_p(u_k).
\end{equation}
Therefore, $u$ is a minimizer of $\mathcal{E}$ and the proof is concluded.
\end{proof}

We notice that the assumptions on the vector field in Theorem \ref{thm:min_funct_E} allow for fields with some vanishing components. Indeed, the transversal vector field $\A$ used in the counterexample to the compactness of $BV^{\A}(\Omega)$, Example \ref{exmp:not_compactness}, satisfies these assumptions as long as we require $|f(x_N)| \ge c > 0$ for $x_N \in (-1, 1)$. We also note that, even though such assumptions imply that $BV^{\A}(\Omega)$ is a Banach space, thanks to point (3) of Proposition \ref{prop:property}, we still need a forcing term in $\mathcal{E}_p$ to achieve the existence of a minimizer, due to the lack of local compactness in $BV^{\A}(\Omega)$ (see Example \ref{exmp:not_compactness}).

\section{A coarea formula} \label{s:coarea}

In analogy with the already established theory of $\lambda$-pairings \cite{CDM}, it is relevant to ask whether our pairing enjoys some type of coarea formula. Even though the problem seems to be non trivial, we are able to present the following result.
 
\begin{theorem}\label{t:coarea}
	Let $\A\in\DM_{\rm loc}(\Omega)$, let $\lambda: \Omega \to [0, 1]$ be a Borel function and let $u\in BV^{\A, \lambda}_{\rm{loc}}(\Omega)$ be such that 
	\begin{equation} \label{ass_1}
	\lambda u^{+}, (1- \lambda)u^- \in L^1_{\rm loc}(\Omega, |\A|) \cap L^1_{\rm loc}(\Omega, |\div \A|)
	\end{equation}
	 and 
	 \begin{equation} \label{ass_2}
	|\A|(N_t) +  |\div \A|(N_t) = 0 \, \text{ for } \, \Leb{1}\text{-a.e. } t \in \R,
	 \end{equation}
	 where $N_t := \{u^- \le t < u^+\} \setminus \{ u > t \}^{1/2}$.
	 Then we have
	\begin{equation}\label{coarea}
	\int_{\Omega} \varphi \, d (\A, Du)_{\lambda} = \int_{- \infty}^{+\infty} \pscal{(\A, D\chiut)_\lambda}{\varphi}\, dt
	\qquad\forall \varphi\in C^1_c(\Omega)
	\end{equation}
	in the sense of distributions.
	In particular, we get
	\begin{equation} \label{coarea_ineq}
	|(\A, Du)_\lambda| \le \int_{-\infty}^{+\infty} |(\A, D\chiut)_\lambda| \, dt
	\end{equation}
	 in the sense of measures, whenever the right hand size is well posed.
\end{theorem}

\begin{proof}
For all \(\varphi\in C^1_c(\Omega)\), by \eqref{eq:further} for every $t\in\R$ we have 	
	\begin{equation*}
	\pscal{(\A, D\chiut)_{\lambda}}{\varphi} := 
	-\int_{\Omega} \chi_{\{u > t\}}^{\lambda} \varphi\, d\Div\A -
	\int_{\Omega} \chi_{\{u > t\}}^{\lambda} \nabla\varphi \cdot d \A = - \int_{\Omega} \chi_{\{u > t\}}^{\lambda} \, d \div(\varphi \A)
	\end{equation*}
	in the sense of distributions. Thanks to Lemma \ref{funzcaracter}, we have that
\begin{equation} \label{eq:chi_lambda_t_repr}
	\chiut[\lambda](x) =(1- \lambda(x))\chi_{\{u^- > t\}}(x) + \lambda(x)\chi_{\{u^+ > t\}}(x) \quad  \forall x\in \Omega\setminus N_t,
\end{equation}
where $N_t := \{u^- \le t < u^+\} \setminus \{ u > t \}^{1/2}$, and we also have 
\begin{equation*}
|\div(\varphi \A)|(N_t) \le \int_{N_t} |\varphi| \, d |\div \A| + \int_{N_t}  |\nabla \varphi| \, d|\A| = 0
\end{equation*} 
for $\Leb{1}$-a.e. $t \in \R$, thanks to \eqref{eq:leibnitzC1} and \eqref{ass_2}. Therefore, \eqref{eq:chi_lambda_t_repr} holds for $|\div(\varphi \A)|$-a.e. $x \in \Omega$. In addition, Lemma \ref{lem:comp_supp_div_A_0} implies that $\div(\varphi \A)(\Omega) = 0$, since $\varphi$ has compact support. 
	Hence, by exploiting these facts together with \eqref{eq:further}, Cavalieri formula and Fubini's Theorem, we get
	\begin{align*}
	\int_{- \infty}^{+\infty} \pscal{(\A, D\chiut)_{\lambda}}{\varphi}\, dt  & =
	-\int_{-\infty}^{+\infty} \int_{\Omega} \chi_{\{u > t\}}^{\lambda} \, d \div(\varphi \A) \, dt	\\
	& = - \int_{0}^{+ \infty} \int_{\Omega} \left ( (1- \lambda)\chi_{\{u^- > t\}} + \lambda \chi_{\{u^+ > t\}} \right ) d \div(\varphi \A) \, dt \\
	& \,\,\,\,\,\, +  \int_{- \infty}^{0} \int_{\Omega} \left (1 - (1- \lambda)\chi_{\{u^- > t\}} - \lambda \chi_{\{u^+ > t\}} \right ) \, d \div(\varphi \A) \, dt \\
	& = - \int_{\Omega} \int_{0}^{+ \infty} (1 - \lambda) \chi_{\{u^- > t\}} \, dt \, d \div(\varphi \A)  \\
	& \,\,\,\,\,\, + \int_{\Omega} \int_{-\infty}^{0} (1 - \lambda) \left ( 1-  \chi_{\{u^- > t\}} \right ) \, dt \, d \div(\varphi \A) \\
	& \,\,\,\,\,\, - \int_{\Omega} \int_{0}^{+ \infty} \lambda \chi_{\{u^+ > t\}} \, dt \, d \div(\varphi \A) +\\
	& \,\,\,\,\,\, + \int_{\Omega} \int_{- \infty}^{0} \lambda (1 - \chi_{\{u^+ > t\}}) \, dt \, d \div(\varphi \A) \\
	& = - \int_{\Omega} (1 - \lambda) u^- \, d \div(\varphi \A) - \int_{\Omega} \lambda u^+ \, d \div(\varphi \A) \\
	& = - \int_{\Omega} u^{\lambda} \, d \div(\varphi \A) = \int_\Omega \varphi \, d (\A, Du)_{\lambda}.
	\end{align*}
	This proves \eqref{coarea}, while \eqref{coarea_ineq} follows immediately by taking the supremum in $\varphi \in C^1_c(\Omega)$ with $\|\varphi\|_{L^\infty(\Omega)} \le 1$.
\end{proof}

\begin{remark} \label{rem:assumption_coarea_classic}
We notice that the assumptions \eqref{ass_1} and \eqref{ass_2} in Theorem \ref{t:coarea} are always satisfied in the case $|\A| \ll \Leb{N}$ and $\div \A \ll \Leb{N}$, since $N_t \subset S_u^*$ and $\Leb{N}(S_u^*) = 0$ by \eqref{eq:S_u_*_negl}; and $u^+(x) = u^-(x) = u(x)$ for $\Leb{N}$-a.e. $x \in \Omega$, so that 
\eqref{ass_1} is implied by $u \in L^1_{\rm loc}(\Omega, |\A|) \cap L^1_{\rm loc}(\Omega, |\div \A|)$, which is one of the conditions for having $u \in BV^{\A, \lambda}_{\rm loc}(\Omega)$. Actually, in such a case, by Proposition \ref{rem:div_A_abs_cont_Leb} all the $\lambda$-pairing coincide, so that, the coarea formula can be rewritten in the following way:
	\begin{equation*}
	\int_{\Omega} \varphi \, d (\A, Du) = \int_{- \infty}^{+\infty} \pscal{(\A, D\chiut)}{\varphi}\, dt,
	\qquad\forall \varphi\in C^1_c(\Omega).
	\end{equation*}
If instead we have $|\A^s|(N_t) + |\div^s \A|(N_t) = 0$, we can also drop \eqref{ass_2}, while we still need \eqref{ass_1}.

Alternatively, if $\A \in L^\infty_{\rm loc}(\Omega; \R^N)$ and we take $u \in BV_{\rm loc}(\Omega)$, then \eqref{ass_2} follows from to the absolute continuity $\div \A \ll \Haus{N-1}$ and the fine properties of $BV$ functions with respect to the $\Haus{N-1}$-measure (see for instance \cite[Lemma 2.2]{DCFV2}). As for \eqref{ass_1}, we see that $u^{\lambda} \in L^1_{\rm loc}(\Omega, |\div \A|)$ implies $u^{\pm} \in L^1_{\rm loc}(\Omega, |\div \A|)$, thanks to \cite[Lemma 3.2]{CDM}, so that we recover \cite[Theorem 5.1]{CDM}, where \eqref{coarea} holds in the sense of Radon measures (that is, for test functions $\varphi \in C_c(\Omega)$). If $u\in BV_{\rm loc}(\Omega)
\cap L^\infty_{\rm loc}(\Omega)$ and
$\A\in\DM^\infty_{\rm loc}(\Omega)$ with $\Div\A\in L^1_{\rm loc}(\Omega)$, then the assumptions \eqref{ass_1} and \eqref{ass_2} are satisfied and \eqref{coarea_ineq} holds as an equality (see \cite[Theorem 4.4]{CD5}).
\end{remark}

\section{The $(\A, \lambda)$-perimeter} \label{s:perimeter}

In this section, in analogy with the classical theory of functions of bounded variation, we consider the particular case in which $u \in X^{\A, \lambda}_{\rm loc}(\Omega)$ is the characteristic function of a Borel set (up to Lebesgue negligible sets). We notice that, given any Borel set $E$ and Borel function $\lambda : \Omega \to [0, 1]$, we get
\begin{equation}\label{chilambda}
\chi_E^\lambda = (1- \lambda) \chi_{E^1} + \lambda \chi_{E^1 \cup \partial^* E} = \chi_{E^1} + \lambda \chi_{\partial^* E},
\end{equation}
since $\chi_E^- = \chi_{E^1}$ and $\chi_{E}^+ = \chi_{E^1 \cup \partial^* E}$.

\begin{definition}
Let  \(\A\in\DM_{\rm loc}(\Omega)\), $\lambda : \Omega \to [0, 1]$ be Borel function and let $E$ be a Borel subset of $\Omega$. We define the {\it $(\A, \lambda)$-perimeter of $E$ in $\Omega$}, denoted by $P_{\A, \lambda}(E,\Omega)$, as the following variation 
\begin{equation}\label{perimeter}
P_{\A, \lambda}(E,\Omega):=\sup\left\{\int_\Omega \chi_E^\lambda \nabla \varphi\cdot d \A +\int_\Omega \chi^{\lambda}_E\varphi\,d\Div\A: \varphi\in C^1_c(\Omega;\R^N),\ \|\varphi\|_\infty\leq 1\right\}.
\end{equation}
We say that $E$ is {\it a set of finite $(\A, \lambda)$-perimeter in $\Omega$} if $P_{\A, \lambda}(E,\Omega)<+\infty$. Moreover, $E$ is {\it a set of locally finite $(\A, \lambda)$-perimeter in $\Omega$} if $P_{\A, \lambda}(E,\Omega^\prime)<+\infty$ for any open set $\Omega^\prime\Subset \Omega$.
\end{definition}

This notion of perimeter will play a crucial role in establishing Gauss-Green and integration by parts formulas in Section \ref{s:gauss-green}.

By the definition of $\lambda$-pairing (Definition \ref{d:lambdapairbis}), it is immediate to see that
$$
P_{\A,\lambda}(E,\Omega)=|(\A,D\chi_E)_\lambda|(\Omega).
$$
Differently from the classical notion of anisotropic perimeters, this type of perimeter can be zero for nontrivial sets, thus being degenerate at least for some choices of the divergence-measure field $\A$. Indeed, let $\A \equiv v$, for some constant vector $v \neq 0$, and let $E_{\nu} = \{ x \in \R^N : x \cdot \nu > 0 \}$ for some other constant vector $\nu \neq 0$ such that $\nu \cdot v = 0$. Hence, Remark \ref{rem:A_nu} implies that
$$(\A, D \chi_{E_{\nu}})_\lambda = D_v \chi_{E_{\nu}} = 0,$$
so that $P_{\A,\lambda}(E,\Omega) = |(\A, D \chi_{E_{\nu}})_\lambda| (\Omega) = 0$ for all open sets $\Omega$.
This degeneracy creates some difficulties when one wants to follow in the footsteps of the classical theory of sets of finite perimeter; however we are still able to prove some properties.

It is an interesting question to ask whether the $(\A, \lambda)$-perimeter of a set is concentrated on some type of generalized boundary of the set. To this purpose, we recall the definition of another type of Lebesgue measure-invariant boundary of a Borel set $E$:
\begin{equation*}
\partial^- E := \{ x \in \R^N : 0 < \Leb{N}(E \cap B_r(x)) < \Leb{N}(B_r(x)) \ \text{ for all } r > 0\}.
\end{equation*}
We point out that $\partial^- E$ is a closed set and $\partial^* E \subseteq \partial^- E \subseteq \partial E$, possibly with strict inclusions.
We prove below that this boundary contains the support of the pairing distribution $(\A,D\chi_E)_\lambda$, thus covering also the case of Borel sets of infinite $(\A, \lambda)$-perimeter.

\begin{proposition} \label{prop:supp_A_lambda_per}
Let $\A \in \DM_{\rm loc}(\Omega)$, $\lambda : \Omega \to [0, 1]$ be Borel function and $E \subset \Omega$ be a Borel set. Then we have ${\rm supp}((\A,D\chi_E)_\lambda) \subseteq \partial^- E$.
\end{proposition}

\begin{proof}
Since $\partial^- E$ is a closed set, then $\Omega \setminus \partial^- E$ is an open set. Therefore, let $\varphi\in C^1_c(\Omega)$ be such that ${\rm supp}(\varphi) \Subset \Omega \setminus \partial^- E$. By \eqref{eq:further} we get
\begin{equation*}
\pscal{\pair{\A,D\chi_E}_\lambda}{\varphi} = - \int_\Omega \chi_E^\lambda \, d \div(\varphi \A) = - \int_{E^1} \, d \div(\varphi \A) -  \int_{\partial^* E} \lambda \, d \div(\varphi \A).
\end{equation*}
Since $\partial^* E \subseteq \partial^- E$ and ${\rm supp}(\div(\varphi \A)) \subset {\rm supp}(\varphi \A) \Subset \Omega \setminus \partial^- E$, we can conclude that the second term must be zero. As for the first one, we notice that 
\begin{equation*}
\Omega \setminus \partial^- E = E^{1,-} \cup E^{0,-},
\end{equation*}
where
\begin{align*}
E^{1,-} & = \{ x \in \R^N : \text{ there exists } r > 0 \text{ such that } \Leb{N}(E \cap B_r(x)) = \Leb{N}(B_r(x)) \},\\
E^{0,-} & = \{ x \in \R^N : \text{ there exists } r > 0 \text{ such that } \Leb{N}(E \cap B_r(x)) = 0 \}.
\end{align*}
It is easy to check that $E^{1,-} \subseteq E^1$, so that $E^1 \cap (\Omega \setminus \partial^- E) = E^{1,-}$. Now, since ${\rm supp}(\varphi) \Subset \Omega \setminus \partial^- E$, then there exists an open set $V \Subset E^{1,-} \cup E^{0,-}$ such that ${\rm supp}(\varphi) \subset V$. Since $E^{1,-} \cap E^{0,-} = \emptyset$, then $V_0 = V \cap E^{0,-}$ and $V_1 = V \cap E^{1,-}$ are open sets satisfying $V_j \Subset E^{j, -}$ for $j = 0, 1$. All in all, we get
\begin{equation*}
\int_{E^1} \, d \div(\varphi \A) = \int_{E^{1,-}} \, d \div(\varphi \A) = \int_{E^{1,-} \cap V} \, d \div(\varphi \A) = \int_{V_1} \, d \div(\varphi \A) = 0,
\end{equation*} 
thanks to Lemma \ref{lem:comp_supp_div_A_0}, since $\varphi \in C^1_c(V_1 \cup V_0)$, and so $\varphi \in C^1_c(V_1)$, in particular. Therefore, 
\begin{equation*}
\pscal{\pair{\A,D\chi_E}_\lambda}{\varphi} = 0
\end{equation*}
for all $\varphi\in C^1_c(\Omega)$ be such that ${\rm supp}(\varphi) \Subset \Omega \setminus \partial^- E$, and this ends the proof.
\end{proof}

\begin{remark}\label{rem:sizesupp}
We point out that the control on the size of the support of the pairing distribution $(\A, D \chi_E)_\lambda$ given in Proposition \ref{prop:supp_A_lambda_per} is in general too large. Indeed, we can find $\A \in \DM_{\rm loc}(\R^N)$ and a Borel set $E$ such that $(\A, D \chi_E) = 0$, while $\partial^- E$ has Hausdorff dimension equal to $N$. To this purpose, we let $N \ge 2$ and consider the set $F \subset \R$ defined in \cite[Example 3.9]{CCDM}, satisfying ${\rm dim}_{\Haus{}}(\partial^* F) = 1$. Arguing analogously as it was done in \cite[Example 3.9]{CCDM}, we define $$E = F \times \R^{N-1},$$ and we conclude that $${\rm dim}_{\Haus{}}(\partial^* E) = N = {\rm dim}_{\Haus{}}(\partial^- E),$$ since $\partial^* E \subseteq \partial^- E$. It is clear that $\chi_E$ is constant in all variables except for $x_1$. Hence, if now we set $\A = (0, 0, \dots, 0, \chi_E(x))$, we see that $\chi_E \A = \A$ and so $\div \A = \div(\chi_E \A) = 0$. Thanks to Proposition \ref{rem:div_A_abs_cont_Leb}, this implies $\chi_E \in BV^{\A}_{\rm loc}(\Omega) = BV^{\A, \lambda}_{\rm loc}(\R^N)$ for all Borel functions $\lambda : \Omega \to [0, 1]$, with $(\A, D \chi_E)_\lambda = 0$. In addition, this provides an example of a Borel set $E$ such that $\chi_E \in BV^{\A, \lambda}_{\rm loc}(\R^N) \setminus BV_{\rm loc}(\R^N)$ and $P_{\A,\lambda}(E,\Omega) = 0$, while still being nonnegligible with respect to the measure $|\A|$; thus showing the degeneracy of the $(\A, \lambda)$-perimeter.
\end{remark}

In analogy with the classical notion of perimeter, it is interesting to check whether the $(\A, \lambda)$-perimeter satisfies locality, additivity and similar properties.

\begin{proposition} \label{prop:per_prop}
Let $\A \in \DM_{\rm loc}(\Omega)$, $E, F \subset \Omega$ be Borel sets, and $\lambda : \Omega \to [0, 1]$ be a Borel function.
\begin{enumerate}
\item If $\Leb{N}(E \Delta F) = 0$, then $(\A, D \chi_E)_\lambda = (\A, D \chi_F)_\lambda$ in the sense of distributions. 
\item In the sense of distributions, we have 
\begin{equation} \label{eq:pairing_complementary}
(\A, D \chi_E)_\lambda = - (\A, D \chi_{\Omega \setminus E})_{1-\lambda},
\end{equation}
so that $\chi_E \in BV^{\A, \lambda}(\Omega)$ if and only if $\chi_{\Omega \setminus E} \in BV^{\A, 1- \lambda}(\Omega)$, with 
$$P_{\A, \lambda}(E, \Omega) = P_{\A, 1-\lambda}(\Omega \setminus E, \Omega).$$
In particular, if $\lambda \equiv \frac{1}{2}$, then we see that $(\A, D \chi_E) = - (\A, D \chi_{\Omega \setminus E})$ in the sense of distributions, so that $\chi_E \in BV^{\A}(\Omega)$ if and only if $\chi_{\Omega \setminus E} \in BV^{\A}(\Omega)$, with 
$$P_{\A, \frac{1}{2}}(E, \Omega) = P_{\A, \frac{1}{2}}(\Omega \setminus E, \Omega).$$
\item If $\Leb{N}(E \cap F) = 0$ and $(\partial^*E) \cap F^1 = (\partial^*F) \cap E^1 = \partial^*E \cap \partial^*F = \emptyset$, then 
\begin{equation} \label{eq:per_prop_1}
(\A, D \chi_{E \cup F} )_\lambda = (\A, D \chi_E)_\lambda + (\A, D \chi_F)_\lambda,
\end{equation}
in the sense of distributions.
In addition, if $\chi_E, \chi_F \in BV^{\A, \lambda}_{\rm loc}(\Omega)$, then $\chi_{E \cup F} \in BV^{\A, \lambda}_{\rm loc}(\Omega)$, with
\begin{equation} \label{eq:per_prop_2}
|(\A, D \chi_{E \cup F} )_\lambda| \le |(\A, D \chi_E)_\lambda| + |(\A, D \chi_F)_\lambda| \ \text{ on } \Omega.
\end{equation}
Furthermore, if $\partial^- E \cap \partial^- F = \emptyset$, then 
\begin{equation} \label{eq:per_prop_3}
|(\A, D \chi_{E \cup F} )_\lambda| = |(\A, D \chi_E)_\lambda| + |(\A, D \chi_F)_\lambda| \ \text{ on } \Omega.
\end{equation}
\end{enumerate}
\end{proposition}

\begin{proof}
Clearly, $\Leb{N}(E \Delta F) = 0$ implies that $\chi_E(x) = \chi_F(x)$ for $\Leb{N}$-a.e. $x \in \Omega$, and so 
\begin{equation*}
\chi_E^{\lambda}(x) = \chi_F^{\lambda}(x) \ \text{ for all } x \in \Omega.
\end{equation*} 
Hence, the equality of the $\lambda$-pairing distributions in point (1) follows immediately from Definition \ref{d:lambdapairbis}. Then, we note that $\chi_{\Omega \setminus E} = 1 - \chi_E$ implies {$$\chi_{\Omega \setminus E}^+ = 1 - \chi_E^- = 1 - \chi_{E^1}$$ and $$\chi_{\Omega \setminus E}^- = 1 - \chi_E^+ = 1 - \chi_{E^1 \cup \partial^* E}.$$ All in all, we get $$\chi_{\Omega \setminus E}^{\lambda} = 1 - \chi_{E^1} - (1- \lambda) \chi_{\partial^* E} = 1 - \chi_E^{1 - \lambda}.$$} 
Therefore, \eqref{eq:pairing_complementary} is an easy consequence of Lemma \ref{eq:further_def_pairing}; and the rest of point (2) follows immediately. Finally, under the assumptions of point (3), we notice that 
$$\partial^* (E \cup F) = \partial^* E \cup \partial^* F \ \text{ and } \ (E \cup F)^1 = E^1 \cup F^1,$$ 
which follow from \cite[Proposition 2.1]{MR1857126} and a straightforward computation. Therefore, we get $\chi_{E \cup F}^{\lambda} = \chi_E^{\lambda} + \chi_F^{\lambda}$, so that \eqref{eq:per_prop_1} follows from Definition \ref{d:lambdapairbis}. Then, if $\chi_E, \chi_F \in BV^{\A, \lambda}_{\rm loc}(\Omega)$, then \eqref{eq:per_prop_1} easily implies \eqref{eq:per_prop_2}, so that $\chi_{E \cup F} \in BV^{\A, \lambda}_{\rm loc}(\Omega)$. Finally, if $\partial^- E \cap \partial^- F = \emptyset$, thanks to Proposition \ref{prop:supp_A_lambda_per} we see that the pairings $(\A, D \chi_E)_\lambda$ and $(\A, D \chi_F)_\lambda$ have disjoint supports; and so, by taking the total variations in \eqref{eq:per_prop_1}, we deduce \eqref{eq:per_prop_3}.
\end{proof}

\begin{remark} \label{rem:not_lin_pairing}
We point out that only the assumption $\Leb{N}(E \cap F) = 0$ is not enough in point (3) of Proposition \ref{prop:per_prop}. Indeed, let us consider $\Omega = \R^N, E = (-1, 0) \times (0, 1)^{N-1}$ and $F = (0, 1)^N$. In this case, we have $E^1 = E$, $F^1 = F$, $\partial^*E = \partial E$, $\partial^* F = \partial F$ and
\begin{equation*}
(E \cup F)^1 = E^1 \cup F^1 \cup L \ \text{ and } \partial^*(E \cup F) = (\partial^* E \cup \partial^* F) \setminus L,
\end{equation*}
where $L = \{0\} \times (0, 1)^{N-1}$; so that we have
\begin{align}
\chi_{E \cup F}^{\lambda} = \chi_{E^1 \cup F^1 \cup L} + \lambda \chi_{(\partial^* E \cup \partial^* F) \setminus L} & = \chi_{E^1} + \chi_{F^1} + \lambda \chi_{\partial^* E} + \lambda \chi_{\partial^* F} + (1 - 2 \lambda) \chi_{L} \nonumber \\
& = \chi_{E}^{\lambda} + \chi_{F}^{\lambda} + (1 - 2 \lambda) \chi_{L}. \label{eq:chi_E_F_lambda_1}
\end{align}
Therefore, in general we cannot have \eqref{eq:per_prop_1} whenever the set $L \cap \left \{ \lambda \neq \frac{1}{2} \right \}$ is not negligible with respect to the measure $|\A| + |\div \A|$. This happens for instance if we choose $$\A(x) = ( \chi_{\{ x_1 > 0 \}}, 0, \dots, 0),$$ in which case we have $\div \A = \Haus{N-1} \res \{ x_1 = 0\}$. Due to the fact that $\A \in \DM^\infty(\R^N)$ and $E, F$ are sets of finite perimeter, we exploit point (3) of Proposition \ref{prop:main_inclusions} to conclude that $\chi_E, \chi_F, \chi_{E \cup F} \in BV^{\A, \lambda}(\R^N)$ for all Borel functions $\lambda : \R^N \to [0, 1]$. Hence, by \eqref{eq:Leibniz_A_vector_field} and \eqref{eq:chi_E_F_lambda_1} we obtain
\begin{align*}
(\A, D \chi_{E \cup F})_\lambda & = - \chi_{E \cup F}^{\lambda} \div \A + \div(\chi_{E \cup F} \A) \\
& = - \left (\chi_{E}^{\lambda} + \chi_{F}^{\lambda}  + (1 - 2 \lambda) \chi_{L} \right ) \div \A + \div(\chi_E \A) + \div(\chi_F \A) \\
& = (\A, D \chi_E)_\lambda + (\A, D \chi_F)_\lambda  - (1 - 2 \lambda) \Haus{N-1} \res L.
\end{align*}
Hence, the $\lambda$-pairing is not additive, as soon as $\lambda(x) \neq \frac{1}{2}$ for $\Haus{N-1}$-a.e. $x \in L$.
It is important to notice that the case $\lambda \equiv \frac{1}{2}$ is not privileged: to see this, we consider $E = (0, 1)^N$ and $F = (-1, 0)^N$. Again, $E^1 = E$, $F^1 = F$, $\partial^*E = \partial E$, $\partial^* F = \partial F$ and
\begin{align}
\chi_{E \cup F}^{\lambda} = \chi_{E^1 \cup F^1} + \lambda \chi_{\partial^* E \cup \partial^* F} & = \chi_{E^1} + \chi_{F^1} + \lambda \chi_{\partial^* E} + \lambda \chi_{\partial^* F} - \lambda \chi_{\{0\}} \nonumber \\
& =  \chi_{E}^{\lambda} + \chi_{F}^{\lambda}  - \lambda \chi_{\{0\}}, \label{eq:chi_E_F_lambda_2}
\end{align}
where we denote by $0$ the origin of $\R^N$. In this case, \eqref{eq:per_prop_1} fails to hold if we choose 
$$\A(x) = \frac{1}{N \omega_N} \frac{x}{|x|^N} \ \text{ and } \ \lambda (0) \neq 0.$$ 
In this case, we have $\A \in \DM^1_{\rm loc}(\R^N)$ with $\div \A = \delta_0$, and so, arguing as in Remark \ref{rem:Dirac_delta} with minor changes, we get $\chi_E, \chi_F, \chi_{E \cup F} \in BV^{\A, \lambda}(\R^N)$ for all Borel functions $\lambda : \R^N \to [0, 1]$. Thus we exploit \eqref{eq:Leibniz_A_vector_field} and \eqref{eq:chi_E_F_lambda_2} to see that
\begin{align*}
(\A, D \chi_{E \cup F})_\lambda & =  - \chi_{E \cup F}^{\lambda} \div \A + \div(\chi_{E \cup F} \A) \\
& = - \left (\chi_{E}^{\lambda} + \chi_{F}^{\lambda}  - \lambda \chi_{\{0\}} \right ) \div \A + \div(\chi_E \A) + \div(\chi_F \A) \\
& = (\A, D \chi_E)_\lambda + (\A, D \chi_F)_\lambda  + \lambda(0) \delta_0.
\end{align*}
\end{remark}

In analogy with point (4) of Proposition \ref{prop:main_inclusions}, we investigate the relation between $(\A, \lambda)$-perimeters for different choices of $\lambda$, as long as $\A$ is a summable divergence-measure field.

\begin{proposition} \label{prop:diff_pairing_lambda_partial_*_E}
If $\A \in \DM^1_{\rm loc}(\Omega)$ and $E \subset \Omega$ is a Borel set, then for every couple of Borel functions $\lambda_1, \lambda_2 : \Omega \to [0, 1]$ we have
\begin{equation} \label{eq:diff_pairing_lambda_partial_*_E}
(\A, D \chi_E)_{\lambda_1} - (\A, D \chi_E)_{\lambda_2} = (\lambda_2 - \lambda_1) \div \A \res \partial^* E
\end{equation}
in the sense of distributions on $\Omega$. In particular, if $|\div \A|(\partial^* E) = 0$, then
$$(\A, D \chi_E)_{\lambda_1} = (\A, D \chi_E)_{\lambda_2}$$ 
in the sense of distributions on $\Omega$. Therefore, if $\A \in \DM^1(\Omega)$, we get 
\begin{equation*}
P_{\A, \lambda_1}(E, \Omega) < + \infty \ \text{ if and only if } \ P_{\A, \lambda_2}(E, \Omega) < + \infty,
\end{equation*}
and, as long as one of these conditions holds, then \eqref{eq:diff_pairing_lambda_partial_*_E} holds in the sense of Radon measures on $\Omega$.
\end{proposition}

\begin{proof}
We apply the definition of the $\lambda$-pairing distribution \eqref{eq:def_lambda_pairing_abs_cont} to $u = \chi_E$ and $\lambda = \lambda_1$, and $u = \chi_E$ and $\lambda = \lambda_2$, respectively. Then, we take the difference between the two formulas and for all $\varphi \in C^1_c(\Omega)$ we get
\begin{equation*}
\pscal{\pair{\A,D\chi_E}_{\lambda_1} - (\A, D \chi_E)_{\lambda_2}}{\varphi} = - \int_{\Omega} \varphi \left ( \chi_E^{\lambda_1} - \chi_E^{\lambda_2} \right ) \, d \div \A.
\end{equation*}
Therefore, \eqref{eq:diff_pairing_lambda_partial_*_E} follows by exploiting \eqref{chilambda}, where the right hand side clearly vanishes as soon as $\partial^* E$ is $|\div \A|$-negligible. Finally, if $\A \in \DM^1(\Omega)$, we see that the right hand side of \eqref{eq:diff_pairing_lambda_partial_*_E} is always a finite Radon measure, so that $(\A, D \chi_E)_{\lambda_1} \in \mathcal{M}(\Omega)$ if and only if $(\A, D \chi_E)_{\lambda_2} \in \mathcal{M}(\Omega)$.
\end{proof}

\begin{remark}
In light of \eqref{eq:diff_pairing_lambda_partial_*_E}, if $\A \in \DM^1(\Omega)$ and $\lambda(x) = \frac{1}{2}$ for $|\div \A|$-a.e. $x \in \partial^* E$, then we obtain $(\A, D \chi_E)_\lambda = (\A, D \chi_E)$. Therefore, by point (2) of Proposition \ref{prop:per_prop}, these assumptions are also sufficient to ensure 
$$(\A, D \chi_{\Omega \setminus E})_\lambda = (\A, D \chi_{\Omega \setminus E}) = - (\A, D \chi_E) = - (\A, D \chi_E)_\lambda,$$
and therefore $P_{\A, \lambda}(E, \Omega) = P_{\A, \lambda}(\Omega \setminus E, \Omega)$.
\end{remark}

In analogy with Proposition \ref{prop:absol}, we list the absolute continuity properties of the $(\A, \lambda)$-perimeter.

\begin{proposition} \label{prop:abs_cont_per}
Let $\A \in \DM^{\infty}_{\rm loc}(\Omega)$ and $E$ be a Borel set.
\begin{enumerate}[i)]
\item If $E$ is a set of locally finite perimeter in $\Omega$, then $\chi_E\in BV^{\A, \lambda}_{\rm loc}(\Omega)$ for every Borel function $\lambda : \Omega \to [0, 1]$, and we have
\begin{equation}\label{tracce}
(\A, D \chi_E)_\lambda = \left ( (1 - \lambda) {\rm Tr}^i (\A,\partial^*E) + \lambda {\rm Tr}^e (\A,\partial^*E) \right ) \Haus{N-1} \res \partial^{*} E,
\end{equation}
where ${\rm Tr}^i (\A,\partial^*E), {\rm Tr}^e (\A,\partial^*E) \in L^\infty_{\rm loc}(\partial^* E, \Haus{N-1})$.

\item If $\chi_E\in BV^{\A, \lambda}_{\rm loc}(\Omega)$, then for every open set $\Omega' \Subset \Omega$ we have
\begin{equation} \label{eq:abs_cont_gen_per_A_lambda}
|(\A,D\chi_E)_\lambda| \le 2 c_{N} \|\A\|_{L^{\infty}(\Omega'; \R^N)} \Haus{N-1} \res \partial^{-} E \ \text{ on } \Omega',
\end{equation}
where $c_N$ is as in Proposition \ref{prop:absol}.
\end{enumerate}
In addition, if $\A\in \DM^{\infty}(\Omega)$, $E$ is a set of finite perimeter in $\Omega$ in (i), and $\chi_E \in BV^{\A,\lambda}(\Omega)$ in (ii), then the respective statements hold true globally.
\end{proposition}

\begin{proof} 
We notice that clearly $\chi_E^\lambda \in L^1_{\rm loc}(\Omega, |\div \A|)$, so that, by point (3) of Proposition \ref{prop:main_inclusions}, we conclude that $\chi_E \in BV^{\A, \lambda}_{\rm loc}(\Omega)$ for every Borel function $\lambda : \Omega \to [0, 1]$. Hence, by applying \eqref{eq:convex_lambda_pairing_0_1} to every open set $\Omega' \Subset \Omega$, we get 
\begin{equation*}
(\A, D \chi_E)_\lambda = (1 - \lambda) (\A, D \chi_E)_0 + \lambda (\A, D \chi_E)_1 \ \text{ on } \Omega.
\end{equation*}
Hence, the result is a consequence of \eqref{eq:def_normal_traces} and \eqref{eq:normal_traces_L_infty}.
As for point (ii), thanks to Proposition \ref{prop:absol} we see that 
$$|(\A,D\chi_E)_\lambda| \le 2 c_{N} \|\A\|_{L^{\infty}(\Omega'; \R^N)} \Haus{N-1} \ \text{ on } \Omega'.$$ 
Hence, it is enough to apply Proposition \ref{prop:supp_A_lambda_per} to conclude.
Finally, it is clear that, under global assumptions, the normal traces are in $L^\infty(\partial^* E, \Haus{N-1})$ in point (i) and \eqref{eq:abs_cont_gen_per_A_lambda} holds on $\Omega$.
\end{proof}

\begin{remark}
We point out that the inclusion $\partial^* E \subseteq \partial^- E$ might be strict, since $\partial^- E$ can have positive Lebesgue measure (see \cite[Proposition 12.19 and Example 12.25]{Maggi}). In addition, $\Haus{N-1} \res \partial^{-} E$ is a Radon measure if and only if $E$ is a set of locally finite perimeter, in which case $\Haus{N-1}(\partial^- E \setminus \partial^* E) = 0$. Indeed, if $\Haus{N-1} \res \partial^{-} E$ is a Radon measure, then for every compact set $K$ we get
\begin{equation*}
\Haus{N-1}(\partial^{*} E \cap K) \le \Haus{N-1}(\partial^{-} E \cap K) < + \infty.
\end{equation*}
Hence, Federer's Theorem \cite[Theorem 5.23]{evans2015measure} implies that $E$ is a set of locally finite perimeter. Therefore, we can apply De Giorgi's and Federer's Theorems to conclude that the perimeter measure satisfies $|D \chi_E| = \Haus{N-1} \res \partial^* E$ (see for instance \cite[Theorem 3.59 and 3.61]{AFP}). Since the perimeter measure is supported on $\partial^- E$, thanks to \cite[Proposition 12.19]{Maggi}, we conclude that $\partial^- E \setminus \partial^* E$ is $\Haus{N-1}$-negligible.
The reverse implication is instead trivial.
\end{remark}

\begin{remark} \label{rem:strange_cube}
It is interesting to notice that a representation analogous to \eqref{tracce} may hold even for sets which do not have locally finite perimeter. 
To this purpose, in the case $N \ge 2$ we provide an example of a field $\A \in \DM^{\infty}_{\rm loc}(\R^N)$ and of a Borel set $F$ satisfying $\chi_F \notin BV_{\rm loc}(\R^N)$ such that $$(\A, D \chi_F)_\lambda = h \, \Haus{N-1} \res L$$ for every Borel function $\lambda : \R^N \to [0,1]$, where $h \in L^{\infty}_{\rm loc}(L, \Haus{N-1})$ and $L \subset \partial^* F$ is a Borel set such that $\Haus{N-1}(L) < + \infty$. We argue similarly as in \cite[Remark 4.9]{ChCoTo}: we consider the set $E$ defined therein; that is, the open bounded set in $\R^2$ whose boundary is given by
\begin{equation*}
\partial E =\big( \{0\} \times [0, 1]\big) \cup \big([0, 1] \times \{0\}\big) \cup \big([0, 1 + \log{2}] \times \{1\}\big) \cup S,
\end{equation*}
where
\begin{align*}
S & = \Big(\left \{1\right \} \times \left [0, \frac{1}{2} \right ]\Big) \bigcup \Big(\left [1, 2 \right ] \times \left \{ \frac{1}{2}\right \}\Big)
 \bigcup \left (\bigcup_{n \ge 1} \left \{ 1 + \sum_{k = 1}^{n} \frac{(-1)^{k - 1}}{k}\right \} \times \left [ 1 - \frac{1}{2^{n}}, 1 - \frac{1}{2^{n + 1}} \right ]\right )  \\
&\quad  \bigcup \left(\bigcup_{n \ge 1} \left [ 1 + \sum_{k = 1}^{2n} \frac{(-1)^{k - 1}}{k} , 1 + \sum_{k = 1}^{2n + 1} \frac{(-1)^{k - 1}}{k} \right ]\times \left \{ 1 - \frac{1}{2^{2n + 1}}\right \}\right) \\
&\quad \bigcup \left(\bigcup_{n \ge 1} \left [ 1 + \sum_{k = 1}^{2n} \frac{(-1)^{k - 1}}{k} , 1 + \sum_{k = 1}^{2n - 1} \frac{(-1)^{k - 1}}{k} \right] \times \left \{ 1 - \frac{1}{2^{2n}} \right \}\right).
\end{align*}
It is clear that $\Haus{1}(S) = + \infty$. Then, we set $$F = E \times (0, 1)^{N-2},$$ which clearly satisfies $\chi_F \notin BV_{\rm loc}(\R^N)$, and, in particular, $\Haus{N-1}(\partial^* F ) = \Haus{N-1}(\partial^- F) = + \infty$. 
\begin{figure} 
  \centering
      \includegraphics[width=0.6\textwidth]{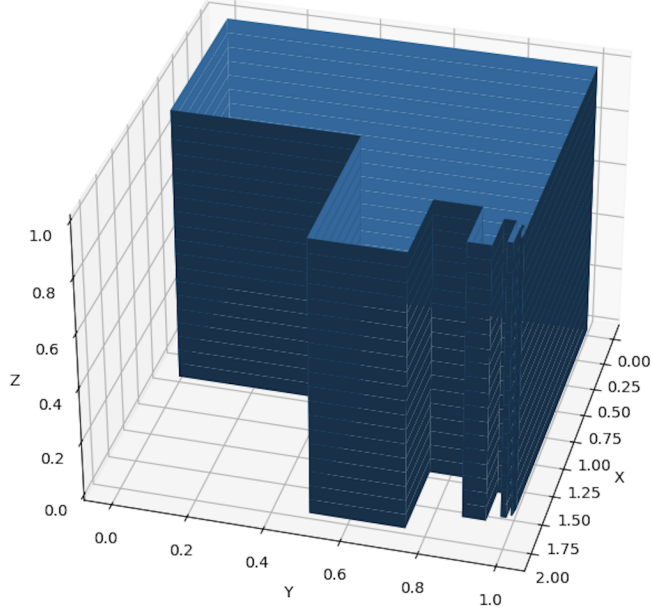}
  \caption{The set $F$ in the case $N = 3$} \label{fig-4.1}
\end{figure}
However, we can exploit the same approach as in \cite[Remark 4.9]{ChCoTo} to show that $D_{x_1} \chi_F \in \mathcal{M}(\R^N)$, with
\begin{align*}
D_{x_{1}} \chi_{F}  = &\, \Haus{N-1} \res \big(\{0\} \times (0, 1)^{N-1}\big) - \Haus{N-1} \res \left (\{1 \} \times \left (0, \frac{1}{2} \right ) \times (0, 1)^{N-2} \right ) \\
& - \Haus{N-1} \res \left ( \bigcup_{n \ge 1}\left\{ 1 + \sum_{k = 1}^{n} \frac{(-1)^{k - 1}}{k} \right \}
        \times \left ( 1 - \frac{1}{2^{n}}, 1 - \frac{1}{2^{n + 1}} \right ) \times (0, 1)^{N-2} \right ).
\end{align*}
Let now $\A(x) = (f(\hat{x}_1) g(x_1), 0, \dots, 0)$, where $\hat{x}_1 =(x_2, x_3, \dots, x_N)$, for some $f \in L^{\infty}_{\rm loc}(\R^{N-1})$ and $g \in C^{1}_{c}(\R)$. It is immediate to see that $\A \in \DM^{\infty}_{\rm loc}(\R^{N})$,
\begin{equation} \label{eq:div_A_example_F}
\div \A = f(\hat{x}_1) g'(x_1) \Leb{N},
\end{equation}
and
\begin{equation*} 
\div(\chi_{F}\A) = D_{x_1}(f(\hat{x}_1) g(x_1) \chi_F(x)) = f(\hat{x}_1) g(x_1) D_{x_1} \chi_{F} + \chi_{F}(x) f(\hat{x}_1) g'(x_1) \Leb{N},
\end{equation*}
since $f$ is constant in $x_1$ and $g \in C^1_c(\R)$ (see \cite[Remark 4.9]{ChCoTo} for more details). Now, since $\div \A \ll \Leb{N}$, thanks to Proposition \ref{rem:div_A_abs_cont_Leb} we know that, for every Borel function $\lambda : \R^N \to [0, 1]$, $BV^{\A, \lambda}(\R^N) = BV^{\A}(\R^N)$ and $\chi_F^\lambda(x) = \chi_F(x)$ for $|\div \A|$-a.e. $x \in \R^N$. Therefore, we get
\begin{equation*}
(\A, D \chi_F)_\lambda = (\A, D \chi_F) = - \chi_F \div \A + \div(\chi_F \A) = f(\hat{x}_1) g(x_1) D_{x_1} \chi_{F} = h \, \Haus{N-1} \res L,
\end{equation*}
where
\begin{align}
L = & \big(\{0\} \times (0, 1)^{N-1}\big) \cup \left (\{1 \} \times \left (0, \frac{1}{2} \right ) \times (0, 1)^{N-2} \right ) \cup \nonumber \\
& \cup \left ( \bigcup_{n \ge 1}\left\{ 1 + \sum_{k = 1}^{n} \frac{(-1)^{k - 1}}{k} \right \}  \times \left ( 1 - \frac{1}{2^{n}}, 1 - \frac{1}{2^{n + 1}} \right ) \times (0, 1)^{N-2} \right ) \label{eq:def_L}
\end{align}
and
\begin{equation} \label{eq:def_h}
h(x) = f(\hat{x}_1) g(x_1) \left ( \chi_{\{0\} \times (0, 1)^{N-1}}(x) - \chi_{L \setminus \left ( \{0\} \times (0, 1)^{N-1} \right )}(x) \right ).
\end{equation}
In particular, we conclude that $\chi_F \in BV^{\A}(\R^N) \setminus BV_{\rm loc}(\R^N)$. One could even define the normal trace of $\A$ on $\partial^* F$ by setting
\begin{equation*}
{\rm Tr} (\A,\partial^*F)(x) := \begin{cases} h(x) & \text{ if } x \in L \\
0 & \text{ if } x \in \partial^*F \setminus L
\end{cases},
\end{equation*}
in this way recovering the representation for the pairing 
\begin{equation} \label{eq:pairing_traccia_A_F}
(\A, D\chi_F) = {\rm Tr} (\A,\partial^*F) \Haus{N-1} \res \partial^* F,
\end{equation}
which is an extension of \eqref{tracce} in the case $\div \A \ll \Leb{N}$, for which interior and exterior normal traces on the measure theoretic boundary of sets of locally finite perimeters coincide (this fact is a simple consequence of \cite[Theorem 4.2]{ComiPayne}, for instance).

Finally, we notice that this example works even for $f \in L^p_{\rm loc}(\R^{N-1})$ for any $p \in [1, + \infty]$, which gives us $\A \in \DM^{p}_{\rm loc}(\R^N)$ and $h \in L^p_{\rm loc}(L, \Haus{N-1})$.
\end{remark}

As an immediate consequence of Proposition \ref{prop:abs_cont_p} and Proposition \ref{prop:supp_A_lambda_per}, we get the following absolute continuity property in the case $\A$ is not essentially bounded.

\begin{proposition} \label{prop:abs_cont_per_unbdd}
Let $N \ge 2$, $p \in \left [\frac{N}{N-1}, + \infty \right)$ and $q \in (1, N]$ be conjugate exponents; that is, satisfying $\frac{1}{p} + \frac{1}{q} = 1$. Let $\A \in \DM^{p}_{\rm loc}(\Omega)$, 
$\lambda : \Omega \to [0, 1]$ be Borel function and $E \subset \Omega$ be a Borel set. If $\chi_E \in BV^{\A, \lambda}_{\rm loc}(\Omega)$, then $|(\A, D\chi_E)_\lambda|(B) = 0$ for every Borel set $B \subset \Omega$ which is $\sigma$-finite with respect to the measure $\Haus{N - q} \res \partial^- E$.
\end{proposition}

As for $BV^{\A, \lambda}$-functions, we consider a notion of $(\A, \lambda)$-convergence for sets of (locally) finite $(\A, \lambda)$-perimeter.

\begin{definition}
\label{d:convergenceper}
Let $\A\in \DM_{\rm loc}(\Omega)$.
We say that a sequence of Borel sets $(E_n)_n$ $(\A, \lambda)$-converges to a Borel set $E$ 
if 
\begin{enumerate}
\item $\chi_{E_n}^\lambda \weakto \chi_E^\lambda$ in $L^1_{\rm{loc}}(\Omega,|\A|)$,
\item $\chi_{E_n}^\lambda  \weakto \chi_E^\lambda$ in $L^1_{\rm{loc}}(\Omega, |\Div \A|)$.
\end{enumerate}
\end{definition}

In the particular case $\A \in \DM^1_{\rm loc}(\Omega)$ with $\div \A \ll \Leb{N}$, we refer to this convergence as $\A$-convergence, as established in Remark \ref{newconverg}.

As an immediate consequence of Theorem \ref{t:lscnuovo11}, we get the following lower semicontinuity property of the $(\A, \lambda)$-perimeter.

\begin{proposition} \label{prop:per_lssc}
The function $E\mapsto P_{\A,\lambda}(E,\Omega)$ is lower semicontinuous with respect to the $(\A, \lambda)$-convergence.
\end{proposition}

The lower semicontinuity of the $(\A, \lambda)$-perimeter naturally suggests the question whether there is any sort of compactness for families of sets with uniformly bounded $(\A, \lambda)$-perimeter. This of course requires first that such perimeter defines some sort of seminorm: due to Corollary \ref{cor:A_Du_seminorm}, we see that this happens as soon as $\A \in \DM^1_{\rm loc}(\Omega)$ satisfies $\div \A \ll \Leb{N}$. However, similarly to the case of $BV^{\A}$-functions (Example \ref{exmp:not_compactness}), in the case of dimension $N \ge 2$, we can find a counterexample to the compactness with respect to the $\A$-convergence.

\begin{example}\label{ex:per_nocompact}
For $k \ge 1$ we set
\begin{equation*}
F_k = \bigcup_{j=0}^{2^{k-1} - 1} \left ( \frac{2j}{2^k}, \frac{2j+1}{2^k} \right ).
\end{equation*}
It is not difficult to see that $\Leb{1}(F_k) = \frac{1}{2}$ for all $k \ge 1$, and that
\begin{equation*}
 \Leb{1} \res F_k \weakto \frac{1}{2} \Leb{1} \res (0, 1) \ \text{ in } \mathcal{M}(\R).
\end{equation*}
In particular, this means that the sequence of sets $(F_k)_{k \ge 1}$ does not admit any subsequence converging in measure.
Let now $N \ge 2$, $\Omega = (-1, 1)^N$ and $\A(x) = (1, 0, \dots, 0)$. It is clear that $\A \in \DM^{\infty}(\Omega)$ and $\div \A = 0$. We set 
$$E_k = (-1, 1)^{N-1} \times F_k,$$ and it is easy to see that $(\A, D \chi_{E_k}) = 0$, so that $P_{\A,\lambda}(E_k,\Omega) = 0$ for all $k \ge 1$. Actually, $\chi_{E_k} \in BV^{\A}(\Omega)= BV^{\A, \lambda}(\Omega)$ for all Borel functions $\lambda : \Omega \to [0, 1]$, due to Proposition \ref{rem:div_A_abs_cont_Leb}. In addition, by Proposition \ref{prop:property}, we know that $BV^{\A}(\Omega)$ is a Banach space, and we see that
\begin{equation*}
\|\chi_{E_k}\|_{BV^{\A}(\Omega)} = \|\chi_{E_k}\|_{L^1(\Omega, |\A|)} = 2^{N-1} \Leb{1}(F_k) = 2^{N-2} \ \text{ for all } k \ge 1.
\end{equation*}
On the other hand, we notice that 
\begin{align*}
\chi_{E_k} |\A| \Leb{N} = \Leb{N-1} \res (-1, 1)^{N-1} \otimes \Leb{1} \res F_k & \weakto \Leb{N-1} \res (-1, 1)^{N-1} \otimes \frac{1}{2} \Leb{1} \res (0, 1) \\
& = \frac{|\A|}{2} \Leb{N} \res (-1, 1)^{N-1} \times (0, 1) \ \text{ in } \mathcal{M}(\Omega).
\end{align*}
Thus, the sequence of sets $(E_k)_{k \ge 1}$ is uniformly bounded in $BV^{\A}(\Omega)$, but it is not weakly compact with respect to the $\A$-convergence.
\end{example}

\section{Gauss-Green and integration by parts formulas} \label{s:gauss-green}

As a remarkable consequence of the Section \ref{subsec:BV_A_lambda} and Section \ref{s:perimeter}, we establish Gauss-Green and integration by parts formulas in our framework.

\begin{theorem} \label{thm:GG}
Let $\A \in \DM_{\rm loc}(\Omega)$ and let $E \Subset \Omega$ be a Borel set. If $\chi_E \in BV^{\A, \lambda}(\Omega)$ for some Borel function $\lambda : \Omega \to [0, 1]$, then we have
\begin{equation}\label{eq:gausslambda}
\div \A(E^1) + \int_{\partial^* E} \lambda \,d\div  \A= - \int_{\partial^- E} \, d (\A, D \chi_E)_\lambda \,.
\end{equation}
If $\chi_E \in BV^{\A, 0}(\Omega)$, then
\begin{equation} \label{eq:GG_i}
\div \A(E^1) = - \int_{\partial^- E} \, d (\A, D \chi_E)_0 \, ;
\end{equation}
if instead $\chi_E \in BV^{\A, 1}(\Omega)
$, then
\begin{equation} \label{eq:GG_e}
\div \A(E^1 \cup \partial^* E) = - \int_{\partial^- E}
 \, d (\A, D \chi_E)_1 \, .
\end{equation}
\end{theorem}

\begin{proof}
It is easy to see that \eqref{eq:gausslambda} follows by applying Lemma \ref{lem:gausslambda_u} to $u = \chi_E$ and exploiting \eqref{chilambda}. Then, \eqref{eq:GG_i} and \eqref{eq:GG_e} are particular cases of \eqref{eq:gausslambda} for $\lambda \equiv 0$ and $\lambda \equiv 1$, respectively.
\end{proof}
 
\begin{corollary} \label{cor:div_A_partial_*_E}
Let $\A \in \DM_{\rm loc}(\Omega)$ and let $E \Subset \Omega$ be a Borel set. If $\chi_E \in BV^{\A, 0}(\Omega) \cap BV^{\A, 1}(\Omega)$, then we have
\begin{equation} \label{eq:div_A_partial_*_E}
\div \A(\partial^* E) = \int_{\partial^- E} \, d \big ( (\A, D \chi_E)_0 - (\A, D \chi_E)_1 \big ) \, .
\end{equation}
In addition, if $\A \in \DM_{\rm loc}^1(\Omega)$ and $E \subset \Omega$ is a Borel set satisfying $\chi_E \in BV^{\A, \lambda}_{\rm loc}(\Omega)$ for some Borel function $\lambda : \Omega \to [0, 1]$, then
\begin{equation} \label{eq:div_A_partial_*_E_equation}
\div \A \res \partial^* E = (\A, D \chi_E)_0 - (\A, D \chi_E)_1 \ \text{ on } \Omega.
\end{equation}
\end{corollary}

\begin{proof}
We obtain \eqref{eq:div_A_partial_*_E} by subtracting \eqref{eq:GG_i} from \eqref{eq:GG_e}. As for \eqref{eq:div_A_partial_*_E_equation}, it is clear that $\A \in \DM^1(\Omega')$ for every open set $\Omega' \Subset \Omega$, therefore the result follows by applying Proposition \ref{prop:diff_pairing_lambda_partial_*_E} to the couple $\lambda_1 = 0$ and $\lambda_2 = 1$.
\end{proof}

\begin{remark}
In the case $\A \in \DM^{\infty}(\Omega)$ and $\chi_E \in BV(\Omega)$, then we know that $\chi_E \in BV^{\A, \lambda}(\Omega)$ for every Borel function $\lambda : \Omega \to [0, 1]$, by Proposition \ref{prop:abs_cont_per}. Hence, $\chi_E \in BV^{\A, 0}(\Omega) \cap BV^{\A, 1}(\Omega)$, and so, taking into account \eqref{tracce}, from the Gauss--Green formulas \eqref{eq:GG_i} and \eqref{eq:GG_e} we retrieve \eqref{eq:classical_Gauss_Green}. Analogously, Corollary \eqref{cor:div_A_partial_*_E} is a generalization of 
\begin{equation*}
\div \A \res \partial^* E = \left ( {\rm Tr}^i (\A,\partial^*E) - {\rm Tr}^e (\A,\partial^*E) \right ) \, \Haus{N-1} \res \partial^* E,
\end{equation*}
for which we refer to \cite[Corollary 3.5 and Theorem 4.2]{ComiPayne}.
 \end{remark}
 
We exploit two examples seen in the previous sections to show some applications of our general Gauss--Green formulas.

\begin{example}
Let $\A(x) = (f(\hat{x}_1) g(x_1), 0, \dots, 0)$, where $\hat{x}_1 =(x_2, x_3, \dots, x_N)$, for some $f \in L^{\infty}_{\rm loc}(\R^{N-1})$ and $g \in C^{1}_{c}(\R)$. Let $F$ be the Borel set in Remark \ref{rem:strange_cube}. Then, we know that $\A \in \DM^\infty_{\rm loc}(\R^N)$ with $\div \A \ll \Leb{N}$, $\chi_F \in BV^{\A}(\R^N) = BV^{\A, \lambda}(\R^N)$ and $(\A, D \chi_F) = (\A, D \chi_F)_\lambda$ for every Borel function $\lambda : \Omega \to [0, 1]$. Thus, we can apply \eqref{eq:gausslambda}, or equivalently \eqref{eq:GG_i} or \eqref{eq:GG_e}, to $\A$ and $F$ to get
\begin{equation} \label{eq:GG_F}
\int_F  f(\hat{x}_1) g'(x_1) \, d x = \div \A(F^1) = - \int_{\partial^- F} \, d (\A, D \chi_F) = - \int_{L} h \, d \Haus{N-1},
\end{equation}
thanks to \eqref{eq:div_A_example_F} and \eqref{eq:pairing_traccia_A_F}, where $L$ is given by \eqref{eq:def_L} and $h$ by \eqref{eq:def_h}.
We point out that \eqref{eq:GG_F} cannot be derived directly from the standard Gauss--Green formula for sets of locally finite perimeter, since $\chi_F \notin BV_{\rm loc}(\R^N)$.
\end{example} 

\begin{example}
Let $N \ge 2$, $\Omega = \R^N$, $\lambda : \R^N \to [0, 1]$ be a Borel function and 
$$\A(x) = \frac{1}{N \omega_N} \frac{x}{|x|^N},$$ 
be as in Remark \ref{rem:Dirac_delta}. We apply \eqref{eq:gausslambda} to $\A$ and $E = (0, 1)^N$, and exploit \eqref{eq:prod_rule_Dirac_delta_lambda} to get
\begin{equation*}
\lambda(0) = \div \A(E^1) + \int_{\partial^*E} \lambda \, d \div \A = - \int_{\partial^- E} \, d (\A, D \chi_E)_\lambda = - \frac{1}{2^N} + \lambda (0) - \overline{(\A, D \chi_{(0, 1)^N})}(\partial (0, 1)^N),
\end{equation*}
where the measure $\overline{(\A, D \chi_{(0, 1)^N})}$ is the one defined in \eqref{eq:def_overline_pairing}. This easily implies
\begin{equation*}
\overline{(\A, D \chi_{(0, 1)^N})}(\partial (0, 1)^N) = - \frac{1}{2^N}.
\end{equation*}
On the other hand, Proposition \ref{prop:supp_A_lambda_per} ensures that $(\A, D \chi_E)_\lambda$ is supported on $\partial^- E = \partial (0, 1)^N$, so that also $\overline{(\A, D \chi_{(0, 1)^N})}$ is supported on $\partial (0, 1)^N$. Hence, by taking $\varphi \equiv 1$ on $(-2, 2)^N$ in \eqref{eq:def_overline_pairing}, we see that
\begin{align*}
\overline{(\A, D \chi_{(0, 1)^N})}(\partial (0, 1)^N)  & = - \frac{1}{N \omega_N} \sum_{j = 1}^N \int_{\partial (0, 1)^N \cap \{ x_j = 1 \} } \frac{1}{(1 + |\hat{x}_j|^2)^{\frac{N}{2}}} \, d \Haus{N-1}(x)  \\
& = - \frac{1}{\omega_N} \int_{(0, 1)^{N-1} } \frac{1}{(1 + |y|^2)^{\frac{N}{2}}} \, d \Leb{N-1}(y),
\end{align*}
since the integrals on each face of the cube are clearly equal.
All in all, by setting $n = N -1$, we deduce directly the following nice identity:
\begin{equation*}
\int_{(0, 1)^{n}} \frac{1}{(1 + |y|^2)^{\frac{n+1}{2}}} \, d y = \frac{\omega_{n+1}}{2^{n+1}} \ \text{ for all } \ n \ge 1.
\end{equation*}
We were not able to find it in literature, and we do believe that a direct computation of such integrals would be a hard task.
\end{example}

As an immediate consequence of Theorem \ref{thm:GG}, we deduce the following general version of the integration by parts formula.
For the reader's convenience, we clarify the notation adopted below. Namely,
$$
BV^{u^{\lambda_1
} \A, \lambda_2}(\Omega)= \left\{
v\in X^{u^{\lambda_1}\A,\lambda_2}(\Omega): (u^{\lambda_1}\A,Dv)_{\lambda_2}\in \radon
\right\}\,,
$$
where $\lambda_1, \lambda_2 : \Omega \to [0, 1]$ are Borel functions and
$$
X^{u^{\lambda_1}
\A,\lambda_2}(\Omega)=\{v \in \Borel : v^{\lambda_2} \in L^1(\Omega, |u^{\lambda_1}\A|) \cap L^1(\Omega, |\Div(u^{\lambda_1}\A)|)
\}\,.
$$

\begin{theorem} \label{thm:IBP_lambda_gen}
Let $\A \in \DM_{\rm loc}(\Omega)$, $\lambda_1, \lambda_2 : \Omega \to [0, 1]$ be Borel functions and $u \in BV^{\A, \lambda_1}(\Omega)$. Let $E \subset \Omega$ be a Borel set such that $\chi_E \in BV^{u^{\lambda_1} \A, \lambda_2}(\Omega)$ and ${\rm supp}(\chi_E^{\lambda_2} u^{\lambda_1} |\A|) \Subset \Omega$. Then we have
\begin{align} \label{eq:IBP_lambda_e}
\int_{E^1} u^{\lambda_1} \, d \div \A + \int_{\partial^* E} \lambda_2 u^{\lambda_1} \, d \div \A  & + \int_{E^1} \, d (\A, Du)_{\lambda_1} + \int_{\partial^*E} \lambda_2 \, d (\A, Du)_{\lambda_1} \\
&  = - \int_{\partial^- E} \, d (u^{\lambda_1} \A, D \chi_E)_{\lambda_2} \, . \nonumber
\end{align}
In particular, if $\lambda : \Omega \to [0, 1]$ is Borel function, $u \in BV^{\A, \lambda}(\Omega)$, $\chi_E \in BV^{u^\lambda \A, 0}(\Omega)$ and they satisfy ${\rm supp}(\chi_E^{-} u^{\lambda} |\A|) \Subset \Omega$, then we have
\begin{equation} \label{eq:IBP_i}
\int_{E^1} u^\lambda \, d \div \A + \int_{E^1} d (\A, Du)_\lambda = - \int_{\partial^- E} \, d (u^\lambda \A, D \chi_E)_0 \, ;
\end{equation}
while, if $u \in BV^{\A, \lambda}(\Omega)$, $\chi_E \in BV^{u^\lambda \A, 1}(\Omega)$ and they satisfy ${\rm supp}(\chi_E^{+} u^{\lambda} |\A|) \Subset \Omega$, then we have
\begin{equation} \label{eq:IBP_e}
\int_{E^1 \cup \partial^*E} u^\lambda \, d \div \A + \int_{E^1 \cup \partial^*E} d (\A, Du)_\lambda = - \int_{\partial^- E} \, d (u^\lambda \A, D \chi_E)_1 \,.
\end{equation}
\end{theorem}

\begin{proof}
Since ${\rm supp}(\chi_E^{\lambda_2} u^{\lambda_1} |\A|) \Subset \Omega$, we get $\div(\chi_E^{\lambda_2} u^{\lambda_1} \A)(\Omega) = 0$ thanks to Lemma \ref{lem:comp_supp_div_A_0}. 
Then, we apply \eqref{f:senseofmeasures2}, to the scalar function $\chi_E$, the field $u^{\lambda_1} \A$ and the Borel function $\lambda_2$, obtaining
\begin{equation*}
\div(\chi_E^{\lambda_2} u^{\lambda_1} \A) = \chi_E^{\lambda_2} \div( u^{\lambda_1} \A) + (u^{\lambda_1} \A, D \chi_E)_{\lambda_2}.
\end{equation*}
Now, we apply again \eqref{f:senseofmeasures2}, this time to the scalar function $u$, the field $\A$ and the Borel function $\lambda_1$, and we get
\begin{equation*}
\div(\chi_E^{\lambda_2} u^{\lambda_1} \A) = \chi_E^{\lambda_2} u^{\lambda_1} \div \A + \chi_E^{\lambda_2} (\A, Du)_{\lambda_1} + (u^{\lambda_1} \A, D \chi_E)_{\lambda_2}.
\end{equation*}
Therefore, by evaluating this identity of measures over $\Omega$ and exploiting \eqref{chilambda}, we obtain \eqref{eq:IBP_lambda_e}. Finally, \eqref{eq:IBP_i} and \eqref{eq:IBP_e} are the cases $\lambda_1 = \lambda$ and $\lambda_2 \equiv 0$ and $\lambda_2 \equiv 1$, respectively.
\end{proof}

\begin{remark}
The formulas \eqref{eq:IBP_i} and \eqref{eq:IBP_e} are generalizations of the integration by parts formulas in \cite[Theorem 6.3]{CDM}.
\end{remark}

\section{The one-dimensional case} \label{sec:1D}

The case $N=1$ presents essential differences in some instances, and therefore we consider it separately. First of all, for all open sets $\Omega \subset \R$ we have the following identifications:
\begin{equation*}
\DM(\Omega) = \DM^1(\Omega) = BV(\Omega),
\end{equation*}
and analogously for the local versions. Indeed, if $N = 1$, then the operator $\div$ reduces only to the first (distributional) derivative, which we denote by $D$; so that the equivalence between $\DM^1(\Omega)$ and $BV(\Omega)$ is trivial. In addition, $\A \in \DM(\Omega)$ if and only if $\A, D \A  \in \mathcal{M}(\Omega)$: this implies that $\A \ll \Leb{1}$, with $\A = \overline{\A} \Leb{1}$ for some $\overline{\A} \in BV(\Omega)$ (see for instance \cite[Exercise 3.2]{AFP}). Clearly, the opposite is also true; that is, for all ${\bf B} \in BV(\Omega)$ the measure ${\bf B} \Leb{1}$ belongs to $\DM(\Omega)$. In other words, there is a bijection between $\DM(\Omega)$ and $BV(\Omega)$. 
Analogously, we see that
\begin{equation*}
\DM^{p}(\Omega) = \{ \A \in L^p(\Omega) : D \A \in \mathcal{M}(\Omega) \} \subset L^p(\Omega) \cap BV_{\rm loc}(\Omega)
\end{equation*}
for all $p \in (1, + \infty]$. Clearly, if $\Leb{1}(\Omega) < + \infty$, then, for all $p \in [1, + \infty]$, we have
\begin{equation*}
\DM^\infty(\Omega) \subseteq \DM^p(\Omega) \subseteq \DM^1(\Omega) = BV(\Omega) \subseteq \DM^\infty(\Omega),
\end{equation*}
due to the embedding $BV(\Omega) \subset L^\infty(\Omega)$, so that we obtain
$$\DM^p(\Omega) = BV(\Omega) \ \text{ for all } p \in [1, + \infty].$$ 
It is also interesting to point out that in general we have $\DM^p(\Omega) \subset L^\infty(\Omega)$ for all $p \in [1, + \infty]$, given that, if $\A \in L^1_{\rm loc}(\Omega)$ is such that $D \A \in \mathcal{M}(\Omega)$, then $\A \in L^\infty(\Omega)$ (see \cite[Sect. 3.2]{AFP}). 
Because of these facts, in this section we shall choose $\A \in BV(\Omega)$.

We gather in the following proposition the main basic properties of $BV^{\A, \lambda}$ in the case $N = 1$.

\begin{proposition} \label{prop:BV_A_emb_1}
Let $\A \in BV(\Omega)$ and $\lambda : \Omega \to [0, 1]$ be a Borel function. Then we have $u \in BV^{\A, \lambda}(\Omega)$ if and only if $u \in X^{\A, \lambda}(\Omega)$ and $u \A \in BV(\Omega)$, in which case we have
\begin{equation*}
D(u \A) =  u^\lambda D \A + (\A, Du)_\lambda \ \text{ on } \Omega
\end{equation*}
and
\begin{equation} \label{eq:BV_BV_A_emb}
\|u \A\|_{BV(\Omega)} \le \|u\|_{L^1(\Omega, |\A| \Leb{1})} + \|u^\lambda\|_{L^1(\Omega, |D \A|)} + |(\A, Du)_\lambda|(\Omega).
\end{equation} 
In particular, if $u \in BV^{\A, \lambda}(\Omega)$, then $u \A \in L^{\infty}(\Omega)$.
\end{proposition}

\begin{proof}
Given that $N = 1$, we know that $\div(u \A) = D(u \A)$, so that the first part of the statement follows from point (2) of Proposition \ref{prop:main_inclusions}. The second part is a consequence of the embedding $BV(\Omega) \subset L^\infty(\Omega)$.
\end{proof}

As in the higher dimensional case (without additional assumptions on $\A$, given that it is always locally essentially bounded), we point out that the inclusion of $BV(\Omega)$ in $BV^{\A, \lambda}(\Omega)$ is strict.

\begin{example}
Let $\Omega = (-1, 1)$, $\A(x) = \chi_{\left ( \frac{1}{2}, 1 \right )}(x)$ and $u(x) = \log{(|x|)}$. Then $u \in BV^{\A, \lambda}(\Omega) \setminus BV_{\rm loc}(\Omega)$ for every Borel function $\lambda : \Omega \to [0, 1]$. Indeed, $D\A = \delta_{\frac{1}{2}}$, and so $u \in X^{\A, \lambda}(\Omega)$ for every Borel function $\lambda : \Omega \to [0, 1]$, given that $u \in C(\Omega \setminus \{0\})$, and so $u^\lambda(x) = u(x)$ for all $x \in \Omega \setminus \{0\}$. In addition, we see that $u \A \in BV(\Omega)$, with
\begin{equation*}
D(u\A) = \log{\left ( \frac{1}{2} \right )} \delta_{\frac{1}{2}} + \frac{1}{x} \Leb{1} \res \left ( \frac{1}{2}, 1 \right ).
\end{equation*}
All in all, by Proposition \ref{prop:BV_A_emb_1}, we obtain $u \in BV^{\A, \lambda}(\Omega)$ with
\begin{equation*}
(\A, Du)_\lambda = - u\left ( \frac{1}{2} \right ) \delta_{\frac{1}{2}} +  \log{\left ( \frac{1}{2} \right )} \delta_{\frac{1}{2}} + \frac{1}{x} \Leb{1} \res \left ( \frac{1}{2}, 1 \right ) = \frac{1}{x} \Leb{1} \res \left ( \frac{1}{2}, 1 \right ).
\end{equation*}
\end{example}

The peculiarity of the one dimensional case lies in the fact that we achieve compactness in a fashion similar to the classical $BV$ space, under suitable assumptions.

\begin{proposition}\label{prop:1_dim_compact}
Let $\A \in BV(\Omega)$ and $\overline{\lambda} : \Omega \to [0, 1]$ be a Borel function. Let $(u_k) \subset BV^{\A, \overline{\lambda}}(\Omega)$ be a sequence of functions such that
\begin{equation} \label{eq:bounded_unif_uno}
\sup_{k \in \N} \|u_k\|_{L^1(\Omega, |\A| \Leb{1})} + \|u^{\overline{\lambda}}_k\|_{L^1(\Omega, |D \A|)} + |(\A, D u_k)_{\overline{\lambda}}|(\Omega) < + \infty.
\end{equation}
Then there exist $u \in L^1(\Omega, |\A| \Leb{1})$ such that $u \A \in BV(\Omega)$ and a subsequence $(u_{k_j})_{j \in \N}$ such that $u_{k_j} \to u \in L^1(\Omega, |\A| \Leb{1})$. In addition, assume that at least one of the following conditions is satisfied:
\begin{enumerate}[i)]
\item there exists $c > 0$ such that $|\A(x)| > c$ for $\Leb{1}$-a.e. $x \in \Omega$,
\item $\A \in W^{1, 1}(\Omega)$,
\item the sequence $(u_k)_{k \in \N}$ is uniformly bounded in $L^\infty(\Omega)$.
\end{enumerate} 
Then $u \in BV^{\A, \lambda}(\Omega)$ for every Borel function $\lambda : \Omega \to [0, 1]$.
Finally, if $\A \in W^{1, 1}(\Omega)$ and the sequence $(u_k)_{k \in \N}$ is uniformly bounded in $L^\infty(\Omega)$, then, possibly up to a further subsequence, $u_{k_j} \to u$ in $L^1(\Omega, |D \A|)$; so that $(u_{k_j})_{j \in \N}$ $\A$-converges to $u$.
\end{proposition}

\begin{proof}
By combining to \eqref{eq:BV_BV_A_emb} and \eqref{eq:bounded_unif_uno}, we see that the sequence $(u_k \A)_{k \in \N}$ is uniformly bounded in $BV(\Omega)$. Thanks to the compactness theorem in $BV$ \cite[Theorem 3.3]{AFP}, we deduce that there exists $w_{\A} \in BV(\Omega)$ and a subsequence $(u_{k_j} \A)_{j \in \N}$ such that $u_{k_j} \A \to w_{\A}$ in $L^1(\Omega)$. With a little abuse of notation, we still denote by $\A$ and $w_{\A}$ the Borel representatives of these two functions. Then, we set
\begin{equation*}
u(x) = \begin{cases} \frac{w_{\A}(x)}{\A(x)} & \text{ if } \A(x) \neq 0 \\ 0 & \text{ otherwise} \end{cases}.
\end{equation*}
Therefore, we clearly get $u \in L^1(\Omega, |\A| \Leb{1})$, $u \A \in BV(\Omega)$ and $u_{k_j} \to u$ in $L^1(\Omega, |\A| \Leb{1})$. We assume now that there exists $c > 0$ such that $|\A(x)| > c$ for $\Leb{1}$-a.e. $x \in \Omega$. Then, $u \in L^\infty(\Omega)$ by definition, given that $w_{\A} \in BV(\Omega) \subset L^\infty(\Omega)$. Hence, we immediately obtain $u^\lambda \in L^1(\Omega, |D \A|)$ for every Borel function $\lambda : \Omega \to [0, 1]$, and so, by Proposition \ref{prop:BV_A_emb_1}, we get $u \in BV^{\A, \lambda}(\Omega)$.
If instead $\A \in W^{1, 1}(\Omega)$, then by Proposition \ref{rem:div_A_abs_cont_Leb} we have $BV^{\A, \lambda}(\Omega) = BV^{\A}(\Omega)$ for every Borel function $\lambda : \Omega \to [0, 1]$. Then we notice that, up to extracting a further subsequence, we have $u_{k_j}(x) \to u(x)$ for $|\A| \Leb{1}$-a.e. $x \in \Omega$, and therefore $u_{k_j}(x) \to u(x)$ for $\Leb{1}$-a.e. $x \in {\rm supp}(|\A|)$. Given that ${\rm supp}(|D \A|) \subset {\rm supp}(|\A|)$, we conclude that $u_{k_j}(x) \to u(x)$ for $\Leb{1}$-a.e. $x \in {\rm supp}(|D\A|)$. Thus, since $D \A \ll \Leb{1}$, this implies 
\begin{equation} \label{eq:subseq_pointwise_conv_DA}
u_{k_j}(x) \to u(x) \ \text{ for } \ |D\A|\text{-a.e. } x \in \Omega,
\end{equation}
and so by Fatou's Lemma we obtain
\begin{equation*}
\int_{\Omega} |u| \, d |D\A| = \int_{\Omega} \liminf_{j \to + \infty} |u_{k_j}| \, d |D \A| \le \liminf_{j \to + \infty} \int_{\Omega} |u_{k_j}| \, d |D \A| < + \infty.
\end{equation*}
All in all, we get $u \in L^1(\Omega, |D \A|)$, and so $u \in BV^{\A}(\Omega)$, by Proposition \ref{prop:BV_A_emb_1}.
As for condition (iii), we claim that it entails $u \in L^\infty(\Omega)$. To see this, we let $C > 0$ be such that $\|u_k\|_{L^\infty(\Omega)} \le C$ for all $k \in \N$. As noted above, up to a further subsequence, we have $u_{k_j}(x) \to u(x)$ for $|\A| \Leb{1}$-a.e. $x \in \Omega$, and so
\begin{equation*}
|u(x)| = \lim_{j \to + \infty} |u_{k_j}(x)| \le \liminf_{j \to + \infty} \|u_{k_j}\|_{L^\infty(\Omega)} \le C \ \text{ for } |\A|\Leb{1}\text{-a.e. } x \in \Omega.
\end{equation*}
This implies that $u \in L^{\infty}(\Omega, |\A| \Leb{1})$ with $\|u\|_{L^{\infty}(\Omega, |\A| \Leb{1})} \le C$, which means
\begin{equation*} 
\int_{\{ |u| > t \}} |\A| \, dx = 0 \ \text{ for all } t > C.
\end{equation*}
By definition, $u(x) = 0$ for $\Leb{1}$-a.e. $x \in \Omega$ such that $\A(x) = 0$. Therefore, by Chebyshev's inequality, for all $t > C$ we obtain
\begin{align*}
\Leb{1}(\{ |u| > t \}) & = \Leb{1}(\{ |u| > t \} \cap \{ |\A| > 0 \}) \\
& = \Leb{1}\left (\{ |u| > t \} \cap \left ( \{ |\A| \ge 1 \} \cup \bigcup_{k = 1}^{+\infty} \left \{ \frac{1}{k} > |\A| \ge \frac{1}{k+1} \right \} \right ) \right ) \\
& \le \Leb{1}\left (\{ |u| > t \} \cap \{ |\A| \ge 1 \} \right ) + \sum_{k = 1}^{+\infty} \Leb{1}\left (\{ |u| > t \} \cap \left \{ \frac{1}{k} > |\A| \ge \frac{1}{k+1} \right \} \right ) \\
& \le \int_{\{ |u| > t \}} |\A| \, dx + \sum_{k = 1}^{+\infty} \int_{\{ |u| > t \}} (k+1) |\A| \, dx = 0.
\end{align*}
All in all, this entails that $\|u\|_{L^\infty(\Omega)} \le C$, and so, as above, we have $u^\lambda \in L^1(\Omega, |D \A|)$ for every Borel function $\lambda : \Omega \to [0, 1]$, which in turn, by Proposition \ref{prop:BV_A_emb_1}, implies $u \in BV^{\A, \lambda}(\Omega)$. Finally, let $\A \in W^{1, 1}(\Omega)$ and $C > 0$ be such that $\|u_k\|_{L^\infty(\Omega)} \le C$. Arguing as above, up to extracting possibly a further subsequence, \eqref{eq:subseq_pointwise_conv_DA} holds, and therefore, given that
\begin{equation*}
|u_{k_j} - u| \le 2 C \in L^1(\Omega, |D\A|),
\end{equation*}
we exploit Lebesgue's Dominated Convergence Theorem to conclude that $u_{k_j} \to u$ in $L^1(\Omega, |D \A|)$, and so $(u_{k_j})_{j \in \N}$ $\A$-converges to $u$.
\end{proof}

\begin{remark}
We point out that in the last part of Proposition \ref{prop:1_dim_compact} we cannot remove the assumption that $A \in W^{1,1}(\Omega)$ and still achieve the convergence of the $\lambda$-representatives in $L^1(\Omega, |D \A|)$. Indeed, let $\Omega = (-1,1)$, $\A(x) = \chi_{(0,1)}(x)$, 
$$u_{k}(x) = \begin{cases} a \arctan{(k x)} & \text{ if } x \ge 0 \\
b \arctan{(k x)} & \text{ if } x < 0 \end{cases}, $$
for some $a, b \ge 0$, and $\lambda : \Omega \to [0, 1]$ be any Borel function. Since
$$D\A = \delta_0 \ \text{ and } \ u_k^\lambda(0) = u_k(0) = 0,$$ 
we see that 
$$D(u_k \A) = \frac{a k}{1+k^2 x^2} \Leb{1} \res (0, 1),$$
and so $u_k \in BV^{\A, \lambda}(\Omega)$, by Proposition \ref{prop:1_dim_compact}, with
\begin{equation*}
(\A, Du_k)_\lambda = - u_k^\lambda(0) \delta_0 + D(u_k \A) = \frac{a k}{1+k^2 x^2} \Leb{1} \res (0, 1).
\end{equation*}
It is clear that the sequence $(u_k)_{k \in \N}$ satisfies \eqref{eq:bounded_unif_uno}, is uniformly bounded in $L^\infty(\Omega)$ and converges to 
$$u(x) = \begin{cases} a \frac{\pi}{2} & \text{ if } x > 0, \\
0 & \text{ if } x = 0, \\
- b \frac{\pi}{2} & \text{ if } x < 0 \end{cases} $$
pointwise and in $L^1(\Omega, |\A| \Leb{1})$. However, 
\begin{equation*}
u^\lambda(0) = (1 - \lambda(0)) \left ( - b \frac{\pi}{2} \right ) + \lambda(0) a \frac{\pi}{2} = ( (a + b) \lambda(0) - b) \frac{\pi}{2},
\end{equation*}
so that $u_k^\lambda$ does not converge to $u^\lambda$ in $L^1(\Omega, |D \A|)$, as long as $(a + b) \lambda(0) \neq b$. Furthermore, under such condition we do not even have any lower semicontinuity with respect to the $L^1(\Omega, |D\A|)$ norm of the $\lambda$-representatives.
\end{remark}

\begin{corollary}\label{cor:1_dim_compact_1}
Let $\A \in BV(\Omega)$ and $\overline{\lambda} : \Omega \to [0, 1]$ be a Borel function. Let $(E_k)$ be a family of Borel sets such that $$\sup_{k \in \N} P_{\A, \overline{\lambda}}(E_k, \Omega) < + \infty.$$ 
Then there exist a Borel set $E$ such that $P_{\A, \lambda}(E, \Omega) < + \infty$ for every Borel function $\lambda : \Omega \to [0, 1]$ and a subsequence $(E_{k_j})_{j \in \N}$ such that $\chi_{E_{k_j}} \to \chi_E$ in $L^1(\Omega, |\A| \Leb{1})$. If $\A \in W^{1, 1}(\Omega)$, we also have $\chi_{E_{k_j}} \to \chi_E$ in $L^1(\Omega, |D \A|)$, possibly up to a further subsequence, which entails the $\A$-convergence.
\end{corollary}

\begin{proof}
Arguing as in the proof of Proposition \ref{prop:1_dim_compact}, we see that there exist $w_{\A} \in BV(\Omega)$ and a subsequence $(\chi_{E_{k_j}} \A)_{j \in \N}$ such that $\chi_{E_{k_j}} \A \to w_{\A}$ in $L^1(\Omega)$. Therefore, $\chi_{E_{k_j}} \to \frac{w_{\A}}{\A}$ in $L^1(\Omega, |\A| \Leb{1})$, and this readily implies that $\frac{w_{\A}(x)}{\A(x)} \in \{0, 1\}$ for $|\A| \Leb{1}$-a.e. $x \in \Omega$. Thus, there exists a Borel set $E$ such that $\frac{w_{\A}(x)}{\A(x)} = \chi_E(x)$ for $|\A| \Leb{1}$-a.e. $x \in \Omega$. It is plain to see that $\chi_E \A \in BV(\Omega)$ and $\chi_E \in L^\infty(\Omega) \subset X^{\A, \lambda}(\Omega)$ for every Borel function $\lambda : \Omega \to [0, 1]$, so that Proposition \ref{prop:BV_A_emb_1} implies $\chi_E \in BV^{\A, \lambda}(\Omega)$ for every Borel function $\lambda : \Omega \to [0, 1]$. Let now $\A \in W^{1, 1}(\Omega)$. Given that $\|\chi_{E_k}\|_{L^\infty(\Omega)} \le 1$, we just need to employ Proposition \ref{prop:1_dim_compact} to conclude that, up to extracting a further subsequence, $\chi_{E_{k_j}} \to \chi_E$ in $L^1(\Omega, |D \A|)$, and so $(E_{k_j})_{j \in \N}$ $\A$-converges to $E$.
\end{proof}


\begin{remark}
Proposition \ref{prop:1_dim_compact} implies that, if $\A \in \DM^1(\Omega)$ with $D \A \ll \Leb{1}$ and $\Leb{1}(\Omega \setminus {\rm supp}(|\A|)) = 0$, we have a weak compactness result in $BV^{\A}(\Omega)$ (which is actually a Banach space, due to Proposition \ref{prop:property}), as well as Corollary \ref{cor:1_dim_compact_1} implies a weak compactness in the class of sets with finite $\A$-perimeter. However, if $N \ge 2$, analogous results fail to hold true due to new degrees of freedom, even if $\div \A  = 0$ (see Examples \ref{exmp:not_compactness} and \ref{ex:per_nocompact}).
\end{remark}

\appendix \section{An alternative definition of $u^\lambda$} \label{rem:vecchio_approccio}

We point out that, if $u \in L^1_{\rm loc}(\Omega)$, we could define the general $\lambda$-pairings by exploiting the representatives of $u$ given by the limit of averages on balls and half-balls; that is, \eqref{eq:approxlim1} and \eqref{f:disc}, respectively.
This is indeed the approach followed in the classical setting of \cite{Anz} for $\lambda \equiv \frac{1}{2}$, and for a general Borel function $\lambda : \Omega \to [0,1]$ in \cite{CDM}: therein, the field $\A$ is $\DM^\infty$ and the scalar function $u$ is $BV$, and therefore it is well known that the function 
\begin{equation} \label{eq:lambda_averages_def}
u^\lambda(x) = \begin{cases} \tilde{u}(x)\,, & \text{ if } x \in \Omega \setminus S_u\,, \\
\\
(1-\lambda(x))\min\{\uint(x),\uext(x)\} +\lambda(x)\max\{\uint(x),\uext(x)\}\,, & \text{ if } x \in J_u,
\end{cases}
\end{equation}
is well defined $|\div \A|$-a.e. (see the argument in the proof of point (3) in Proposition \ref{prop:main_inclusions}). Indeed, we notice that, if $u \in BV(\Omega)$, then $\Haus{N-1}(Z_u) \le \Haus{N-1}(S_u \setminus J_u) = 0$, so that our definition \eqref{f:pr} coincides with \eqref{eq:lambda_averages_def} up to an $\Haus{N-1}$-negligible set, due to \eqref{eq:u_+_-_tilde} and \eqref{eq:u_+_-_i_e}. However, as soon as $u \notin BV_{\rm loc}(\Omega)$ or $\A \in \DM_{\rm loc}(\Omega) \setminus L^{\infty}_{\rm loc}(\Omega; \R^N)$, in general the function $u^\lambda$ given by \eqref{eq:lambda_averages_def} is not well defined $|\A|$-a.e. or $|\div \A|$-a.e., since we have no control, a priori, on the set $S_u \setminus J_u$, except that $\Leb{N}(S_u \setminus J_u) = 0$, obviously. An alternative possible approach would be to consider functions $u \in L^1_{\rm loc}(\Omega)$ such that 
\begin{equation} \label{eq:S_u_J_u_A_negligible}
|\A^s|(S_u \setminus J_u) = |\Div^s \A|(S_u\setminus J_u)=0
\end{equation} 
and $u^\lambda \in L^1_{\rm loc}(\Omega, |\A|) \cap L^1_{\rm loc}(\Omega, |\div \A|)$. Under these assumptions, the $\lambda$-pairing can be defined as done in \eqref{eq:def_lambda_pairing}, and therefore we can provide alternative definitions of $X^{\A, \lambda}(\Omega)$ and $BV^{\A, \lambda}(\Omega)$ as the classes
\begin{equation*}
\widetilde{X}^{\A, \lambda}(\Omega) := \left\{ u \in L^1(\Omega): (|\A^s| + |\Div^s \A|)(S_u\setminus J_u)=0, u^\lambda  \in L^1(\Omega, |\A|) \cap L^1(\Omega, |\div \A|) \right \},
\end{equation*}
\begin{equation*}
\widetilde{BV}^{\A, \lambda}(\Omega) := \left\{ u \in \widetilde{X}^{\A, \lambda}(\Omega): (\A,Du)_{\lambda} \in \radon \right\}
\end{equation*}
and analogously for the corresponding local versions. Due to the additional condition on $S_u \setminus J_u$, there are no obvious inclusions between $X^{\A, \lambda}(\Omega)$ and $BV^{\A, \lambda}(\Omega)$ and these alternative versions, unless $(|\A| + |\div \A|) \ll \Leb{N}$, in which case $\widetilde{X}^{\A, \lambda}(\Omega) \subseteq X^{\A, \lambda}(\Omega)$ and $\widetilde{BV}^{\A, \lambda}(\Omega) \subseteq BV^{\A, \lambda}(\Omega)$.
Furthermore, in the particular case $\lambda \equiv \frac{1}{2}$, we could completely identify $u^{\frac{1}{2}}$ with the precise representative $u^*$, in such way removing the problematic related to the set $S_u \setminus J_u$. Therefore, we could drop \eqref{eq:S_u_J_u_A_negligible}, and replace it with the requirement that the limit \eqref{def:precise_repr} exists for $(|\A| + |\div \A|)$-a.e. $x \in \Omega$. This is indeed coherent with the fact that the pointwise limit of any mollification $(u \ast \rho_{\eps})(x)$ is $u^*(x)$, whenever the precise representative is well defined (see Theorem \ref{t:lscnuovo}). This approach takes also in account the possibility of having a precise representative which does not necessarily coincide with $\frac{u^+(x) + u^{-}(x)}{2}$ for $x \in S_u \setminus J_u$: consider for instance $u = \chi_{(0, 1)^N}$ and $x = 0$.

Nevertheless, such an analysis would have been out the scope of the present paper, so we preferred to exploit the fact that the approximate liminf and limsup $u^{\pm}$ are always well defined, in order to circumvent all these difficulties involving the well-posedness of the representatives and to avoid any assumption on the set $S_u \setminus J_u$. Furthermore, in the classical setting there are no substantial differences. We defer the exploration of this alternative approach to a future research.


\begin{bibdiv}
\begin{biblist}

\bib{AmbCriMan}{article}{
      author={Ambrosio, {L.}},
      author={Crippa, {G.}},
      author={Maniglia, {S.}},
       title={Traces and fine properties of a {$BD$} class of vector fields and
  applications},
        date={2005},
        ISSN={0240-2963},
     journal={Ann. Fac. Sci. Toulouse Math. (6)},
      volume={14},
      number={4},
       pages={527\ndash 561},
         url={http://afst.cedram.org/item?id=AFST_2005_6_14_4_527_0},
}

\bib{AFP}{book}{
      author={Ambrosio, {L.}},
      author={Fusco, {N.}},
      author={Pallara, {D.}},
       title={Functions of bounded variation and free discontinuity problems},
      series={Oxford Mathematical Monographs},
   publisher={The Clarendon Press Oxford University Press},
     address={New York},
        date={2000},
        ISBN={0-19-850245-1},
}

\bib{MR1814993}{article}{
   author={Andreu, F.},
   author={Ballester, C.},
   author={Caselles, V.},
   author={Maz\'{o}n, J. M.},
   title={The Dirichlet problem for the total variation flow},
   journal={J. Funct. Anal.},
   volume={180},
   date={2001},
   number={2},
   pages={347--403},
   doi={10.1006/jfan.2000.3698},
}
	

\bib{AVCM}{book}{
      author={Andreu-Vaillo, {F.}},
      author={Caselles, {V.}},
      author={Maz\'on, {J. M.}},
       title={Parabolic quasilinear equations minimizing linear growth
  functionals},
      series={Progress in Mathematics},
   publisher={Birkh\"auser Verlag, Basel},
        date={2004},
      volume={223},
        ISBN={3-7643-6619-2},
         url={http://dx.doi.org/10.1007/978-3-0348-7928-6},
}

\bib{Anz}{article}{
      author={Anzellotti, {G.}},
       title={Pairings between measures and bounded functions and compensated
  compactness},
        date={1983},
        ISSN={0003-4622},
     journal={Ann. Mat. Pura Appl. (4)},
      volume={135},
       pages={293\ndash 318 (1984)},
         url={http://dx.doi.org/10.1007/BF01781073},
}


\bib{Anz2}{misc}{
      author={Anzellotti, {G.}},
       title={Traces of bounded vector--fields and the divergence theorem},
        date={1983},
        note={Unpublished preprint},
}

\bib{BuffaComiMira}{article}{
   author={Buffa, V.},
   author={Comi, G. E.},
   author={Miranda, M. Jr.},
   title={On BV functions and essentially bounded divergence-measure fields
   in metric spaces},
   journal={Rev. Mat. Iberoam.},
   volume={38},
   date={2022},
   number={3},
   pages={883--946},
   doi={10.4171/rmi/1291},
}
	

\bib{Cas}{article}{
      author={Caselles, V.},
       title={On the entropy conditions for some flux limited diffusion
  equations},
        date={2011},
        ISSN={0022-0396},
     journal={J. Differential Equations},
      volume={250},
      number={8},
       pages={3311\ndash 3348},
         url={http://dx.doi.org/10.1016/j.jde.2011.01.027},
}

\bib{ChCoTo}{article}{
      author={Chen, {G.-Q.}},
      author={Comi, {G. E.}},
      author={Torres, {M.}},
       title={Cauchy fluxes and {G}auss-{G}reen formulas for divergence-measure
  fields over general open sets},
        date={2019},
        ISSN={0003-9527},
     journal={Arch. Ration. Mech. Anal.},
      volume={233},
      number={1},
       pages={87\ndash 166},
         url={https://doi.org/10.1007/s00205-018-01355-4},
}



\bib{ChenFrid}{article}{
      author={Chen, {G.-Q.}},
      author={Frid, {H.}},
       title={Divergence-measure fields and hyperbolic conservation laws},
        date={1999},
        ISSN={0003-9527},
     journal={Arch. Ration. Mech. Anal.},
      volume={147},
      number={2},
       pages={89\ndash 118},
         url={http://dx.doi.org/10.1007/s002050050146},
}

\bib{ChFr1}{article}{
      author={Chen, {G.-Q.}},
      author={Frid, {H.}},
       title={Extended divergence-measure fields and the {E}uler equations for
  gas dynamics},
        date={2003},
        ISSN={0010-3616},
     journal={Comm. Math. Phys.},
      volume={236},
      number={2},
       pages={251\ndash 280},
         url={http://dx.doi.org/10.1007/s00220-003-0823-7},
}


\bib{ChIrTo}{article}{
      author={Chen, {G.-Q.}},
      author={Irving, {C.}},
      author={Torres, {M.}},
       title={Extended divergence-measure fields, the Gauss-Green formula, and Cauchy fluxes},
        date={2025},
        ISSN={},
     journal={Archive for Rational Mechanics and Analysis},
      volume={249},
      number={6},
       pages={79},
         url={},
}


\bib{ChDCMe1}
{article}{
author={Chiad\`o Piat, {V.}},
author = {De Cicco, {V.}},
author = {Melchor Hernandez, {A.}},
title={Relaxation for a degenerate functional with linear growth in the onedimensional case},
journal={J. Convex Anal.},
volume={33},
number={2},
year={2026}
}

\bib{ChDCMe2}
{article}{
author={Chiad\`o Piat, {V.}},
author = {De Cicco, {V.}},
author = {Melchor Hernandez, {A.}},
title={Relaxation for degenerate nonlinear functionals in the onedimensional case},
journal={Nonlinear Differ. Equ. Appl. },
volume={32},
number={4},
year={2025}
}


\bib{Comi}{article}{
   author={Comi, G. E.},
   title={Refined Gauss--Green formulas and evolution problems for Radon measures},
   journal={Scuola Normale Superiore, Pisa},
   date={2020},
 note={Ph.D.\ Thesis. Available at \href{https://cvgmt.sns.it/paper/4579/}{cvgmt.sns.it/paper/4579/}}
}

\bib{CCDM}{article}{
  title={Representation formulas for pairings between divergence-measure fields and {B}{V} functions},
  author={Comi, G. E.}
  author={Crasta, G.}, 
  author={De Cicco, V.}, 
  author={Malusa, A.},
  journal={J. Funct. Anal.},
  volume={286},
  number={1},
  pages={110192},
  year={2024},
}

\bib{ComiLeo}{article}{
author = {Comi, {G. E.}},
author={Leonardi, {G. P.}},
title={Measures in the dual of $BV$: perimeter bounds and relations with divergence-measure fields},
journal={Nonlinear Analysis},
volume={251},
pages={113686}
year={2025}
}

\bib{ComiMagna}{article}{
   author={Comi, G. E.},
   author={Magnani, V.},
   title={The Gauss-Green theorem in stratified groups},
   journal={Adv. Math.},
   volume={360},
   date={2020},
   pages={106916, 85},
   doi={10.1016/j.aim.2019.106916},
}

\bib{ComiPayne}{article}{
   author={Comi, G. E.},
   author={Payne, K. R.},
   title={On locally essentially bounded divergence measure fields and sets
   of locally finite perimeter},
   journal={Adv. Calc. Var.},
   volume={13},
   date={2020},
   number={2},
   pages={179--217},
   doi={10.1515/acv-2017-0001},
}


\bib{CD3}{article}{
      author={Crasta, {G.}},
      author={De~Cicco, {V.}},
       title={Anzellotti's pairing theory and the {G}auss--{G}reen theorem},
        date={2019},
        ISSN={0001-8708},
     journal={Adv. Math.},
      volume={343},
       pages={935\ndash 970},
        url={https://doi.org/10.1016/j.aim.2018.12.007},
}

\bib{crasta2019extension}{article}{
  title={An extension of the pairing theory between divergence-measure fields and {B}{V} functions},
  author={Crasta, G.}, 
  author={De Cicco, V.},
  journal={J. Funct. Anal.},
  volume={276},
  number={8},
  pages={2605--2635},
  year={2019},
  publisher={Elsevier}
}


\bib{CD5}{article}{
author={Crasta, {G.}},
      author={De~Cicco, {V.}},
       title={On the variational nature of the Anzellotti pairing},
journal={Adv. Calc. Var.},
volume={18},
  number={3},
  pages={755--771},
  year={2025},
  publisher={De Gruyter}
url={https://doi.org/10.1515/acv-2024-0067}
}

\bib{CDM}{article}{
      author={Crasta, {G.}},
      author={De~Cicco, {V.}},
      author={Malusa, {A.}},
       title={Pairings between bounded divergence-measure vector fields and {BV} functions},
        date={2022},
     journal={Adv. Calc. Var.},
       pages={787-810},
             volume={15},
            number={4},
         url={https://doi.org/10.1515/acv-2020-0058},
}


\bib{DCFV2}{article}{
      author={De~Cicco, V.},
      author={Fusco, N.},
      author={Verde, A.},
       title={A chain rule formula in {$BV$} and application to lower
  semicontinuity},
        date={2007},
        ISSN={0944-2669},
     journal={Calc. Var. Partial Differential Equations},
      volume={28},
      number={4},
       pages={427\ndash 447},
         url={http://dx.doi.org/10.1007/s00526-006-0048-7},
}

\bib{MR3939259}{article}{
   author={De Cicco, V.},
   author={Giachetti, D.},
   author={Segura de Le\'{o}n, S.},
   title={Elliptic problems involving the 1-Laplacian and a singular lower
   order term},
   journal={J. Lond. Math. Soc. (2)},
   volume={99},
   date={2019},
   number={2},
   pages={349--376},
   doi={10.1112/jlms.12172},
}


\bib{DGMM}{article}{
      author={Degiovanni, {M.}},
      author={Marzocchi, {A.}},
      author={Musesti, {A.}},
       title={Cauchy fluxes associated with tensor fields having divergence
  measure},
        date={1999},
        ISSN={0003-9527},
     journal={Arch. Ration. Mech. Anal.},
      volume={147},
      number={3},
       pages={197\ndash 223},
         url={http://dx.doi.org/10.1007/s002050050149},
}

\bib{DELNIN}{article}{
      author={Del Nin, {G.}},
             title={Rectifiability of the jump set of locally integrable functions},
        date={2021},
     journal={Ann. Sc. Norm. Super. Pisa Cl. Sci.},
      volume={XXII},
      number={5},
       pages={1233\ndash 1240},
         url={http://dx.doi.org/10.2422/2036-2145.202002_006},
}

\bib{evans2015measure}{book}{
   author={Evans, L. C.},
   author={Gariepy, R. F.},
   title={Measure theory and fine properties of functions},
   series={Textbooks in Mathematics},
   edition={Revised edition},
   publisher={CRC Press, Boca Raton, FL},
   date={2015},
   pages={xiv+299},
}

\bib{FED}{book}{
      author={Federer, H.},
       title={Geometric measure theory},
      series={Die Grundlehren der mathematischen Wissenschaften, Band 153},
   publisher={Springer-Verlag New York Inc., New York},
        date={1969},
}

\bib{Frid}{book}{
author={Frid, H.},
title = {Divergence-measure fields on domains with Lipschitz boundary},
series = {Hyperbolic conservation laws and related analysis with applications},
publisher={Springer Proc. Math. Stat., Springer, Heidelberg,},
volume={49},
pages = {207--225},
year={(2014)}
}

\bib{MR4381316}{article}{
   author={G\'{o}rny, Wojciech},
   title={Local and nonlocal 1-Laplacian in Carnot groups},
   journal={Ann. Fenn. Math.},
   volume={47},
   date={2022},
   number={1},
   pages={427--456},
   doi={10.54330/afm.114742},
}

\bib{gorny2023anzellotti}{article}{
  author={G\'{o}rny, Wojciech},
  author={Maz\'{o}n, Jos\'{e} M.},
  title={The Anzellotti-Gauss-Green Formula and Least Gradient Functions in Metric Measure Spaces},
  journal={Communications in Contemporary Mathematics},
  year={2023},
  publisher={World Scientific}
}


\bib{Ir}{article}{
      author={Irving, {C.}},
       title={On the normal trace space of extended divergence-measure fields},
        date={2025},
        ISSN={},
     journal={arXiv:2503.09536},
      volume={},
      number={},
       pages={},
         url={},
}


\bib{MR2348842}{article}{
   author={Kawohl, B.},
   author={Schuricht, F.},
   title={Dirichlet problems for the 1-Laplace operator, including the
   eigenvalue problem},
   journal={Commun. Contemp. Math.},
   volume={9},
   date={2007},
   number={4},
   pages={515--543},
   doi={10.1142/S0219199707002514},
}

\bib{LeoComi}{article}{
author={Leonardi, {G. P.}},
author = {Comi, {G. E.}},
title={The prescribed mean curvature measure equation in non-parametric form},
journal={Ann. Sc. Norm. Super. Pisa Cl. Sci.},
year={2025}
}


\bib{LeoSar}{article}{
      author={Leonardi, {G. P.}},
      author={Saracco, {G.}},
       title={The prescribed mean curvature equation in weakly regular
  domains},
        date={2018},
     journal={Nonlinear Differ. Equ. Appl. },
      volume={25},
      number={2},
       pages={Art. 9, 29},
         url={https://doi.org/10.1007/s00030-018-0500-3},
}

\bib{LeoSar2}{article}{
      author={Leonardi, {G.P.}},
      author={Saracco, {G.}},
       title={Rigidity and trace properties of divergence-measure vector fields},
          date={2022},
        journal={Adv. Calc. Var.},
   volume={15},
   number={1},
   pages={133--149},
        doi={10.1515/acv-2019-0094},
}

\bib{Maggi}{book}{
   author={Maggi, F.},
   title={Sets of finite perimeter and geometric variational problems},
   series={Cambridge Studies in Advanced Mathematics},
   volume={135},
   note={An introduction to geometric measure theory},
   publisher={Cambridge University Press, Cambridge},
   date={2012},
   pages={xx+454},
   doi={10.1017/CBO9781139108133},
}

\bib{MR1857126}{article}{
   author={Marzocchi, Alfredo},
   author={Musesti, Alessandro},
   title={Decomposition and integral representation of Cauchy interactions
   associated with measures},
   journal={Contin. Mech. Thermodyn.},
   volume={13},
   date={2001},
   number={3},
   pages={149--169},
   doi={10.1007/s001610100046},
}


\bib{MR2502520}{article}{
   author={Mercaldo, A.},
   author={Segura de Le\'{o}n, S.},
   author={Trombetti, C.},
   title={On the solutions to 1-Laplacian equation with $L^1$ data},
   journal={J. Funct. Anal.},
   volume={256},
   date={2009},
   number={8},
   pages={2387--2416},
   doi={10.1016/j.jfa.2008.12.025},
}

\bib{MR3676052}{article}{
   author={Phuc, Nguyen Cong},
   author={Torres, Monica},
   title={Characterizations of signed measures in the dual of $BV$ and
   related isometric isomorphisms},
   journal={Ann. Sc. Norm. Super. Pisa Cl. Sci. (5)},
   volume={17},
   date={2017},
   number={1},
   pages={385--417},
}

\bib{Rudin}{book}{
  title={Functional Analysis},
  author={Rudin, W.},
  series={International series in pure and applied mathematics},
  year={1974},
  publisher={Tata McGraw-Hill}
}

\bib{MR3501836}{article}{
   author={Scheven, C.},
   author={Schmidt, T.},
   title={BV supersolutions to equations of 1-Laplace and minimal surface
   type},
   journal={J. Differential Equations},
   volume={261},
   date={2016},
   number={3},
   pages={1904--1932},
   doi={10.1016/j.jde.2016.04.015},
}

\bib{MR3813962}{article}{
   author={Scheven, C.},
   author={Schmidt, T.},
   title={On the dual formulation of obstacle problems for the total
   variation and the area functional},
   journal={Ann. Inst. H. Poincar\'{e} C Anal. Non Lin\'{e}aire},
   volume={35},
   date={2018},
   number={5},
   pages={1175--1207},
   doi={10.1016/j.anihpc.2017.10.003},
}
	

\bib{Schu}{article}{
      author={Schuricht, {F.}},
       title={A new mathematical foundation for contact interactions in
  continuum physics},
        date={2007},
        ISSN={0003-9527},
     journal={Arch. Ration. Mech. Anal.},
      volume={184},
      number={3},
       pages={495\ndash 551},
         url={http://dx.doi.org/10.1007/s00205-006-0032-6},
}

\bib{Silh}{article}{
      author={\v{S}ilhav\'{y}, {M.}},
       title={Divergence measure fields and {C}auchy's stress theorem},
        date={2005},
        ISSN={0041-8994},
     journal={Rend. Sem. Mat. Univ. Padova},
      volume={113},
       pages={15\ndash 45},
}

\bib{MR2532602}{article}{
   author={\v{S}ilhav\'{y}, M.},
   title={The divergence theorem for divergence measure vectorfields on sets
   with fractal boundaries},
   journal={Math. Mech. Solids},
   volume={14},
   date={2009},
   number={5},
   pages={445--455},
   doi={10.1177/1081286507081960},
}


\bib{Silhavy19}{article}{
   author={\v{S}ilhav\'{y}, M.},
   title={The Gauss-Green theorem for bounded vector fields with divergence
   measure on sets of finite perimeter},
   journal={Indiana Univ. Math. J.},
   volume={72},
   date={2023},
   number={1},
   pages={29--42},
}

\end{biblist}
\end{bibdiv}

\def\cprime{$'$}

\end{document}